\documentclass[12pt,draftcls, onecolumn]{IEEEtran}

\usepackage{cite}

\ifCLASSINFOpdf
  \usepackage[pdftex]{graphicx}
  \graphicspath{{./}}
  \DeclareGraphicsExtensions{.jpeg,.png}
\else
\fi
\usepackage{amsmath}
\usepackage{amssymb}
\usepackage{xcolor}
\usepackage{amsfonts}
\usepackage{amssymb}
\usepackage{subcaption}

\newtheorem{theorem}{Theorem}
\newtheorem{definition}{Definition}
\newtheorem{remark}{Remark}
\newtheorem{corollary}{Corollary}

\begin{document}
%
\title{The Eigenvectors of  Single-spiked Complex Wishart Matrices: Finite and Asymptotic Analyses}
%
%
%

\author{
  Prathapasinghe Dharmawansa,~\IEEEmembership{Member,~IEEE,} Pasan~Dissanayake,
  and~Yang~Chen
  \thanks{P. Dharmawansa and P. Dissanayake  are with the Department
		of Electronic and Telecommunication Engineering, University of Moratuwa, Moratuwa 10400, Sri Lanka (e-mail: prathapa@uom.lk, pasandissanayake@gmail.com).}%
  \thanks{Y. Chen is with the Department of Mathematics, Faculty of Science and Technology, University of Macau, Macau, P. R. China (e-mail:
 yayangchen@umac.mo).}
}
\maketitle
\vspace{-4mm}
\begin{abstract}
Let $\mathbf{W}\in\mathbb{C}^{n\times n}$ be a {\it single-spiked} Wishart matrix in the class $\mathbf{W}\sim \mathcal{CW}_n(m,\mathbf{I}_n+ \theta \mathbf{v}\mathbf{v}^\dagger) $ with $m\geq n$,  where $\mathbf{I}_n$ is the $n\times n$ identity matrix, $\mathbf{v}\in\mathbb{C}^{n\times 1}$ is an arbitrary vector with unit Euclidean norm, $\theta\geq 0$ is a non-random parameter, and $(\cdot)^\dagger$ represents the conjugate-transpose operator. Let $\mathbf{u}_1$ and $\mathbf{u}_n$ denote the eigenvectors corresponding to the smallest and the largest eigenvalues of $\mathbf{W}$, respectively.
This paper investigates the probability density function (p.d.f.) of the random quantity $Z_{\ell}^{(n)}=\left|\mathbf{v}^\dagger\mathbf{u}_\ell\right|^2\in(0,1)$ for $\ell=1,n$. In particular, we derive a finite dimensional closed-form p.d.f. for $Z_{1}^{(n)}$ which is amenable to asymptotic analysis as $m,n$ diverges with $m-n$ fixed. It turns out that, in this asymptotic regime, the scaled random variable $nZ_{1}^{(n)}$ converges in distribution to $\chi^2_2/2(1+\theta)$, where $\chi_2^2$ denotes a chi-squared random variable with two degrees of freedom. This reveals that $\mathbf{u}_1$ can be used to infer information about the spike. 
On the other hand, the finite dimensional p.d.f. of  $Z_{n}^{(n)}$ is expressed as a double integral in which the integrand contains a determinant of a square matrix of dimension $(n-2)$. Although a simple solution to this double integral seems intractable, for special configurations of $n=2,3$, and $4$, we obtain closed-form expressions.

\end{abstract}


%
\IEEEpeerreviewmaketitle

\section{Introduction}
\IEEEPARstart{C}{onsider} a set of $m (\geq n)$ independent noisy observation vectors  $\mathbf{x}_j\in\mathbb{C}^{n\times 1}$ modeled as
\begin{align}
    \mathbf{x}_j=\sqrt{\theta}\mathrm{s}_j \mathbf{v}+ \mathbf{n}_j,\;\;j=1,2,\ldots,m,
\end{align}
where $\theta\geq 0$ is a non-random parameter, $\mathrm{s}_j\sim\mathcal{CN}(0,1)$ is the signal, $\mathbf{v}\in\mathbb{C}^{n\times 1}$ is an unknown non-random vector with $||\mathbf{v}||^2=1$,  $\mathbf{n}_j\sim\mathcal{CN}_n(\mathbf{0},\mathbf{I}_n)$ denotes the additive white Gaussian noise vector which is independent of $\mathrm{s}_j$, $||\cdot||$ denotes the Euclidean norm, and $\mathbf{I}_n$ stands for the $n\times n$ identity matrix. Given the above observations, one would seek to infer reliable information about the vector $\mathbf{v}$. To this end, it is customary to consider the sample covariance matrix $\mathbf{S}=\displaystyle \frac{1}{m}\sum_{j=1}^m \mathbf{x}_j\mathbf{x}_j^\dagger$ which is the simplest estimator for the population covariance matrix $\boldsymbol{\Sigma}=\mathbf{I}_n+\theta \mathbf{vv}^\dagger$, where $(\cdot)^\dagger$ denotes the conjugate transpose operator. Since $\mathbf{S}$ admits the eigen-decomposition $\mathbf{S}=\mathbf{U}\boldsymbol{\Lambda}_s\mathbf{U}^\dagger$, where $\mathbf{U}=\left(\mathbf{u}_1\; \mathbf{u}_2\;\ldots\; \mathbf{u}_n\right)\in\mathbb{C}^{n\times n}$ is a unitary matrix and $\boldsymbol{\Lambda}_s=\text{diag}\left(\lambda_1/m,\lambda_2/m,\ldots,\lambda_n/m\right)$ is a diagonal matrix with $0< \lambda_1\leq \lambda_2\leq \ldots\leq \lambda_n<\infty$ denoting the ordered eigenvalues of $m\mathbf{S}$, the squared modulus of the projection of the true vector $\mathbf{v}$ (i.e., spike) onto $\mathbf{u}_n$ (i.e., $|\mathbf{v}^\dagger\mathbf{u}_n|^2$) is commonly used to infer information about the latent vector $\mathbf{v}$ \cite{ref:paul,ref:bGeorges,ref:bloemendal,ref:haokai,ref:bao,ref:wWang}. On the other hand, in the finite dimensional setting, one would expect $\mathbf{u}_1$ to carry the least amount of information about $\textbf{v}$ (i.e., via leaked signal energy into the noise subspace). 

{\color{blue}The utility of the general random quantity $\left|\mathbf{v}^\dagger \mathbf{u}_\ell\right|,\; \ell=1,2,\ldots,n$, in various application domains has been highlighted in \cite{ref:vantrees, ref:harryLee, ref:denis, ref:fengFeng, ref:luYan, ref:johnPaul, ref:arthurLu, ref:huiZou, ref:liaoMahoney, ref:johnsonMestre, ref:yangEdgar, ref:fumito, ref:donohoGavishJohn,ref:donohoBehrooz,ref:remiMonasson, ref:bunPotters, ref:donohoeGavish, ref:huiJiaLi, ref:rish, ref:rCoulliet, ref:oliver, ref:nadler, ref:nobel, ref:rajPCInfo, ref:leeb, ref:rajDenoise, ref:gavishSing, ref:baoWang, ref:leebDeno, ref:wolf, ref:tuladhar,ref:wage,ref:sandeep,ref:gog1,ref:gog2,ref:gog3}. In particular, depending on the availability of the knowledge of $\mathbf{v}$, these application domains can be broadly classified into two groups: applications involving detection and estimation.  In this respect, to use $\left|\mathbf{v}^\dagger \mathbf{u}_\ell\right|$ as a test statistic for the detection of signals with certain orientation, one requires to have the full knowledge of $\mathbf{v}$. For instance, in array signal processing,  $\mathbf{v}$ is referred to as the nominal array steering vector which is known to the receiver \cite{ref:vantrees}. As such, to detect a desired signal arriving from the direction of the array steering vector, based on the eigenvectors of the sample covariance matrix, the test statistic $\left|\mathbf{v}^\dagger \mathbf{u}_n\right|$ has been used in the literature (see e.g., \cite[Section 7.8.2] {ref:vantrees}, \cite{ref:harryLee, ref:denis, ref:fengFeng, ref:luYan}, and references therein). 

Be that as it may, from the estimation perspective, the metric $\left|\mathbf{v}^\dagger \mathbf{u}_n\right|$ has been frequently used to measure the closeness of  the population eigenvector $\mathbf{v}$ to the sample eigenvector $\mathbf{u}_n$ (or a measure of informativeness of $\mathbf{u}_n$ with respect to the latent vector $\mathbf{v}$ \cite{ref:rajPCInfo}) when the {\it unknown} $\mathbf{v}$ is approximated by $\mathbf{u}_n$ in various studies, see e.g., \cite{ref:bloemendal, ref:johnPaul, ref:arthurLu}, \cite{ref:yangEdgar}, \cite{ref:nadler, ref:nobel, ref:rajPCInfo, ref:leeb, ref:rajDenoise, ref:gavishSing, ref:leebDeno, ref:tuladhar,ref:wage,ref:sandeep,ref:gog1,ref:gog2,ref:gog3}, and references therein. For instance, a concrete example in this respect is the  principal component analysis (PCA) in which the eigenvectors of the {\it unknown} population covariance matrix is approximated by the eigenvectors of the sample covariance matrix (i.e., empirical principal components). This metric has further been used in the covariance estimation based on the optimal shrinkage of the eigenvalues of the sample covariance matrix in the high dimensional setting when the unobserved population covariance matrix assumes the spiked structure \cite{ref:johnsonMestre}, \cite{ref:donohoBehrooz}, \cite{ref:donohoeGavish}, \cite{ref:leeb}, \cite{ref:wolf}. In particular, in the spiked population covariance setting, optimal shrinkage functions, which depend on the former metric, corresponding to various orthogonally invariant loss functions, such as the operator loss, Frobenius loss, entropy loss, Fr\'echet loss, and Stein's loss, have been derived in \cite{ref:donohoGavishJohn}. Recently, Yang et al. \cite{ref:yangEdgar} have used this metric to asses the performance of doing PCA after taking various random projections (a.k.a. sketches) of the observed noisy data having spiked population covariance structure in the high dimensional setting. A more recent extension of these concepts can be found in \cite{ref:fumito}.
Various versions of such {\it combined}-algorithms have been proposed in the early works of \cite{ref:rokhlin, ref:halkoa, ref:halkob}. It is worth mentioning that the random projection based techniques have been studied related to certain problems arising in linear regression \cite{ref:raskutti, ref:liuDobriban, ref:drineas}, ridge regression \cite{ref:liuDobriban}, \cite{ref:wangGittens}, classification \cite{ref:cannings}, and convex optimization \cite{ref:pilancia, ref:pilancib, ref:pilancic}. Asymptotically optimal bias corrections for sample eigenvalues in certain covariance estimation problems have been evaluated in \cite{ref:oliver} using the metrics $\left|\mathbf{v}^\dagger \mathbf{u}_\ell\right|,\; \ell=1,2,\ldots,n$, and more general extensions of these results can be found in \cite{ref:bunPotters}. The application of $\left|\mathbf{v}^\dagger \mathbf{u}_n\right|$ in matrix denoising, again in the high dimensional setting, can be found in  \cite{ref:donohoeGavish}, \cite{ref:nobel}, \cite{ref:leeb, ref:rajDenoise, ref:gavishSing, ref:baoWang, ref:leebDeno}. Moreover, the correlation between the sample and population eigenvectors in the high dimensional setting  has been instrumental in certain radar/array signal processing applications, see e.g., \cite{ref:tuladhar,ref:wage,ref:sandeep,ref:gog1,ref:gog2,ref:gog3} and references therein.

It is common to estimate the {\it unobservable} population spike $\mathbf{v}$ (i.e., the first population PC) by its empirical counterpart $\mathbf{u}_n$ (i.e., the leading eigenvector of the sample covariance matrix or the first sample PC). Since $\mathbf{v}$ and $\mathbf{u}_n$ span a plane, it is trivial that some signal energy has been leaked into the noise subspace. Therefore, the lower sample eigenvectors (i.e., $\mathbf{u}_{n-1},\ldots,\mathbf{u}_1$) are informative about $\mathbf{v}$ as well. For instance, a systematic asymptotic analysis of this phenomena based on the metric $\left|\mathbf{v}^\dagger \mathbf{u}_\ell\right|,\;\ell=1,2,\ldots,n$, has used information-theoretic and random matrix tools in \cite{ref:remiMonasson}. To shed some light into this issue, noting that $\mathbf{u}_\ell,\; \ell=1,2,\ldots,n$, forms an orthonormal basis, let us write 
\begin{align}
\mathbf{v}=\sum_{\ell=1}^n c_\ell \mathbf{u}_\ell
\end{align}
from which we obtain
   $ c_\ell=\mathbf{u}_\ell^\dagger \mathbf{v}$.
Therefore, the fraction of signal energy contained in the signal subspace (i.e., first sample PC) is $\left|\mathbf{u}_n^\dagger \mathbf{v}\right|^2$, whereas the noise subspace (i.e., other sample PCs) carries $1-\left|\mathbf{u}_n^\dagger \mathbf{v}\right|^2$ of the total energy. Of $1-\left|\mathbf{u}_n^\dagger \mathbf{v}\right|^2$, the worst eigenmode (i.e., last sample PC) carries $\left|\mathbf{u}_1^\dagger \mathbf{v}\right|^2$ of the total signal energy. Clearly, the distributions of the random quantities $\left|\mathbf{u}_n^\dagger \mathbf{v}\right|^2$ and $\left|\mathbf{u}_1^\dagger \mathbf{v}\right|^2$ can shed some light into how the spiked energy is statistically distributed between the signal and the noise subspaces in non-asymptotic setting. These facts further highlight the utility of $\left|\mathbf{u}_n^\dagger \mathbf{v}\right|^2$ and $\left|\mathbf{u}_1^\dagger \mathbf{v}\right|^2$ in the presence of an {\it unknown} spike $\mathbf{v}$.}

Since the high dimensional statistical characteristics of $|\mathbf{v}^\dagger \mathbf{u}_\ell|^2$, for $m\mathbf{S}\sim \mathcal{CW}_n\left(m, \mathbf{I}_n+\theta \mathbf{vv}^\dagger\right)$ with  $\ell=1,n$, have been analyzed extensively in the literature, see e.g., \cite{ref:paul,ref:bGeorges,ref:bloemendal,ref:haokai,ref:bao,ref:wWang} and references therein, in this study, we mainly focus on the  finite dimensional distributions of these random quantities when the population spike $\mathbf{v}$ is unknown.


Real and complex Wishart matrices play a central role in various scientific disciplines including multivariate analysis and high dimensional statistics \cite{ref:andersonIntroMultivariate,ref:muirhead,ref:johnstoneStatisticalChallenges}, random matrix theory \cite{ref:anderson,ref:bai,ref:edelman,ref:jack1,ref:jack2,ref:jack3, ref:zanella}, physics \cite{ref:mehta,ref:forresterLogGases}, and finance \cite{ref:potters} with numerous engineering applications \cite{ref:akemannOxfordHandbook}, particularly in signal processing and communications \cite{ref:couillet, ref:telatar, ref:tulino, ref:asendorf}. In this respect, of prominent interest are the eigenvalues and eigenvectors (more generally spectral projectors) of a Wishart matrix, especially with a certain nontrivial covariance structure. Among various population covariance structures, the spiked model\footnote{These spikes arise in various practical settings in different scientific disciplines. For instance, they correspond to the first few dominant factors in factor models arising in financial economics \cite{ref:fan,ref:ontaski}, the number of clusters in gene expression data \cite{ref:ke}, and the number of signals in detection and estimation \cite{ref:couillet,ref:onatskiSignal}.} introduced by Johnstone \cite{ref:johnstone} has been widely used in the literature to analyze the effects of having a few dominant trends or correlations in the covariance matrix. Capitalizing on this model for a single dominant correlation (i.e., $\boldsymbol{\Sigma}=\mathbf{I}_n+\theta \mathbf{vv}^\dagger$ or the rank-one perturbation of the identity), in the high dimensional setting, as $m,n\to\infty$ such that $n/m\to \gamma\in(0,1]$,  Baik et al. \cite{ref:baikPhaseTrans} proved that if $\theta>\sqrt{\gamma}$ (the supercritical regime) then the maximum eigenvalue of $\mathbf{S}$ converges to a Gaussian distribution, whereas it converges to a Tracy-Widom distribution \cite{ref:tracyLevelSp,ref:tracySymplectic} if $\theta<\sqrt{\gamma}$ (i.e., the subcritical regime). {\color{blue} To be specific, in the supercritical regime, the population spike generates the top sample eigenvalue of $\mathbf{S}$ outside the bulk, which follows the Marchenko-Pastur or ``quarter-circle" law \cite{ref:marchenko}, supported on a single interval $\left[\left(1-\sqrt{\gamma}\right)^2, \left(1+\sqrt{\gamma}\right)^2\right]$. Moreover, the  top sample eigenvalue converges almost surely to a limiting position $(1+\theta)\left(1
+\frac{\gamma}{\theta}\right)$ outside the upper edge of the bulk, thereby demonstrating an upward top sample eigenvalue bias. In contrast, in the subcritical regime, the top sample eigenvalue converges almost surely to the upper support of the bulk (i.e., $(1+\sqrt{\gamma})^2$).} This phenomenon is commonly known as the Baik, Ben Arous, P\'ech\'e (BBP) phase transition\footnote{The signal processing analogy of this phenomenon is known as the ``subspace swap" \cite{ref:johnsonMestre},\cite{ref:thomas}, \cite{ref:tuft}.} because of their seminal contribution in \cite{ref:baikPhaseTrans}. Subsequently, this analysis has been extended to various other random matrix ensembles, see e.g., \cite{ref:paul,ref:paulRandomMtrx,ref:bGeorges,ref:viragLimits2,ref:viragLimits1,ref:peche,ref:nadakuditi,ref:dWang} and references therein. These results reveal that, in the high dimensional setting,  the maximum eigenvalue of the sample covariance matrix can be used to detect a supercritical spike, whereas it cannot be used to detect a subcritical spike. However, this high dimensional behavior is not necessarily true in the finite dimensional setting.

Despite its utility, the high dimensional characteristics of the eigenvectors of Wishart matrices with spiked covariance matrices have received less attention in the early literature with notable exceptions in \cite{ref:paul,ref:hoyle,ref:nadler,ref:bGeorges}. In particular, it was shown in \cite{ref:bGeorges}\footnote{This has been first proved for real spiked-Wishart matrices in \cite{ref:paul}.} that, in the high dimensional setting as outlined above, $|\mathbf{v}^\dagger \mathbf{u}_n|^2$ converges almost surely to $\frac{1-\frac{\gamma}{\theta^2}}{1+\frac{\gamma}{\theta}}$ when $\theta>\sqrt{\gamma}$, whereas it converges almost surely to $0$ when $\theta<\sqrt{\gamma}$. This phase transition result reveals that no information about a subcritical spike can be inferred from the leading sample eigenvector. In a sharp contrast to this observation, a key recent result in \cite{ref:bloemendal} has established that a subcritical spike very close to the critical threshold can even cause $\textbf{u}_n$ to have a bias of small order towards $\mathbf{v}$. Capitalizing on this fact, the authors in \cite{ref:bloemendal} have concluded that a subcritical spike can be observed using the dominant eigenvector of the sample covariance matrix. {\color{blue}The asymptotic behavior of the extreme eigenvectors of not necessarily the Wishart type matrices has been studied on the level of the first order limit in \cite{ref:bGeorges, ref:georgeFob, ref:capitaineFo, ref:dingFo}. A comprehensive study on the eigenvector behavior of the sample covariance matrix in the full subcritical and supercritical regimes can be found in \cite{ref:bloemendal}.  In particular, as shown in \cite{ref:bloemendal}, the eigenvector distribution is similar to (up to appropriate scaling) that of the bulk and edge regimes of Wigner matrices without spikes; see \cite{ref:bourChi, ref:knowChi, ref:taoChi, ref:beniChi, ref:marChi} and references therein for instance. Moreover, in the subcritical regime, the limiting distribution of the square of the projected maximum eigenvector component (after proper scaling) is given by a Chi squared distribution, which implies the asymptotic Gaussianity of the eigenvector components \cite{ref:bloemendal}. Notwithstanding that the result was established for sample covariance matrix only in \cite{ref:bloemendal}, it can well be extended to deformed Wigner random matrices without essential difference \cite{ref:dWang}. In the supercritical regime, the fluctuation of the eigenvectors was recently studied in \cite{ref:bao}, \cite{ref:baoWang}, \cite{ref:capitaine}, \cite{ref:baoChi} for various other random matrix models. It is noteworthy that the aforementioned studies are restricted either to subcritical or to supercritical regime only, thereby leaving the critical regime unattended. To fill in this gap, recently, Bao and Wang \cite{ref:dWang} have focused on the critical regime of the  BBP phase transition and established the distribution\footnote{Here they have exploited the so called {\it eigenvector-eigenvalue} identity \cite{ref:terry} to derive the asymptotic distribution.} of the eigenvectors associated with the leading eigenvalues, however, for the unitary Gaussian ensemble with spiked external source only. In a sharp contrast, the behavior of the leading eigenvector, for $\gamma\to\infty$, has been analyzed in \cite{ref:feldman}}.
 A different kind of asymptotic analysis based on small noise perturbation approach \cite{ref:johnstoneRoys} is used in \cite{ref:prathapRoyRoot} to derive accurate stochastic approximations to $|\mathbf{v}^\dagger \mathbf{u}_n|^2$ for various   
complex Wishart matrices. 

Whereas the above technical results assume either diverging matrix dimensions (i.e., $m,n\to \infty$) or large signal-to-noise ratio (i.e., small noise power), these quantities assume relatively small values in many practical applications. In this respect, the above mentioned results pertaining to the sample eigenvectors may provide a poor approximation to the true behavior of the population eigenvectors. For instance, as depicted in \cite[Fig. 9.5]{ref:couillet}, the convergence rate of the empirical average of $|\mathbf{v}^\dagger \mathbf{u}_n|^2$ is very slow below the phase transition. {\color{blue}Therefore, despite the fact that, below the phase transition, $\mathbf{u}_n$ is asymptotically orthogonal to the population spike $\mathbf{v}$ (i.e., uninformative first order behavior of $\mathbf{u}_n$ with respect to $\mathbf{v}$), there exists a large range of values of $m$ and $n$ for which this orthogonality result does not hold.}
Against this backdrop, in the finite dimensional setting, it is plausible that $\boldsymbol{u}_1$ is also informative about the spike; no matter how small the amount of information is. Having motivated with these observations, in this paper, we characterize  the finite dimensional distributions of the leading (i.e., $\mathbf{u}_n$) and least eigenvectors (i.e., $\mathbf{u}_1$)  of the {\color{blue}matrix $m\mathbf{S}=\mathbf{W}\sim \mathcal{CW}_n\left(m,\mathbf{I}_n+\theta \mathbf{vv}^\dagger\right)$} via the distributions of the squared modulus of the eigen-projectors $|\mathbf{v}^\dagger \mathbf{u}_n|^2$ and $|\mathbf{v}^\dagger \mathbf{u}_1|^2$. Since a closed-form analytical probability density function (p.d.f.) for the joint density of the eigenvectors of $\mathbf{W}$ seems intractable, here we adopt a moment generating function (m.g.f.) based approach which requires the joint density of the eigenvalues and eigenvectors of $\mathbf{W}$ instead. The resultant matrix integral over the unitary manifold is further simplified using a contour integral approach due to \cite{ref:wang}. This key step transformed our problem into an equivalent form involving only the eigenvalues of $\mathbf{W}$ which is amenable to further analysis.

In particular, we leverage the powerful contour integral representation of unitary integrals \cite{ref:wang} and orthogonal polynomial techniques developed in \cite{ref:mehta} to derive a closed-form expression for the p.d.f. of $|\mathbf{v}^\dagger \mathbf{u}_1|^2$ which is valid for arbitrary $m,n$ and $\theta$. This result further indicates that the least eigenvector contains a certain amount of information about the spike. The resultant p.d.f. expression consists of a determinant of a square matrix whose dimensions depend on the relative difference between $m$ and $n$ (i.e., $m-n$). This key feature further facilitates the asymptotic analysis of $|\mathbf{v}^\dagger \mathbf{u}_1|^2$ in the regime in which $m,n\to \infty$ such that $m-n$ is fixed\footnote{This is also known as the microscopic limit in the literature of theoretical physics \cite{ref:forresterLogGases,ref:ghur1,ref:ghur2,ref:akemann}.}. It turns out that, in this particular regime, $|\mathbf{v}^\dagger \mathbf{u}_1|^2$ scales on the order $1/n$; therefore, $n |\mathbf{v}^\dagger \mathbf{u}_1|^2$ converges in distribution to $\displaystyle \frac{\chi_2^2}{2(1+\theta)}$, with $\chi_2^2$ denoting a chi-squared random variable with two degrees of freedom. This simple stochastic convergence result reveals that the least eigenvector contains information about the spike in this particular asymptotic domain irrespective of the value of $\theta$. Moreover, our numerical experiments indicate that this stochastic convergence result compares favourably with finite values of $m$ and $n$. Apart from these outcomes, we also derive an exact finite dimensional p.d.f. expression for $|\mathbf{v}^\dagger \mathbf{u}_n|^2$. This key p.d.f. result is expressed in terms of a double integral in which the integrand contains a determinant of a square matrix of size $n-2$. Although, an analytical closed-form solution to this integral seems intractable for general $n$, we have derived closed-form expressions for the special configurations of $n=2,3$, and $n=4$. Nevertheless, the double integral form of the p.d.f. can be evaluated numerically for arbitrary $m,n$, and $\theta$. To further illustrate the utility of the m.g.f. machinery developed in this manuscript beyond the statistical characterization of the extreme eigenvectors, we also present the analytical p.d.f. of $|\mathbf{v}^\dagger\mathbf{u}_2|^2$, where $\mathbf{u}_2$ denotes the eigenvector corresponding to the {\it second smallest} eigenvalue of $\mathbf{W}$.

{\color{blue}This paper also extends the former analytical framework to the real and singular Wishart cases as well.  Various application domains of the spectral characteristics of the real and singular Wishart matrices include signal processing, statistics, and econometrics, see e.g., \cite{ref:muirhead,ref:anderson,ref:AOSing,ref:srivastava,ref:ratnarajahSing,ref:ratnarajahMulti,ref:mallik,ref:shafi,ref:garcia,ref:goodall,ref:uhlig} and references therein. The main technical challenge pertaining to these extensions is that the resultant contour integrals, in general,  do not admit tractable simple closed-form expressions which are amenable to further analysis. For instance, in case of singular Wishart, the corresponding contour integral evaluates to an expression containing zonal polynomials 
However, as our detailed analysis reveals,  they can be expressed in closed-form for a few special configurations of $m$ and $n$.}

The remainder of this paper is organized as follows. Section II provides some key preliminary results required in the subsequent sections. The new exact p.d.f.s of  $|\mathbf{v}^\dagger \mathbf{u}_1|^2$ and  $|\mathbf{v}^\dagger \mathbf{u}_n|^2$ are derived in Section III, whereas a sketch of the proof of the p.d.f. of $|\mathbf{v}^\dagger \mathbf{u}_2|^2$ has also been given therein. Moreover, a detailed asymptotic analysis (i.e., the microscopic limit) of   $|\mathbf{v}^\dagger \mathbf{u}_1|^2$ is also provided in Section III. Further extension of the m.g.f. framework to real and singular Wishart scenarios is given in Section IV. Finally, conclusive remarks are made in Section IV.

\section{Preliminaries}

Here we present some fundamental results pertaining to the finite dimensional representation of correlated Wishart matrices, related unitary integrals and generalized Laguerre polynomials which are instrumental in our main derivations. 

\begin{definition}
  Let $\mathbf{X}\in\mathbb{C}^{n\times m}$ ($m\geq
  n$) be distributed as $ \mathcal{CN}_{n,m}\left(\mathbf{0},\boldsymbol{\Sigma}\otimes\mathbf{I}_m\right)$. Then the matrix $\mathbf{W}=\mathbf{X} \mathbf{X}^\dagger$ is said to follow a complex correlated Wishart distribution $\mathbf{W}\sim\mathcal{CW}_n(m,\boldsymbol{\Sigma})$ with p.d.f. \cite{ref:james}
  \begin{align}
  \label{eq wishart}
  f(\mathbf{W}) {\rm d}\mathbf{W}=
  \frac{{\det}^{m-n}(\mathbf{W})}{\tilde{\Gamma}_n(m) \det^m(\boldsymbol{\Sigma})}\;
\text{etr}\left(-\boldsymbol{\Sigma}^{-1}\mathbf{W}\right) {\rm d}\mathbf{W}
  \end{align}
\end{definition}
where {\color{blue} ${\rm d}\mathbf{W}$ is Lebesgue measure on the space of $n\times n$ Hermitian matrices, identifiable unambiguously with $\mathbb{R}^{n^2}$($=\mathbb{R}^n\times \mathbb{C}^{n(n-1)/2}$), defined by taking on-or-above diagonal entries as coordinates \cite{ref:anderson}}, $\tilde{\Gamma}_n(a)=\displaystyle \pi^{\frac{1}{2}n(n-1)}\prod_{j=1}^n \Gamma(a-j+1)$ represents the complex multivariate gamma function with $\Gamma(\cdot)$ denoting the classical gamma function, $\det(\cdot)$ is the determinant of a square matrix, and $\text{etr}(\cdot)\triangleq e^{\text{tr}(\cdot)}$ with $\text{tr}(\cdot)$ denoting the trace of a square matrix.
\begin{definition}
  Let $\mathbf{V}\in\mathbb{C}^{n\times n }$ and $\mathbf{T}\in\mathbb{C}^{n\times n}$ be two Hermitian non-negative definite matrices. Then the hypergeometric function of two matrix arguments is defined as \cite{ref:james}
   \begin{align*}
    {}_0\widetilde {\mathsf{F}}_0\left(\mathbf{V},\mathbf{T}\right)=\sum_{k=0}^\infty \frac{1}{k!} \sum_{\kappa} \frac{C_\kappa(\mathbf{V})C_\kappa(\mathbf{T})}{ C_\kappa(\mathbf{I}_n)}
  \end{align*} 
  where $C_\kappa(\cdot)$ is the complex zonal polynomial which is a symmetric, homogeneous polynomial of degree $k$ in the eigenvalues of the argument matrix\footnote{The exact algebraic definition of the zonal polynomial is tacitly avoided here, since it is not required in the subsequent analysis. More details of the zonal polynomials can be found in \cite{ref:james,ref:takemura}.}, $\kappa=(k_1,\ldots,k_n)$, with $k_i$'s being non-negative integers, is a partition of $k$ such that $k_1\geq\cdots\geq k_n\geq 0$ and $\sum_{i=1}^nk_i=k$.
\end{definition}
Moreover, {\color{blue}the Harish-Chandra–Itzykson-Zuber integral \cite{ref:harish, ref:itzykson, ref:james}} can be written as
\begin{align}
\label{eq split}
    {}_0\widetilde {\mathsf{F}}_0\left(\mathbf{V},\mathbf{T}\right)=\int\limits_{\mathcal{U}_n} \text{etr}\left(\mathbf{V}\mathbf{U}\mathbf{TU}^\dagger\right) {\rm d}\mathbf{U}
\end{align}
where ${\rm d}\mathbf{U}$ is the invariant measure on the unitary group $\mathcal{U}_n$ normalized to make the total measure unity (i.e., $\int_{\mathcal{U}_n} {\rm d}\mathbf{U}=1$). If at least one of the argument matrices assumes a rank-one structure, the hypergeometric function of two matrix arguments simplifies as shown in the following theorem which is due to \cite{ref:wang}\footnote{ A more generalized contour integral representation in this respect can be found in \cite{ref:peterJ}, \cite{ref:peter22}. However, the above form is sufficient for our work in this paper.}.
\begin{theorem}
\label{hcizint}
    Let $\mathbf{V}=\mathbf{v}\mathbf{v}^\dagger$ with $\mathbf{v}\in\mathbb{C}^{n\times 1}$ and $||\mathbf{v}||=1$. Then we have \cite{ref:wang}
    \begin{align}
    \label{eq cont integral}
        {}_0\widetilde {\mathsf{F}}_0\left(\mathbf{vv}^\dagger,\mathbf{T}\right)=\frac{(n-1)!}{2\pi \mathrm{i}} 
        \oint\limits_{\mathcal{C}}
        \frac{e^\omega}{\displaystyle \prod_{j=1}^n \left(\omega-\tau_j\right)} {\rm d}\omega
    \end{align}
    where $\tau_1,\tau_2,\ldots,\tau_n$ are the eigenvalues of $\mathbf{T}$, the contour $\mathcal{C}$ is large enough so that all $\tau_j$'s are in its interior, and $\rm i=\sqrt{-1}$.
\end{theorem}
  The following statistical characterization of the eigen-decomposition of $\mathbf{W}\sim\mathcal{CW}_n(m,\boldsymbol{\Sigma})$ is also useful in the sequel.
  \begin{theorem}
  Let the Hermitian positive definite matrix $\mathbf{W}$ assume the eigen-decomposition $\mathbf{W}=\mathbf{U}\boldsymbol{\Lambda} \mathbf{U}^\dagger$, where  $\boldsymbol{\Lambda}=\text{diag}\left(\lambda_1,\lambda_2,\ldots, \lambda_n\right)$ with $0<\lambda_1<\lambda_2<\ldots<\lambda_n<\infty$ denoting the ordered eigenvalues of $\mathbf{W}$ and $\mathbf{U}=\left(\mathbf{u}_1\; \mathbf{u}_2\; \ldots\; \mathbf{u}_n\right)$ contains the corresponding eigenvectors\footnote{The uniqueness of this representation and related implications are discussed in \cite{ref:anderson}.}. {\color{blue}Then the Jacobian of matrix transformation can be written as \cite{ref:james} ${\rm d}\mathbf{W}=\frac{\pi^{n(n-1)}}{\tilde{\Gamma}_n(n)}\Delta_n^2(\boldsymbol{\lambda}) {\rm d}\boldsymbol{\Lambda}{\rm d}\mathbf{U}$, where $\Delta_n(\boldsymbol{\lambda})=\displaystyle \prod_{1\leq i<j\leq n } \left(\lambda_j-\lambda_i\right)$ denotes the Vandermonde determinant and ${\rm d}\boldsymbol{\Lambda}=\displaystyle \prod_{j=1}^n {\rm d}\lambda_j$.} Therefore, we obtain
  \begin{align}
  \label{eq prod measure}
      f(\mathbf{W}){\rm d}\mathbf{W}=\frac{\pi^{n(n-1)}}{\tilde{\Gamma}_n(n)}\Delta_n^2(\boldsymbol{\lambda}) f\left(\mathbf{U}\boldsymbol{\Lambda}\mathbf{U}^\dagger\right) {\rm d}\boldsymbol{\Lambda}{\rm d}\mathbf{U}. 
  \end{align}
\end{theorem}

\begin{definition}
  For $\rho>-1$, the generalized Laguerre polynomial of degree $M$, $L^{(\rho)}_M(z)$, is given by \cite{ref:szego}
  \begin{equation}
   \label{lagdef}
    L^{(\rho)}_M(z)=\frac{(\rho+1)_M}{M!}\sum_{j=0}^{M}\frac{(-M)_j}{(\rho+1)_j}\frac{z^j}{j!},
  \end{equation}
 with its $k$\textsuperscript{th} derivative satisfying
  \begin{align}
    \label{lagderi}
    \frac{{\rm d}^k}{{\rm d}z^k}L^{(\rho)}_M(z)=(-1)^kL^{(\rho+k)}_{M-k}(z),
  \end{align}
  where $(a)_j=a(a+1)\ldots(a+j-1)$ with  $(a)_0=1$ denotes the Pochhammer symbol. 
\end{definition}
It is also noteworthy that the Pochhammer symbol admits, for $M\in\mathbb{Z}_+$,
\begin{align}
    (-M)_j=\left\{
    \begin{array}{cc}
    \frac{(-1)^j M!}{(M-j)!} & \text{if $0\leq j\leq M$}\\
    0 & \text{if $j>M$}.
    \end{array}
    \right.
\end{align}
Finally, we use the following compact notation to represent the
determinants of $N\times N$ block matrices:
\begin{equation}
 \begin{split}
    \det\left[a_{i}\;\; b_{i,j-1}\right]_{\substack{i=1,\ldots,N\\
    j=2,\dots,N}}&=\left|\begin{array}{ccccc}
    a_{1} & b_{1,1}& b_{1,2}& \ldots & b_{1,N-1}\\
      a_{2} & b_{2,1}& b_{2,2}& \ldots & b_{2,N-1}\\
      \vdots &\vdots & \vdots &\ddots & \vdots \\
      a_{N} & b_{N,1}& b_{N,2}& \ldots & b_{N,N-1}
    \end{array}\right|,
  \end{split}
\end{equation}
\begin{equation}
 \begin{split}
    \det\left[a_{i}\;\; b_{i,j-1} \;\; c_{i,k-3}\right]_{\substack{i=1,\ldots,N\\
    j=2,3\\k=4,\dots,N}}&=\left|\begin{array}{cccccc}
    a_{1} & b_{1,1}& b_{1,2}& c_{1,1} & \ldots & c_{1,N-1}\\
      a_{2} & b_{2,1}& b_{2,2}& c_{2,1} & \ldots & c_{2,N-1}\\
      \vdots &\vdots & \vdots & \vdots &\ddots & \vdots \\
      a_{N} & b_{N,1}& b_{N,2}& c_{N,1} & \ldots & c_{N,N-1}
    \end{array}\right|,
  \end{split}
\end{equation}
and $\det[a_{i,j}]_{i,j=1,\ldots,N}$ denotes the determinant of an $N\times N$ square matrix with its $(i,j)^{\text{th}}$ element given by $a_{i,j}$.

\section{Finite Dimensional Results for the Distributions of the Eigenvectors}
 Here we adopt an m.g.f. based approach to determine the p.d.f.s of the random variables $Z^{(n)}_n$ and $Z^{(n)}_1$. To this end, by definition, the m.g.f. can be written as
 \begin{align}
     \mathcal{M}_{Z^{(n)}_\ell}(s)=\mathbb{E}\left\{e^{-s|\mathbf{v}^\dagger \mathbf{u}_\ell|^2}\right\},\;\; \ell=1,n
 \end{align}
where $\mathbb{E}\{\cdot\}$ denotes the mathematical expectation evaluated with respect to the density of $\mathbf{u}_\ell$. To facilitate further analysis, noting that $|\mathbf{v}^\dagger \mathbf{u}_\ell|^2=\text{tr}\left(\mathbf{vv}^\dagger \mathbf{U}\mathbf{e}_\ell \mathbf{e}_\ell^T\mathbf{U}^\dagger\right)$ with $\mathbf{e}_\ell$ denoting the $\ell$th column of the $n\times n$ identity matrix, we may rewrite the m.g.f. as
\begin{align}
    \label{eq mgf Udep}
    \mathcal{M}_{Z^{(n)}_\ell}(s)=\mathbb{E}\left\{\text{etr}\left(-s\mathbf{vv}^\dagger \mathbf{U}\mathbf{e}_\ell \mathbf{e}_\ell^T\mathbf{U}^\dagger\right)\right\}
\end{align}
where the expectation is now taken with respect to the joint distribution of $\mathbf{U}$ and $\boldsymbol{\Lambda}$. Following (\ref{eq prod measure}), the above expectation assumes
\begin{align}
    \mathcal{M}_{Z^{(n)}_\ell}(s)=\frac{\pi^{n(n-1)}}{\tilde{\Gamma}_n(n)} \int\limits_{\mathcal{R}} \Delta_n^2(\boldsymbol{\lambda})\int\limits_{\mathcal{U}_n}
    f\left(\mathbf{U}\boldsymbol{\Lambda}\mathbf{U}^\dagger\right) \text{etr}\left(-s\mathbf{vv}^\dagger \mathbf{U}\mathbf{e}_\ell \mathbf{e}_\ell^T\mathbf{U}^\dagger\right){\rm d}\mathbf{U} \;{\rm d}\boldsymbol{\Lambda}
\end{align}
where $\mathcal{R}=\{0<\lambda_1<\lambda_2<\ldots<\lambda_n<\infty\}$. Since, for a single-spiked Wishart matrix,  $\boldsymbol{\Sigma}=\mathbf{I}_n+\theta \mathbf{vv}^\dagger$  with  $\boldsymbol{\Sigma}^{-1}=\mathbf{I}_n-\frac{\theta}{\theta+1}\mathbf{vv}^\dagger$, 
we use (\ref{eq wishart}) to further simplify the above matrix integral as
\begin{align}
\label{eq mgf cont}
   \mathcal{M}_{Z^{(n)}_\ell}(s)=K_{n,\alpha} (1-\beta)^{n+\alpha}
   \int\limits_{\mathcal{R}}
   \Delta_n^2(\boldsymbol{\lambda}) \prod_{j=1}^n \lambda_j^\alpha e^{-\lambda_j} \Psi_\ell(\boldsymbol{\lambda},s) {\rm d}\boldsymbol{\Lambda}
\end{align}
where  $K_{n,\alpha}=\displaystyle 1/\displaystyle\prod_{j=1}^n (n-j)!(n+\alpha-j)!$, $\alpha=m-n$, and
\begin{align}
\label{eq unitary comp}
    \Psi_\ell(\boldsymbol{\lambda}, s)=\int\limits_{\mathcal{U}_n} \text{etr}\left\{\mathbf{vv}^\dagger \mathbf{U}\left(\beta \boldsymbol{\Lambda}-s\mathbf{e}_\ell \mathbf{e}_\ell^T\right) \mathbf{U}^\dagger\right\} {\rm d} \mathbf{U}
\end{align}
with $\beta=\theta/(1+\theta)$. In view of the fact that the matrix  $\beta \boldsymbol{\Lambda}-s\mathbf{e}_\ell \mathbf{e}_\ell^T$ assumes two different diagonal structures depending on the value of $\ell$ as
\begin{align}
    \beta \boldsymbol{\Lambda}-s\mathbf{e}_\ell \mathbf{e}_\ell^T=\Biggl\{\begin{array}{cc}
    \text{diag}\left(\beta \lambda_1-s, \beta\lambda_2,\ldots, \beta \lambda_n\right) & \text{for $\ell=1$}\\
    \text{diag}\left(\beta \lambda_1, \beta\lambda_2,\ldots, \beta \lambda_n-s\right) & \text{ for $\ell=n$},
    \end{array}\Biggr.
\end{align}
we find it convenient to consider the two cases separately from this point onward.
\subsection{The P.D.F. of $Z^{(n)}_1$}
For $\ell=1$, (\ref{eq unitary comp}) specializes to 
\begin{align}
    \Psi_1(\boldsymbol{\lambda},s)=\int\limits_{\mathcal{U}_n} \text{etr}\left\{\mathbf{vv}^\dagger \mathbf{U}\left(\beta \boldsymbol{\Lambda}-s\mathbf{e}_1 \mathbf{e}_1^T\right) \mathbf{U}^\dagger\right\} {\rm d} \mathbf{U},
\end{align}
which can be simplified using (\ref{eq split}) and   (\ref{eq cont integral}) to yield
\begin{align}
    \Psi_1(\boldsymbol{\lambda},s)=\frac{(n-1)!}{2\pi \mathrm{i}} 
        \oint\limits_{\mathcal{C}}
        \frac{e^\omega}{\displaystyle \left(s+\omega-\beta\lambda_1\right)\prod_{j=2}^n \left(\omega-\beta \lambda_j\right)} {\rm d}\omega.
\end{align}
This in turn enables us to express (\ref{eq mgf cont}) as
\begin{align}
    \mathcal{M}_{Z^{(n)}_1}(s)=C^\beta_{n,\alpha}
   \int\limits_{\mathcal{R}}
   \Delta_n^2(\boldsymbol{\lambda}) \prod_{j=1}^n \lambda_j^\alpha e^{-\lambda_j} \frac{1}{2\pi \mathrm{i}} 
        \oint\limits_{\mathcal{C}}
        \frac{e^\omega}{\displaystyle \left(s+\omega-\beta\lambda_1\right)\prod_{j=2}^n \left(\omega-\beta \lambda_j\right)} {\rm d}\omega\; {\rm d}\boldsymbol{\Lambda}
\end{align}
where $C^\beta_{n,\alpha}=(n-1)!K_{n,\alpha} (1-\beta)^{n+\alpha}$. Now we take the inverse Laplace transform of both sides to yield
\begin{align}
\label{eq lap inv count int}
    f^{(\alpha)}_{Z^{(n)}_1}(z)=C^\beta_{n,\alpha}
   \int\limits_{\mathcal{R}} e^{\beta \lambda_1 z}
   \Delta_n^2(\boldsymbol{\lambda}) \prod_{j=1}^n \lambda_j^\alpha e^{-\lambda_j}
   \frac{1}{2\pi \mathrm{i}} 
        \oint\limits_{\mathcal{C}}
        \frac{e^{(1-z)\omega}}{\displaystyle \prod_{j=2}^n \left(\omega-\beta \lambda_j\right)} {\rm d}\omega\; {\rm d}\boldsymbol{\Lambda}
\end{align}
in which the innermost contour integral can be evaluated, for $n\geq 3$, with the help of the residue theorem to obtain
\begin{align}
   f^{(\alpha)}_{Z^{(n)}_1}(z)=\frac{C^\beta_{n,\alpha}}{\beta^{n-2}}
   \int\limits_{\mathcal{R}} 
   e^{\beta \lambda_1 z}
   \Delta_n^2(\boldsymbol{\lambda}) \prod_{j=1}^n \lambda_j^\alpha e^{-\lambda_j}
   \sum_{k=2}^{n}
   \frac{e^{(1-z)\beta \lambda_k}}{\displaystyle \prod_{\substack{i=2\\i\neq k}}^{n} \left(\lambda_k-\lambda_i\right)}\;
  {\rm d}\boldsymbol{\Lambda}.
\end{align}
Now let us rewrite the $n$-fold integral, keeping the integration with respect to $\lambda_1$ last, as
\begin{align}
\label{eq z1pdf int}
  f^{(\alpha)}_{Z^{(n)}_1}(z)=\frac{C^\beta_{n,\alpha}}{\beta^{n-2}}  
  \int_{0}^\infty \lambda_1^\alpha e^{-(1-\beta z)\lambda_1}
  \Phi(\lambda_1,z) {\rm d}\lambda_1
\end{align}
where
\begin{align}
    \Phi(\lambda_1,z)=\int\limits_{\mathcal{S}} \sum_{k=2}^{n}
   \frac{e^{(1-z)\beta \lambda_k}}{\displaystyle \prod_{\substack{i=2\\i\neq k}}^{n} \left(\lambda_k-\lambda_i\right)}
   \Delta_{n-1}^2(\boldsymbol{\lambda})
   \prod_{j=2}^n \lambda_j^\alpha e^{-\lambda_j}
   \left(\lambda_j-\lambda_1\right)^2 {\rm d}\lambda_j
    \end{align}
 in which we have used the decomposition $\Delta_n^2(\boldsymbol{\lambda})=\prod_{j=2}^n (\lambda_j-\lambda_1)^2 \Delta_{n-1}^2(\boldsymbol{\lambda})$ and $\mathcal{S}=\{\lambda_1<\lambda_2<\ldots<\lambda_n<\infty\}$. For convenience, we relabel the variables as $\lambda_1=x$ and $x_j=\lambda_{j-1}$, $j=2,3,\ldots,n$, to arrive at  
\begin{align}
\label{eq pdf cvx}
  f^{(\alpha)}_{Z^{(n)}_1}(z)=\frac{C^\beta_{n,\alpha}}{\beta^{n-2}}  
  \int_{0}^\infty x^\alpha e^{-(1-\beta z)x}
  \Phi(x,z) {\rm d}x
\end{align}
where
\begin{align}
    \Phi(x,z)=\int\limits_{\mathcal{S}_x} \sum_{k=1}^{n-1}
   \frac{e^{(1-z)\beta x_k}}{\displaystyle \prod_{\substack{i=1\\i\neq k}}^{n-1} \left(x_k-x_i\right)}
   \Delta_{n-1}^2(\boldsymbol{x})
   \prod_{j=1}^{n-1} x_j^\alpha e^{-x_j}
   \left(x_j-x\right)^2 {\rm d}x_j
    \end{align}
with $\mathcal{S}_x=\{x<x_1<x_2<\ldots<x_{n-1}<\infty\}$. Since the integrand in the above $(n-1)$-fold integral is symmetric in $x_1,x_2,\ldots,x_{n-1}$, we may remove the ordered region of integration to obtain
\begin{align}
    \Phi(x,z)=\frac{1}{(n-1)!}\int\limits_{(x,\infty)^{n-1}} \sum_{k=1}^{n-1}
   \frac{e^{(1-z)\beta x_k}}{\displaystyle \prod_{\substack{i=1\\i\neq k}}^{n-1} \left(x_k-x_i\right)}
   \Delta_{n-1}^2(\boldsymbol{x})
   \prod_{j=1}^{n-1} x_j^\alpha e^{-x_j}
   \left(x_j-x\right)^2 {\rm d}x_j.
    \end{align}
Consequently, we can observe that the each term in the above summation evaluates to the same amount. Therefore, capitalizing on that observation, we may simplify the above $(n-1)$-fold integral to yield

    \begin{align}
    \Phi(x,z)=\frac{1}{(n-2)!}\int\limits_{(x,\infty)^{n-1}}
   \frac{e^{(1-z)\beta x_1}}{\displaystyle \prod_{{i=2}}^{n-1} \left(x_1-x_i\right)}
   \Delta_{n-1}^2(\boldsymbol{x})
   \prod_{j=1}^{n-1} x_j^\alpha e^{-x_j}
   \left(x_j-x\right)^2 {\rm d}x_j.
    \end{align}
To facilitate further analysis, we introduce the variable transformations $y_j=x_j-x$, $j=1,2,\ldots,n-1$, to  the above integral to arrive at
\begin{align}
    \Phi(x,z)=\frac{e^{-(n+\beta z-\beta-1)x}}{(n-2)!}\int\limits_{(0,\infty)^{n-1}}
    \frac{e^{(1-z)\beta y_1}}{\displaystyle \prod_{{i=2}}^{n-1} \left(y_1-y_i\right)}
    \Delta_{n-1}^2(\boldsymbol{y})
    \prod_{j=1}^{n-1} (y_j+x)^\alpha y_j^2 e^{-y_j}
    {\rm d}y_j,
\end{align}
from which we obtain by keeping the integration with respect to $y_1$ last
\begin{align}
\label{Phi int form}
    \Phi(x,z)=\frac{e^{-(n+\beta z-\beta-1)x}}{(n-2)!}\int_0^\infty
    e^{-(1-(1-z)\beta) y_1} y_1^2 (y_1+x)^\alpha
    \mathcal{Q}_{n-2}(y_1,x) {\rm d}y_1
\end{align}
where
{\color{blue}
\begin{align}
    \mathcal{Q}_{n-2}(y_1,x)=
    \int\limits_{(0,\infty)^{n-2}}
    \frac{\Delta_{n-1}^2(\boldsymbol{y})}{\displaystyle \prod_{{i=2}}^{n-1} \left(y_1-y_i\right)}
    \prod_{j=2}^{n-1} (y_j+x)^\alpha y_j^2 e^{-y_j}
    {\rm d}y_j.
\end{align}
}
To further simplify $\mathcal{Q}_{n-2}(y_1,x)$, noting the decomposition {\color{blue} $\Delta^2_{n-1}(\boldsymbol{y})=\prod_{j=2}^{n-1}(y_j-y_1)^2\Delta_{n-2}^2(\boldsymbol{y})$, we may relabel the variables as $z_{j-1}=y_{j}, j=2,3,\ldots,n-1$,} with some algebraic manipulation to arrive at
\begin{align}
\label{Qintdef}
    \mathcal{Q}_{n-2}(y_1,x)=
    \int\limits_{(0,\infty)^{n-2}}
    \Delta_{n-2}^2(\boldsymbol{z}) 
    \prod_{j=1}^{n-2} (y_1-z_j)(z_j+x)^\alpha z_j^2 e^{-z_j}
    {\rm d}z_j.
\end{align}
The above multiple integral can be solved following the orthogonal polynomial approach due to Mehta \cite{ref:mehta}, as shown in (\ref{RQ}) and (\ref{Rans}) of Appendix A, to yield
\begin{align}
\label{eq Q solu}
   \mathcal{Q}_{n-2}(y_1,x)=
   \frac{(-1)^n \hat{K}_{n-2,\alpha}}{(x+y_1)^\alpha} 
   \det\left[L^{(2)}_{n+i-3}(y_1)\;\;\; L^{(j)}_{n+i-j-1}(-x)\right]_{\substack{i=1,\ldots,\alpha+1\\
    j=2,\dots,\alpha+1}}
\end{align}
where 
\begin{align*}\hat{K}_{n-2,\alpha}=\prod_{j=1}^{\alpha+1}(n+j-3)! \prod_{j=0}^{n-3}(j+1)!(j+2)!/\prod_{j=0}^{\alpha-1} j!.
\end{align*}
{\color{blue}It is noteworthy that, when $\alpha=0$, we interpret the above determinant as $L^{(2)}_{n-2}(y_1)$, since the $(\alpha+1) \times (\alpha+1)$ matrix degenerates to a scalar. Moreover, an empty product is interpreted as unity. } 
Now we substitute $\mathcal{Q}_{n-2}(y_1,x)$ in (\ref{eq Q solu}) into (\ref{Phi int form}) with some algebraic manipulation to obtain
\begin{align}
\label{Phi int form 1}
    \Phi(x,z)=\frac{(-1)^n \hat{K}_{n-2,\alpha}}{(n-2)!}  e^{-(n+\beta z-\beta-1)x}\int_0^\infty
   & e^{-(1-(1-z)\beta) y_1} y_1^2\nonumber\\
    & \times \det\left[L^{(2)}_{n+i-3}(y_1)\;\;\; L^{(j)}_{n+i-j-1}(-x)\right]_{\substack{i=1,\ldots,\alpha+1\\
    j=2,\dots,\alpha+1}} {\rm d}y_1.
\end{align}
Since only the first column of the determinant depends on $y_1$, we can conveniently rewrite the above integral as
\begin{align}
\label{Phi int form 2}
    \Phi(x,z)=&\frac{(-1)^n \hat{K}_{n-2,\alpha}}{(n-2)!}  e^{-(n+\beta z-\beta-1)x} \nonumber\\
    & \qquad \qquad \times 
    \det\left[\int_0^\infty
    e^{-(1-(1-z)\beta) y_1}y_1^2L^{(2)}_{n+i-3}(y_1){\rm d}y_1\;\;\; L^{(j)}_{n+i-j-1}(-x)\right]_{\substack{i=1,\ldots,\alpha+1\\
    j=2,\dots,\alpha+1}}.
\end{align}
{\color{blue}Now we use \cite[Eq. 7.414.8]{ref:gradshteyn} to evaluate the above integral and substitute back the resultant $\Phi(x,z)$ into (\ref{eq z1pdf int}) along with some algebraic manipulations to arrive at
\begin{align}
\label{eq min thm 1step}
     f^{(\alpha)}_{Z^{(n)}_1}(z)=\frac{(n-1)!}{(n+\alpha-1)!}(1-\beta)^{n+\alpha}
    \int_0^\infty
    e^{-(n-\beta)x} x^\alpha 
    \det\left[\zeta_i(z,\beta)\;\;\; L^{(j)}_{n+i-j-1}(-x)\right]_{\substack{i=1,\ldots,\alpha+1\\
    j=2,\dots,\alpha+1}} {\rm d}x
\end{align}
where 
\begin{align}
\label{def zeta mod}
    \zeta_i(z,\beta)=
    \frac{(n+i-1)!}{(n+i-3)!} \frac{(-\beta)^{i-1}(1-z)^{n+i-3}}
    {\left[1-\beta(1-z)\right]^{n+i}}.
\end{align}
To facilitate further analysis, noting that the columns denoted by $j=2,3,\ldots, \alpha+1$, depend only on $x$, 
we use (\ref{lagdef}) to rewrite the Laguerre polynomials in the determinant with some algebraic manipulation as
\begin{align}
\label{det lag decom}
    &\det\left[\zeta_i(z,\beta)\;\;\; L^{(j)}_{n+i-j-1}(-x)\right]_{\substack{i=1,\ldots,\alpha+1\\
    j=2,\dots,\alpha+1}}\nonumber\\
    & =\det\left[\zeta_i(z,\beta)\;\;\; \frac{(n+i-1)!}{(n+i-j-1)! j!}\sum_{k_j=0}^{n+i-j-1}\frac{(-n-i+j+1)_{k_j}}{(j+1)_{k_j}}\frac{(-x)^{k_j}}{k_j!} \right]_{\substack{i=1,\ldots,\alpha+1\\
    j=2,\dots,\alpha+1}}.
\end{align}
Further manipulation in this form is highly undesirable due to the dependence of the summation upper limits on $i$ and $j$. To circumvent this difficulty, in view of making the summation upper limits depend only on $j$, we use the decomposition
\begin{align}
    (-n-\alpha+j)_{k_j}\frac{(-n-i+j+1)_{k_j}}{(-n-\alpha+j)_{k_j}}=(-n-\alpha+j)_{k_j}\frac{(n+i-j-1)! (n+\alpha-j-k_j)!}{(n+i-j-1-k_j)!(n+\alpha-j)!},
\end{align}
in (\ref{det lag decom}) with some algebraic manipulation to obtain
\begin{align}
    &\det\left[\zeta_i(z,\beta)\;\;\; L^{(j)}_{n+i-j-1}(-x)\right]_{\substack{i=1,\ldots,\alpha+1\\
    j=2,\dots,\alpha+1}}\nonumber\\
    & =\det\left[\zeta_i(z,\beta)\;\;\; \frac{(n+i-1)!}{(n+\alpha-j)! j!}\sum_{k_j=0}^{n+\alpha-j}\frac{(-n-\alpha+j)_{k_j}}{(j+1)_{k_j}}\frac{(-x)^{k_j}}{k_j!} \frac{ (n+\alpha-j-k_j)!}{(n+i-j-1-k_j)!} \right]_{\substack{i=1,\ldots,\alpha+1\\
    j=2,\dots,\alpha+1}}.
\end{align}
Now we may take the common factors out from the determinant and apply some algebraic manipulation to yield
\begin{align}
\label{min det simp}
   &\det\left[\zeta_i(z,\beta)\;\;\; L^{(j)}_{n+i-j-1}(-x)\right]_{\substack{i=1,\ldots,\alpha+1\\
    j=2,\dots,\alpha+1}}\nonumber\\
    & =  \frac{(n+\alpha-1)! (n+\alpha)!}{(n-1)!}
    \sum_{k_2=0}^{n+\alpha-2} \sum_{k_3=0}^{n+\alpha-3}\ldots
     \sum_{k_{\alpha+1}=0}^{\alpha-1}
     \prod_{j=2}^{\alpha+1} \frac{(n+\alpha-j)!}{(j+k_j)! k_j!} x^{\sum_{j=2}^{\alpha+1}k_j}\nonumber\\
     &\hspace{5.5cm}\times 
     \det\left[\frac{\zeta_i(z,\beta)}{(n+i-1)!}\;\;\; \frac{ 1}{\Gamma(n+i-j-k_j)} \right]_{\substack{i=1,\ldots,\alpha+1\\
    j=2,\dots,\alpha+1}}.
\end{align}
Finally, we use (\ref{min det simp}) and (\ref{def zeta mod}) in (\ref{eq min thm 1step}) and perform the integration with respect to $x$ followed by shifting the indices from $i,j$ to $i-1, j-1$ with some algebraic manipulation to obtain the p.d.f. of $Z_1^{(n)}$, for $n\geq 3$, as shown in Theorem \ref{thm min vec}.

For $n=2$  case, (\ref{eq lap inv count int}) specializes with the help of residue theorem to 
\begin{align}
    f^{(\alpha)}_{Z^{(2)}_1}(z)=C^\beta_{2,\alpha}
   \int\limits_{\lambda_1<\lambda_2<\infty} e^{\beta \lambda_1 z}
   \Delta_2^2(\boldsymbol{\lambda}) \prod_{j=1}^2 \lambda_j^\alpha e^{-\lambda_j}
        e^{(1-z)\beta \lambda_2\omega} {\rm d}\lambda_1{\rm d}\lambda_2
\end{align}
from which we obtain after interchanging the order of integration followed by some algebraic manipulations
\begin{align}
\label{min n2 double}
    f^{(\alpha)}_{Z^{(2)}_1}(z)=C^\beta_{2,\alpha} \int_0^\infty \lambda_1^\alpha e^{-(1-\beta z)\lambda_1}
    \int_{\lambda_1}^\infty
    (\lambda_2-\lambda_1)^2 \lambda_2^\alpha e^{-(1-(1-z)\beta)\lambda_2} {\rm d}\lambda_2 {\rm d}\lambda_1.
\end{align}
Now we can evaluate the inner integral to obtain
\begin{align}
    \int_{\lambda_1}^\infty
    (\lambda_2-\lambda_1)^2 \lambda_2^\alpha e^{-(1-(1-z)\beta)\lambda_2} {\rm d}\lambda_2=
    \sum_{k=0}^\alpha \frac{(k+2)! \alpha !}{(\alpha-k)! k!}
    \frac{\lambda_1^{\alpha-k}e^{-\lambda_1(1-(1-z)\beta)}}{\left[1-\beta(1-z)\right]^{k+3}},
\end{align}
which upon substituting into (\ref{min n2 double}) followed by integration with respect to $\lambda_2$ with some algebraic manipulation gives the p.d.f. of $Z_1^{(n)}$, for $n=2$, as shown in Theorem \ref{thm min vec}.}

\begin{theorem}\label{thm min vec}
    Let $\mathbf{W}\sim\mathcal{CW}_n\left(m,\mathbf{I}_n+\theta \mathbf{v}\mathbf{v}^\dagger\right)$ with $||\mathbf{v}||=1$ and $\theta\geq 0$. Let $\mathbf{u}_1$ be the eigenvector corresponding to the smallest eigenvalue of $\mathbf{W}$. Then the p.d.f. of $Z^{(n)}_1=|\mathbf{v}^\dagger \mathbf{u}_1|^2\in(0,1)$ is given by
    \begin{align}
    \label{thm main}
        f^{(\alpha)}_{Z^{(n)}_1}(z)=\left\{
        \begin{array}{ll}
        (1-\beta)^{n+\alpha} 
       \displaystyle \frac{(1-z)^{n-2}}{[1-\beta(1-z)]^{n+1}}
       \mathcal{H}_n^{(\alpha, \beta)}(z)
          & \text{if $n\geq 3$}\\
       \displaystyle \frac{(1-\beta)^{2+\alpha} }{(\alpha+1)!}
        \displaystyle\sum_{k=0}^\alpha 
        \frac{(k+2)!(2\alpha-k)!}{k!(\alpha-k)!} \frac{1}{(2-\beta)^{2\alpha-k+1}[1-\beta(1-z)]^{k+3}} & \text{if $n=2$}
        \end{array}\right.
    \end{align}
    where 
    \begin{align}
    \label{thm main sup}
         \mathcal{H}_n^{(\alpha, \beta)}(z)=\sum_{k_1=0}^{n+\alpha-2}\sum_{k_2=0}^{n+\alpha-3}\ldots \sum_{k_{\alpha}=0}^{n-1} &
        \frac{\left(\alpha+\sum_{j=1}^{\alpha}k_j\right)!}{(n-\beta)^{\alpha+\sum_{j=1}^{\alpha}k_j+1}}
        \prod_{j=1}^{\alpha} \frac{(n+\alpha-j-1)!}{(j+k_j+1)! k_j!}\nonumber\\
        & \times 
        \det\left[\frac{(n+\alpha)!(-\beta)^i(1-z)^i}{(n+i-2)![1-\beta(1-z)]^i}\;\;\;\;a_{i,j}(k_j)
        \right]_{\substack{i=0,\ldots,\alpha\\
    j=1,\dots,\alpha}}
    \end{align}
    with $a_{i,j}(\ell)=1/\Gamma\left(n+i-j-\ell\right)$ and $\beta=\theta/(\theta+1)$.
\end{theorem}

The above p.d.f. expression depends on $\theta$ via $\beta$ as we clearly see. This observation further reveals that, in the finite dimensional setting,  the least eigenvector also contains a certain amount of information about the spike $\mathbf{v}$. As we shall see below in Corollary \ref{cor min asy}, the least eigenvector retains this ability even in a particular asymptotic domain as well.

It is noteworthy that, since  the number of nested summations in Theorem \ref{thm min vec} depends only on $\alpha$, it provides an efficient way of evaluating the p.d.f. of $Z^{(n)}_1$, particularly for small values of $\alpha$. As such, 
for some small values of $\alpha$, (\ref{thm main}) admits the following simple forms.
\begin{corollary}
The exact p.d.f.s of $Z^{(n)}_1$ corresponding to $\alpha=0$ and $\alpha=1$, for $n\geq 3$, are given, respectively, by 
\begin{align}
    f^{(0)}_{Z^{(n)}_1}(z)&=
        \displaystyle \frac{n(n-1)(1-\beta)^n(1-z)^{n-2}}{(n-\beta)[1-\beta(1-z)]^{n+1}}\\
     f^{(1)}_{Z^{(n)}_1}(z) &= 
     \frac{n(n^2-1)(1-\beta)^{n+1}(1-z)^{n-2}}{2(n-\beta)^2[1-\beta(1-z)]^{n+1}} \left[
     {}_2F_1\left(-n+1,2;3;-\frac{1}{n-\beta}\right)
     \right.\nonumber\\
    &\hspace{6cm}\qquad + \left.
    \frac{\beta (1-z)}{1-\beta(1-z)}
    {}_2F_1\left(-n+2,2;3;-\frac{1}{n-\beta}\right)\right]
\end{align}
where ${}_2F_1\left(a,b;c;z\right)$ is the Gauss hypergeometric function \cite{ref:erdelyi}.
\end{corollary}
\begin{remark}
  The p.d.f. corresponding to $\theta=0$ (i.e., $\beta=0$) can be obtained from (\ref{thm main}) after some tedious algebraic manipulations. Alternatively, for $\beta=0$, noting that $\zeta_i(z,0)=0$, $i=2,3,\ldots,\alpha+1$, (\ref{eq min thm 1step}) can be simplified to yield
  \begin{align}
      f^{(\alpha)}_{Z^{(n)}_1}(z)=\frac{(n-1)n!}{(n+\alpha-1)!} (1-z)^{n-2}
    \int_0^\infty
    e^{-nx} x^\alpha 
    \det\left[ L^{(j)}_{n+i-j-1}(-x)\right]_{i,j=2,\ldots,\alpha+1} {\rm d}x,
  \end{align}
  from which we obtain after relabeling the indices $k=i-1$ and $\ell=j-1$
  \begin{align}
      f^{(\alpha)}_{Z^{(n)}_1}(z)=\frac{(n-1)n!}{(n+\alpha-1)!} (1-z)^{n-2}
    \int_0^\infty
    e^{-nx} x^\alpha 
    \det\left[ L^{(\ell+1)}_{n+k-\ell-1}(-x)\right]_{k,\ell=1,\ldots,\alpha} {\rm d}x.
  \end{align}
  Since the p.d.f. of $\lambda_1$, for $\beta=0$, assumes the form \cite[Eq. 3.8]{ref:prathapSIAM}
  \begin{align}
      f_{\lambda_1}(x)=\frac{n!}{(n+\alpha-1)!} e^{-nx} x^\alpha 
    \det\left[ L^{(\ell+1)}_{n+k-\ell-1}(-x)\right]_{k,\ell=1,\ldots,\alpha},
  \end{align}
  we easily obtain the desired answer $f^{(\alpha)}_{Z_1^{(n)}}(z)=(n-1)(1-z)^{n-2}$ as expected\footnote{It is well known that, for $\theta=0$, $\mathbf{U}$ is Haar distributed (i.e., uniformly distributed over the unitary manifold). Consequently, some algebraic manipulations will  establish the result $f^{(\alpha)}_{Z_1^{(n)}}(z)=(n-1)(1-z)^{n-2}$. It is also noteworthy that, since any permutation matrix is orthogonal, the same p.d.f. result holds for a general class of random variables of the form $Z_\ell^{(n)}=|\mathbf{v}^\dagger \mathbf{u}_\ell|^2, \; \ell=1,2,\ldots,n$.  }.
\end{remark}

Figure 1 compares the analytical p.d.f. result for $Z^{(n)}_1$ with simulated data. Analytical curves are generated based on Theorem \ref{thm min vec}. As can be seen from the figure, our analytical results match with the simulated data, thus validating our theorem. Quantile-Quantile (Q-Q) plots have been provided in Fig. \ref{fig_qq_mineigvec_different_n} to further compare the analytical and simulated data. Moreover, the effect of $\theta$ on the p.d.f. of $Z^{(n)}_1$ is depicted in Fig. \ref{fig_mineigvec_different_theta}, while the corresponding Q-Q plots are given in Fig. \ref{fig_qq_mineigvec_different_theta}.

\begin{figure}[htb]
    \centering
    \includegraphics[width=0.9\textwidth]{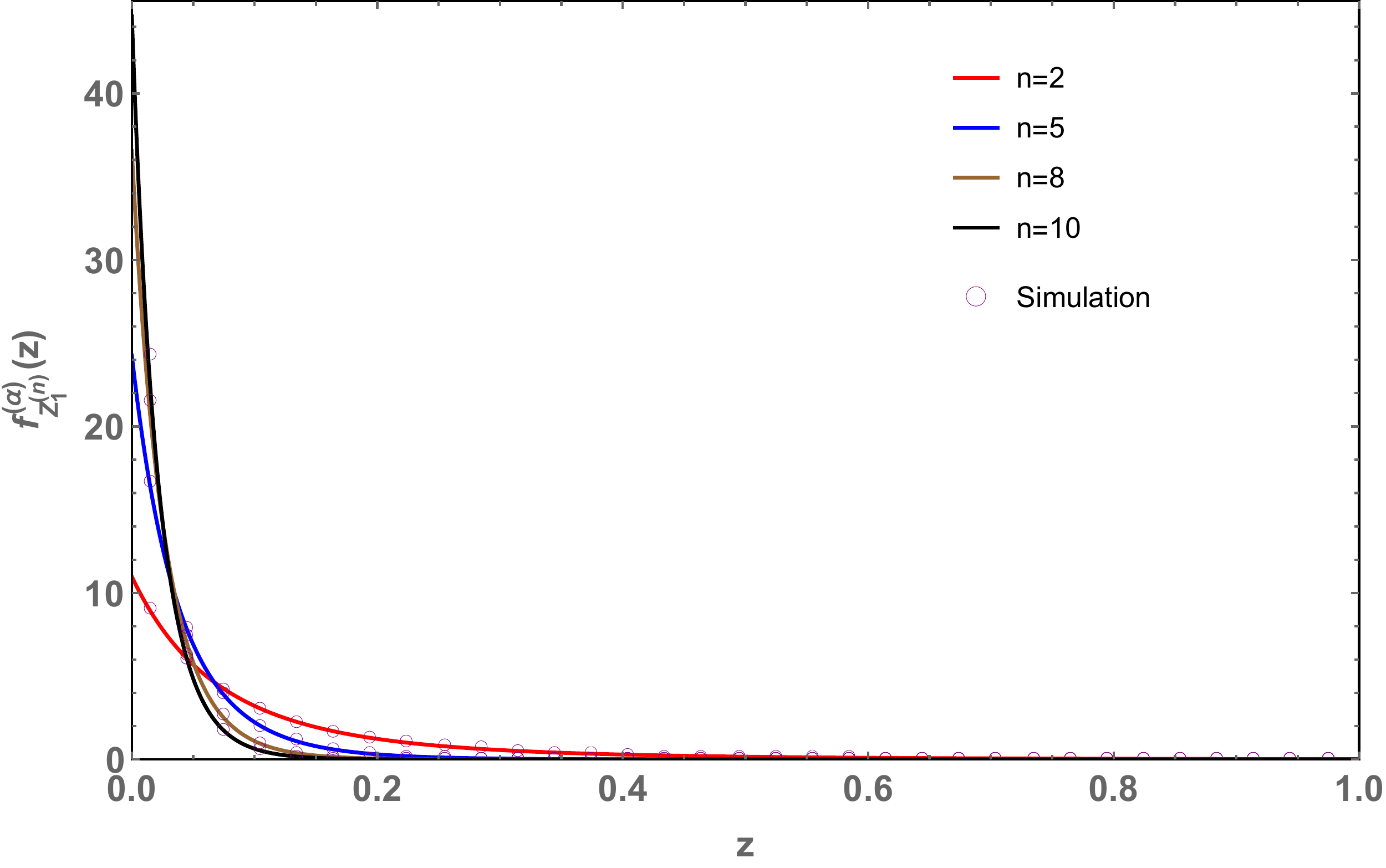}
    \caption{Comparison of simulated data points and the analytical p.d.f. $f^{(\alpha)}_{Z^{(n)}_1}(z)$ for different values of $n$ with $\alpha=2$ and $\theta=3$.}
    \label{fig_mineigvec_different_n}
\end{figure}

\begin{figure}
    \centering
    \begin{subfigure}[b]{.45\textwidth}
    \includegraphics[width=\textwidth]{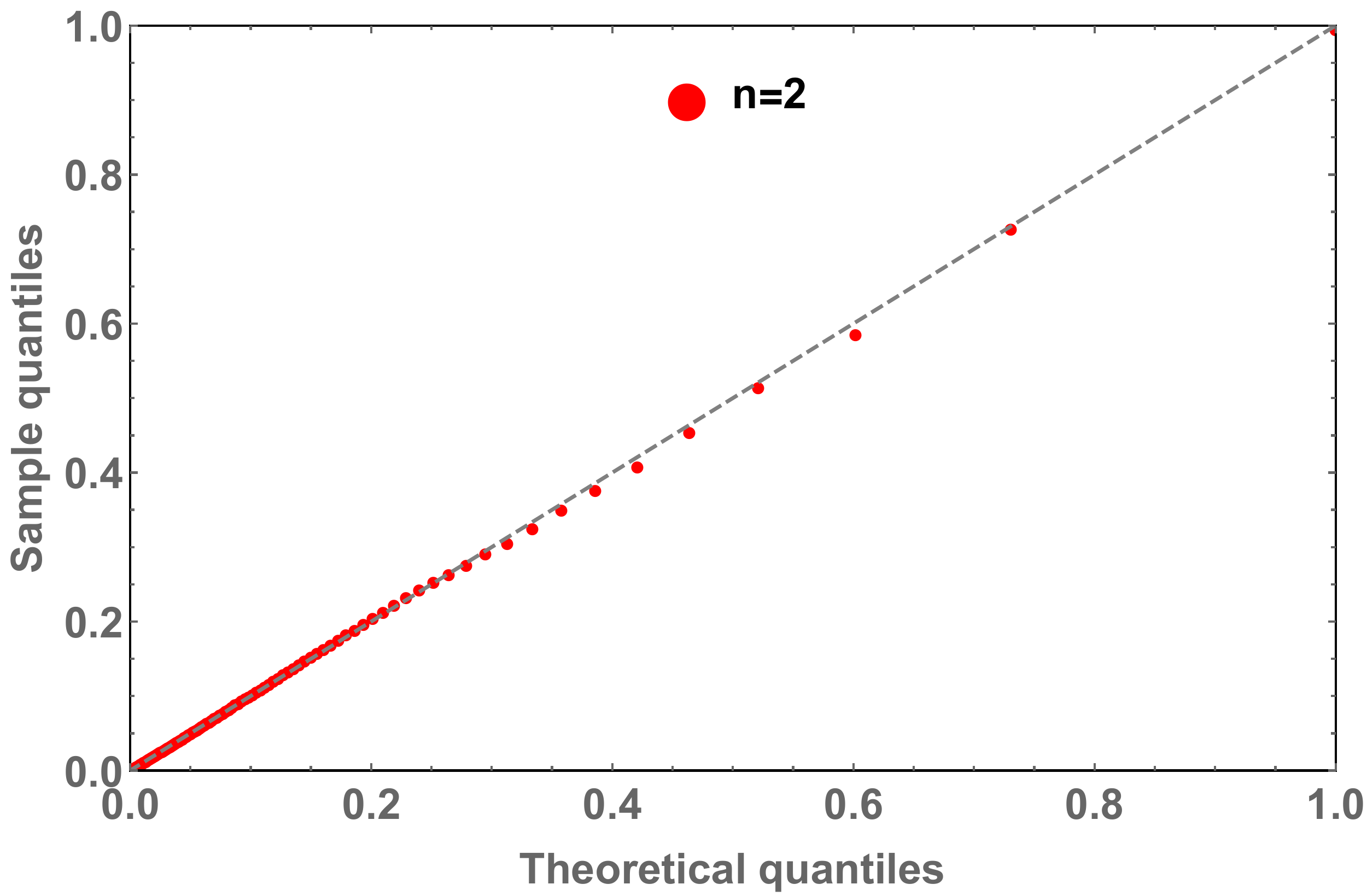}
    \end{subfigure}
    \begin{subfigure}[b]{.45\textwidth}
    \includegraphics[width=\textwidth]{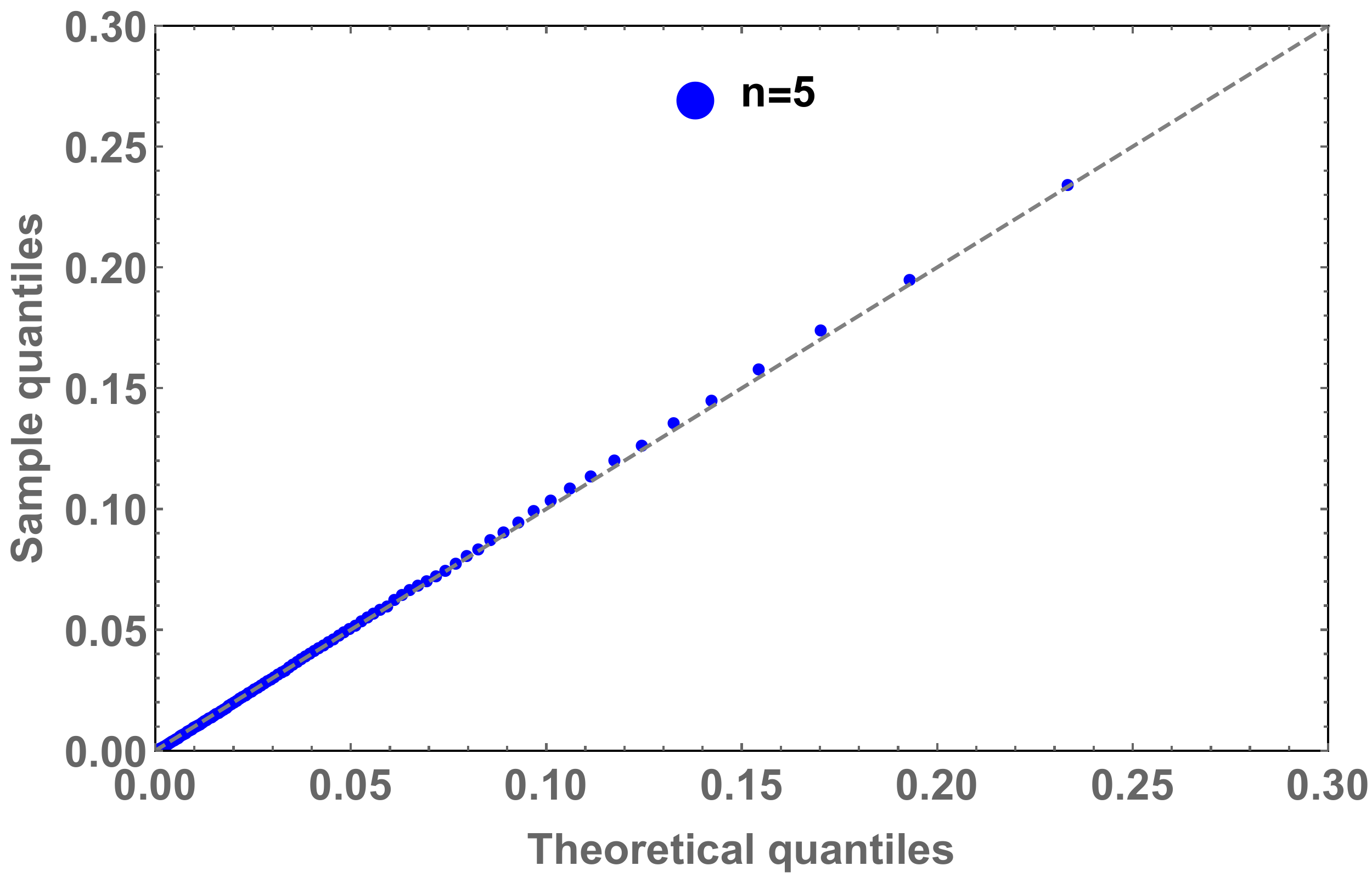}
    \end{subfigure}
    
    \begin{subfigure}[b]{.45\textwidth}
    \includegraphics[width=\textwidth]{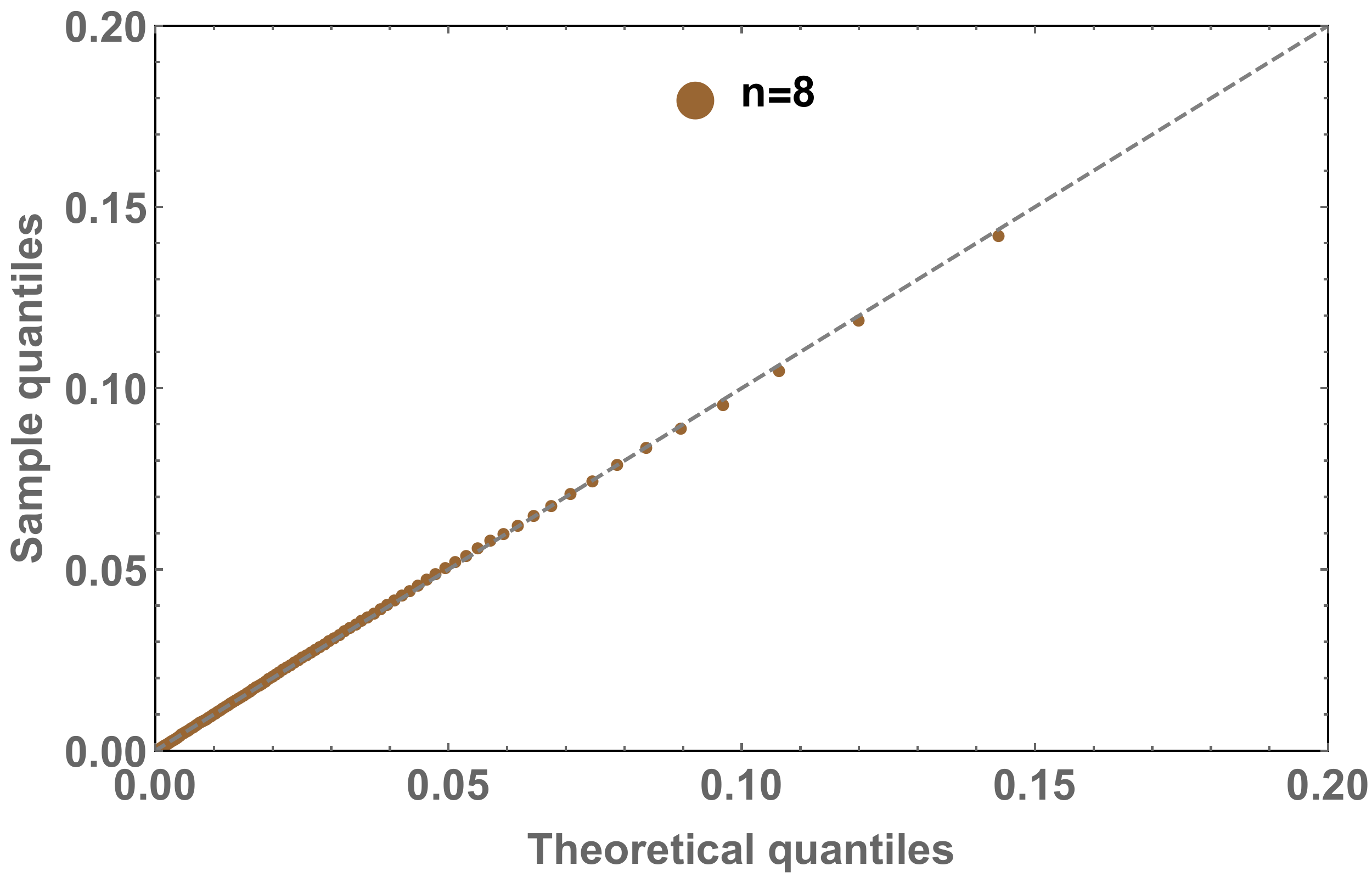}
    \end{subfigure}
    \begin{subfigure}[b]{.45\textwidth}
    \includegraphics[width=\textwidth]{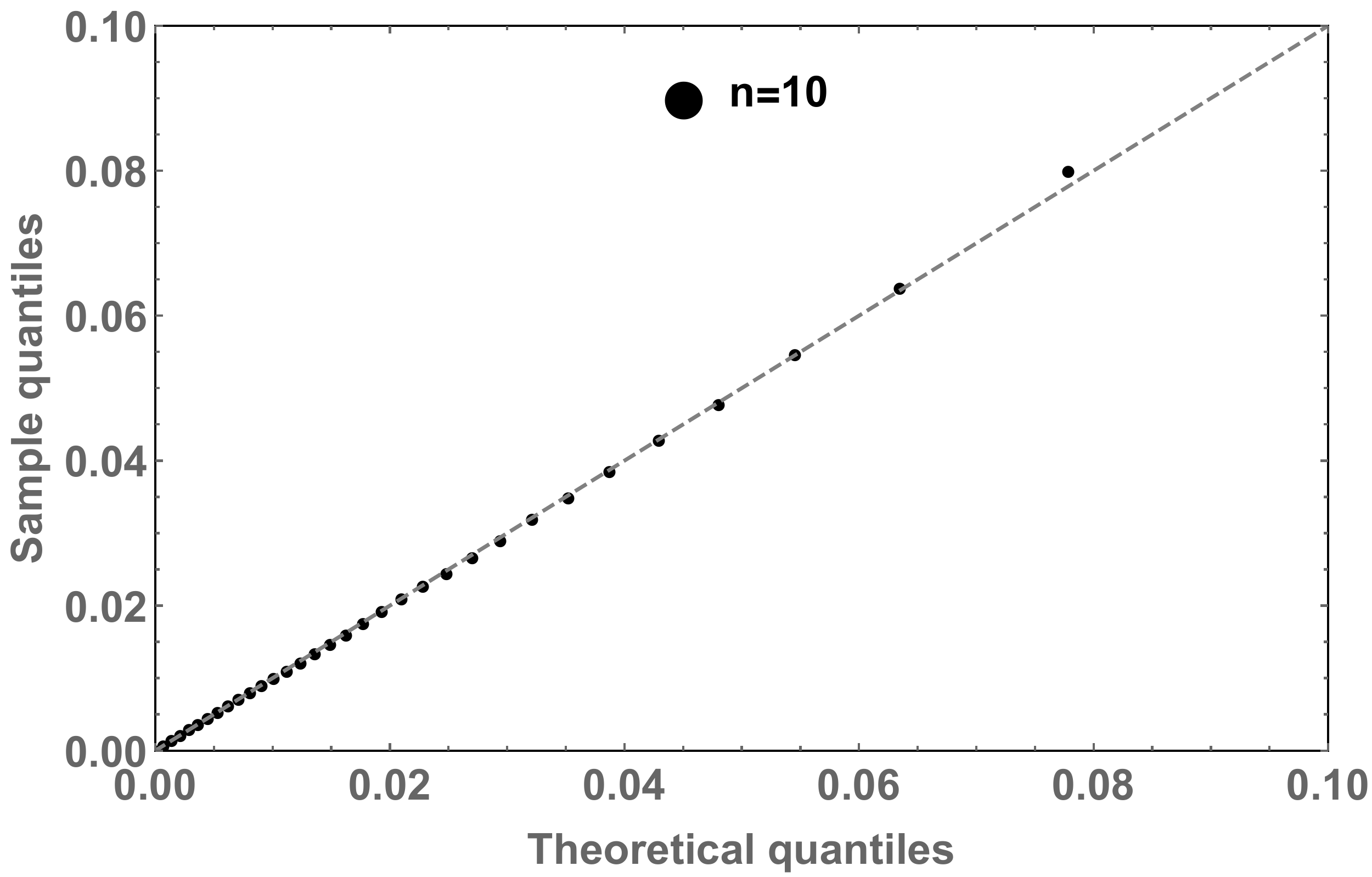}
    \end{subfigure}
    
    \caption{Quantile-Quantile plots of simulated data of $Z^{(n)}_1$ drawn from  $f^{(\alpha)}_{Z^{(n)}_1}(z)$ for different values of $n$ with $\alpha=2$ and $\theta=3$.}
     \label{fig_qq_mineigvec_different_n}
\end{figure}

\begin{figure}
    \centering
    \includegraphics[width=0.9\textwidth]{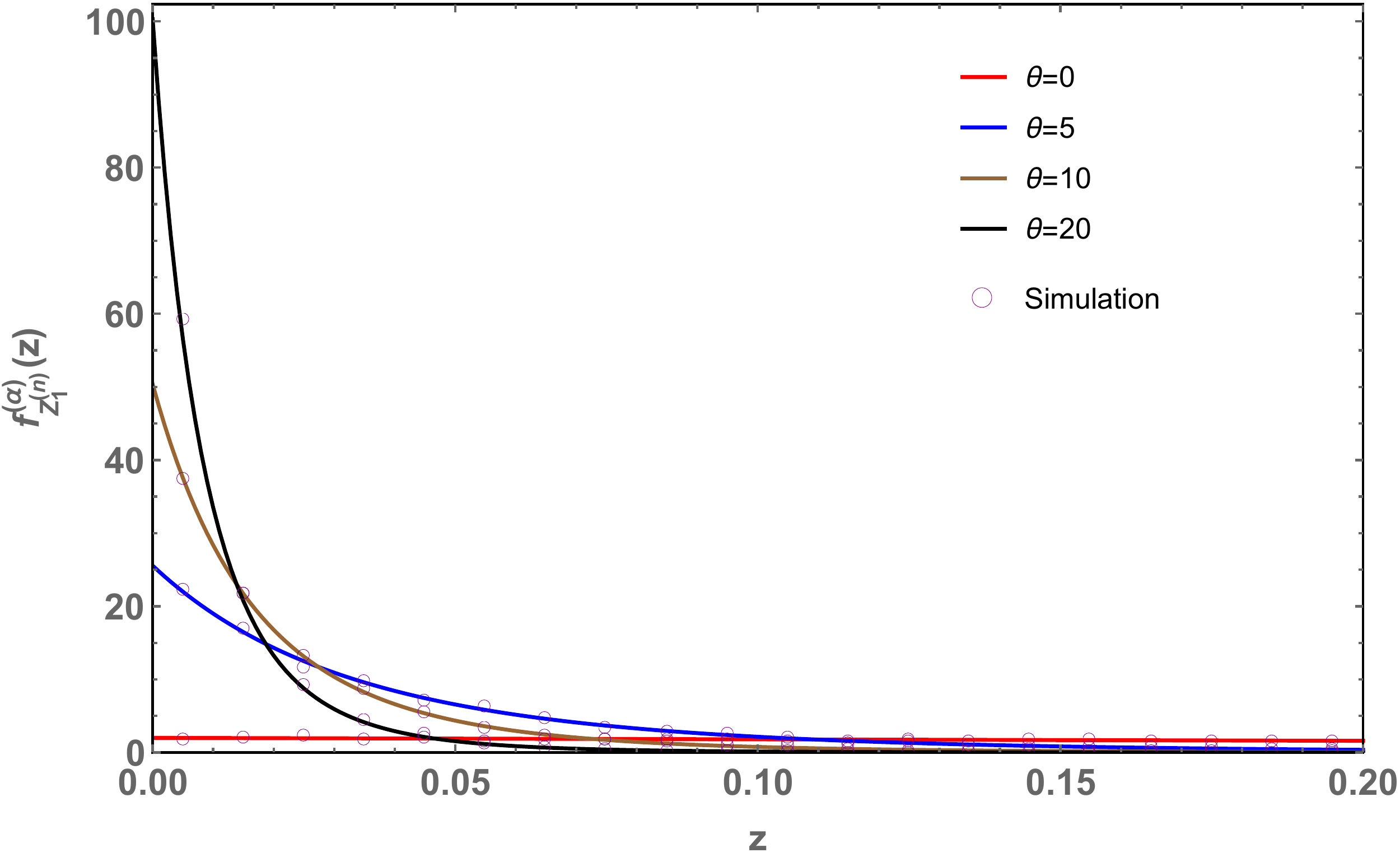}
    \caption{Comparison of simulated data points and the analytical p.d.f. $f^{(\alpha)}_{Z^{(n)}_1}(z)$ for different values of $\theta$ with $\alpha=2$ and $n=3$.}
    \label{fig_mineigvec_different_theta}
\end{figure}

\begin{figure}
    \centering
    \begin{subfigure}[b]{.45\textwidth}
    \includegraphics[width=\textwidth]{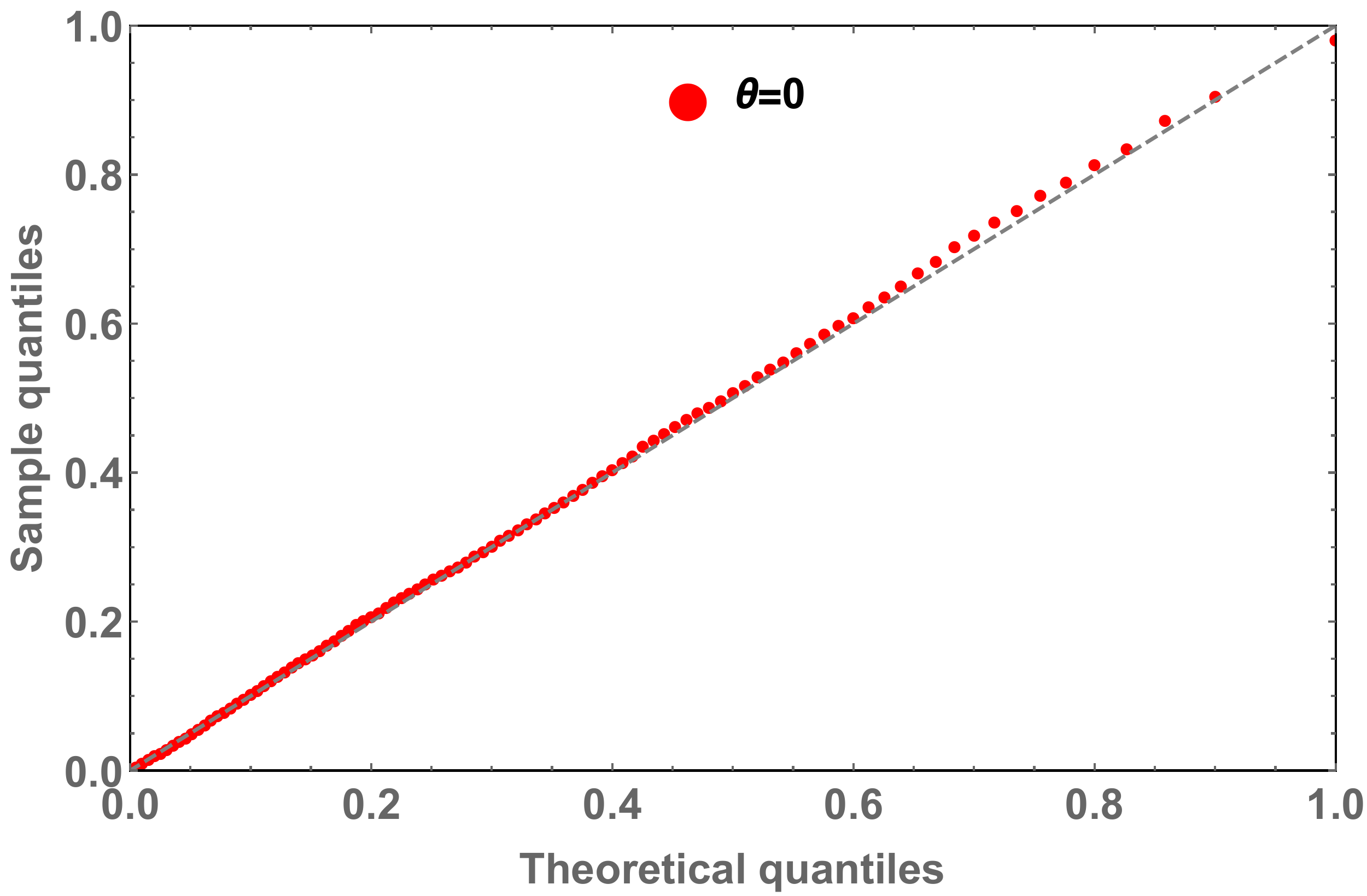}
    \end{subfigure}
    \begin{subfigure}[b]{.45\textwidth}
    \includegraphics[width=\textwidth]{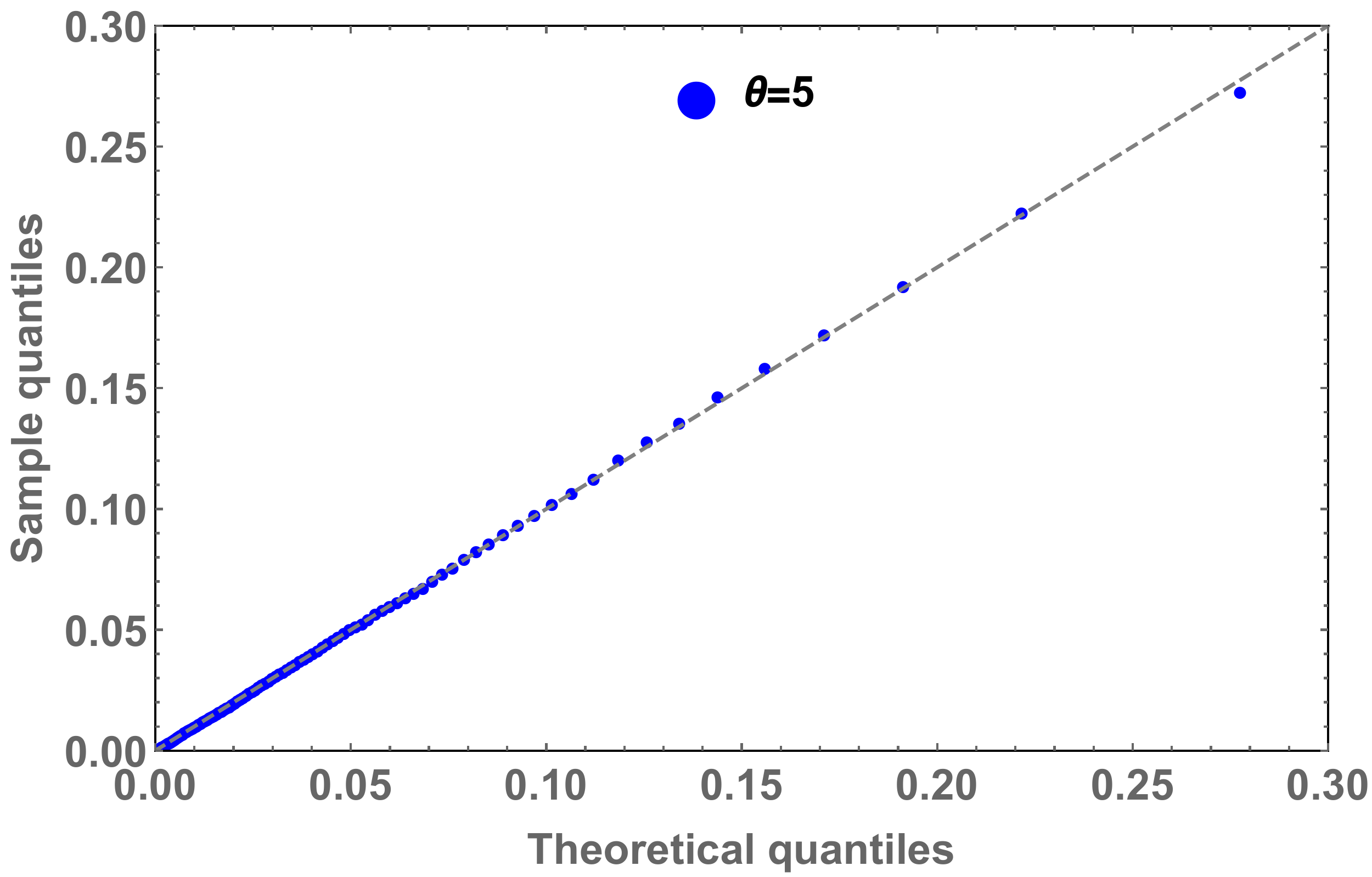}
    \end{subfigure}
    
    \begin{subfigure}[b]{.45\textwidth}
    \includegraphics[width=\textwidth]{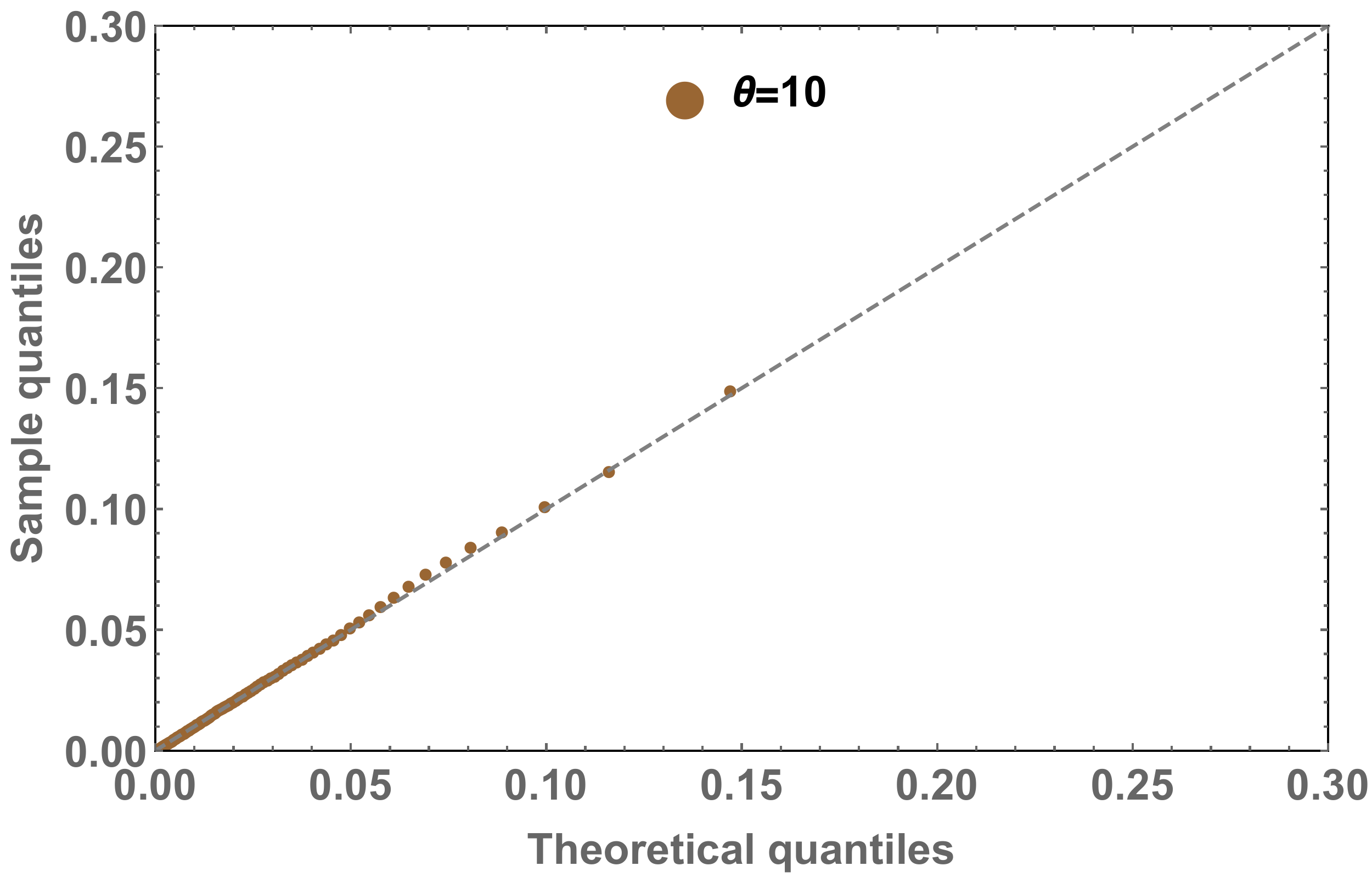}
    \end{subfigure}
    \begin{subfigure}[b]{.45\textwidth}
    \includegraphics[width=\textwidth]{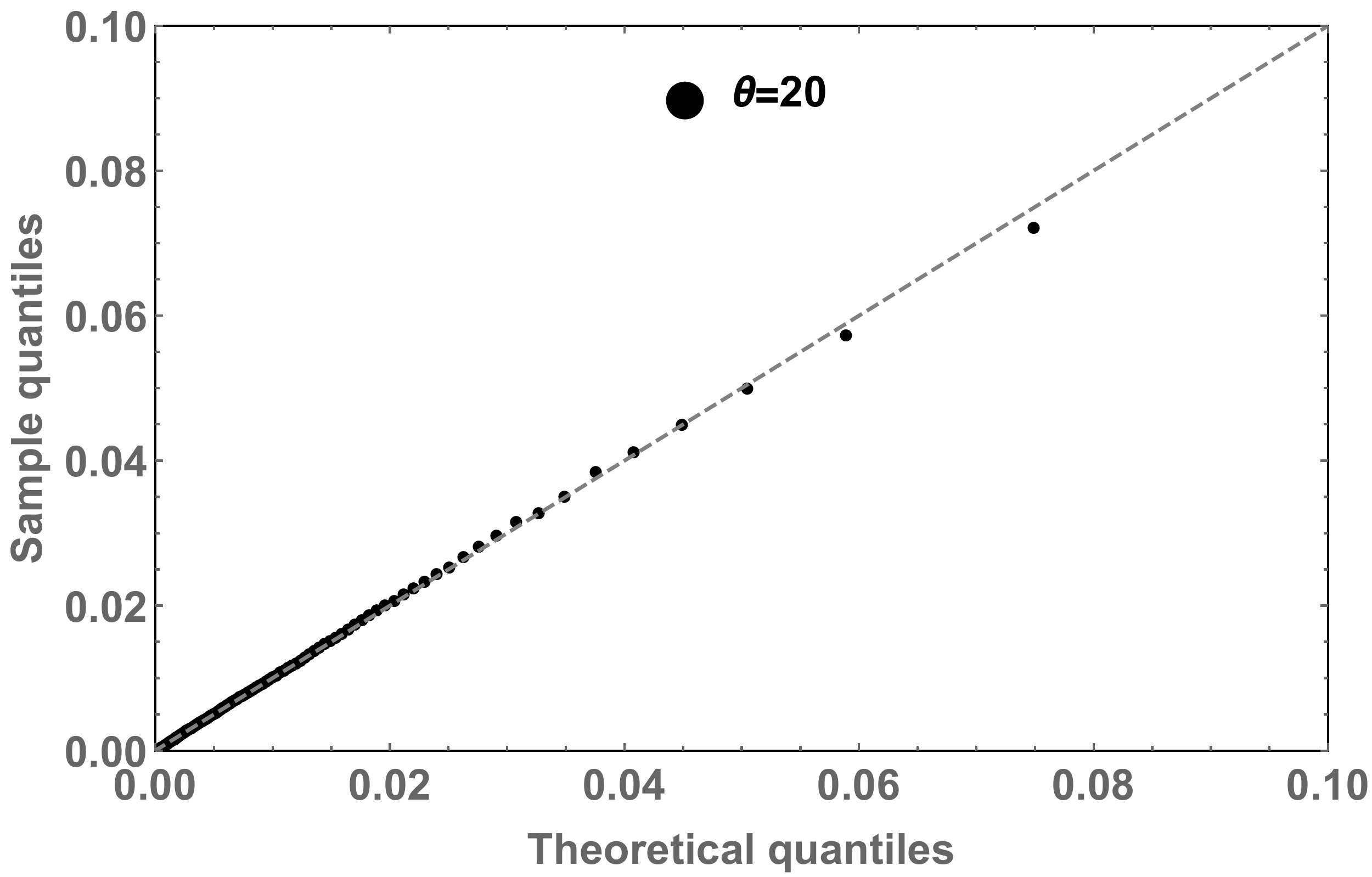}
    \end{subfigure}
    \caption{Quantile-Quantile plots of simulated data of $Z^{(n)}_1$ drawn from  $f^{(\alpha)}_{Z^{(n)}_1}(z)$ for different values of $\theta$ with $\alpha=2$ and $n=3$.}
    \label{fig_qq_mineigvec_different_theta}
\end{figure}

To further illustrate the utility of Theorem \ref{thm min vec}, we now focus on the asymptotic behavior of scaled $Z^{(n)}_1$ as $m$ and $n$ grow large, but their difference does not. The following corollary establishes the stochastic convergence result in this respect.

\begin{corollary}
\label{cor min asy}
As $m,n\to \infty$ such that $m-n$ is fixed\footnote{This is also known as the microscopic limit in the literature of theoretical physics \cite{ref:forresterLogGases,ref:ghur1}.} (i.e., $\alpha$ is fixed), the scaled random variable $nZ^{(n)}_1=n|\mathbf{v}^\dagger\mathbf{u}_1|^2$ {\it converges in distribution} to $\displaystyle \frac{\chi^2_2}{2(1+\theta)}$, where $\chi_2^2$ is a chi-squared random variable with two degrees of freedom (i.e., sum of squares of two independent standard normal random variables).
\end{corollary}
\begin{IEEEproof}
See Appendix \ref{appendix_corr2}.
\end{IEEEproof}
\begin{remark}
  It is noteworthy that, since the chi-squared random variable $V=\displaystyle \frac{\chi^2_2}{2(1+\theta)}$ has a continuous c.d.f., we have
  \begin{align}
      \lim_{n\to\infty} \Pr\left\{nZ_1\leq z\right\}=\Pr\left\{V\leq z\right\}=1-\exp\left(-(1+\theta)z\right)
  \end{align}
  in which the convergence is uniform in $z$ \cite[Lemma 2.11]{ref:vaart}.
\end{remark}
{\color{blue}This interesting result reveals that, in this particular asymptotic regime (i.e., $m,n\to\infty$ such that $m-n$ is fixed), irrespective of the value of $\theta$ (i.e., whether $\theta$ is either below or above the phase transition threshold $1$), the eigenvector corresponding to the least eigenvalue is informative with respect to the latent spike $\mathbf{v}$.} In contrast, as we are well aware of, when $m,n\to \infty$ such that $n/m\to1$, {\color{blue}even the eigenvector corresponding to the largest eigenvalue of a sample covariance matrix shows an  uninformative first order behavior with respect to a subcritical population spike (i.e., $\theta<1$) \cite{ref:paul, ref:bGeorges}. This observation brings into sharp focus the detection of subcritical spike (also known as a weak signal) with the eigenvalues of the sample covariance matrix. To be precise, below the phase transition, the maximum eigenvalue does not have discrimination power, since irrespective of the signal availability the maximum eigenvalue (after proper centering and scaling) stochastically converges to a Tracy-Widom distribution \cite{ref:baikPhaseTrans, ref:paul, ref:johnPaul}. To circumvent this ambiguity, informative tests based on {\it all the eigenvalues} instead of the largest eigenvalue have been proposed in \cite{ref:onatskiSphe,ref:onatskiSpiked,ref:dobriban} to detect subcritical spikes.
Although not directly related to Corollary \ref{cor min asy}, $\chi^2$ type stochastic convergence results for the spectral projectors of various matrix ensembles have been established in \cite{ref:bloemendal}, \cite{ref:bao}, \cite{ref:baoWang}, \cite{ref:capitaine}, \cite{ref:baoChi}.}

\begin{figure}
    \centering
    \includegraphics[width=0.9\textwidth]{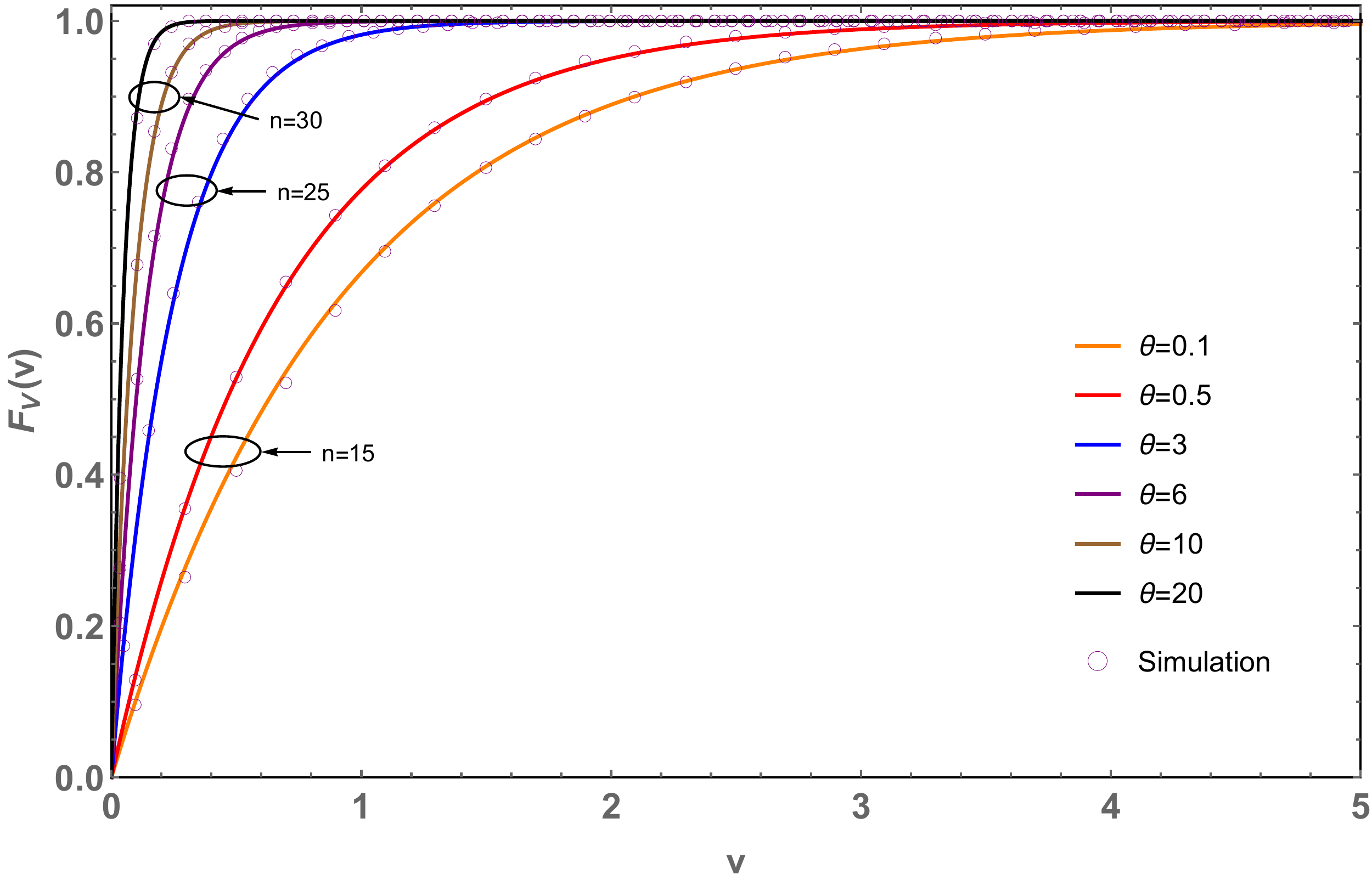}
    \caption{Comparison of simulated data points and the asymptotic analytical c.d.f. $F_{V}(v)$ for different values of $n, \theta$ with $\alpha=2$.}
    \label{min_asymp_cdf}
\end{figure}

The advantage of the asymptotic formula presented in Corollary \ref{cor min asy} is that it provides a simple expression which compares favorably with finite $n$ values. To further highlight this salient point, in Fig. \ref{min_asymp_cdf}, we compare the analytical asymptotic c.d.f. of $V$ with simulated data points corresponding to $\alpha=2$ with $n=15,25$ and $30$ for various values of $\theta$. The close agreement is clearly apparent from the figure.

 Having statistically characterized the eigenvector corresponding to the minimum eigenvalue, let us now focus on the eigenvector corresponding to the maximum eigenvalue (i.e., $\ell=n$).

\subsection{The P.D.F. of $Z^{(n)}_n$}
For $\ell=n$, (\ref{eq unitary comp}) specializes to 
\begin{align}
    \Psi_n(\boldsymbol{\lambda},s)=\int\limits_{\mathcal{U}_n} \text{etr}\left\{\mathbf{vv}^\dagger \mathbf{U}\left(\beta \boldsymbol{\Lambda}-s\mathbf{e}_n \mathbf{e}_n^T\right) \mathbf{U}^\dagger\right\} {\rm d} \mathbf{U},
\end{align}
which can be simplified using (\ref{eq split}) and   (\ref{eq cont integral}) to yield
\begin{align}
    \Psi_n(\boldsymbol{\lambda},s)=\frac{(n-1)!}{2\pi \mathrm{i}} 
        \oint\limits_{\mathcal{C}}
        \frac{e^\omega}{\displaystyle \left(s+\omega-\beta\lambda_n\right)\prod_{j=1}^{n-1} \left(\omega-\beta \lambda_j\right)} {\rm d}\omega.
\end{align}
This in turn enables us to express (\ref{eq mgf cont}) as
\begin{align}
    \mathcal{M}_{Z^{(n)}_n}(s)=C^\beta_{n,\alpha}
   \int\limits_{\mathcal{R}}
   \Delta_n^2(\boldsymbol{\lambda}) \prod_{j=1}^n \lambda_j^\alpha e^{-\lambda_j} \frac{1}{2\pi \mathrm{i}} 
        \oint\limits_{\mathcal{C}}
        \frac{e^\omega}{\displaystyle \left(s+\omega-\beta\lambda_n\right)\prod_{j=1}^{n-1} \left(\omega-\beta \lambda_j\right)} {\rm d}\omega\; {\rm d}\boldsymbol{\Lambda}.
\end{align}
 Now we take the inverse Laplace transform of both sides to yield
\begin{align}
    f^{(\alpha)}_{Z^{(n)}_n}(z)=C^\beta_{n,\alpha}
   \int\limits_{\mathcal{R}} e^{\beta \lambda_n z}
   \Delta_n^2(\boldsymbol{\lambda}) \prod_{j=1}^n \lambda_j^\alpha e^{-\lambda_j}
   \frac{1}{2\pi \mathrm{i}} 
        \oint\limits_{\mathcal{C}}
        \frac{e^{(1-z)\omega}}{\displaystyle \prod_{j=1}^{n-1} \left(\omega-\beta \lambda_j\right)} {\rm d}\omega\; {\rm d}\boldsymbol{\Lambda}
\end{align}
in which the innermost contour integral can be evaluated with the help of the residue theorem to obtain
\begin{align}
   f^{(\alpha)}_{Z^{(n)}_n}(z)=\frac{C^\beta_{n,\alpha}}{\beta^{n-2}}
   \int\limits_{\mathcal{R}} 
   e^{\beta \lambda_n z}
   \Delta_n^2(\boldsymbol{\lambda}) \prod_{j=1}^n \lambda_j^\alpha e^{-\lambda_j}
   \sum_{k=1}^{n-1}
   \frac{e^{(1-z)\beta \lambda_k}}{\displaystyle \prod_{\substack{i=1\\i\neq k}}^{n-1} \left(\lambda_k-\lambda_i\right)}\;
  {\rm d}\boldsymbol{\Lambda}.
\end{align}
Now let us rewrite the $n$-fold integral, keeping the integration with respect to $\lambda_n$ last, as
\begin{align}
\label{eq znpdf int}
  f^{(\alpha)}_{Z^{(n)}_n}(z)=\frac{C^\beta_{n,\alpha}}{\beta^{n-2}}  
  \int_{0}^\infty \lambda_n^\alpha e^{-(1-\beta z)\lambda_n}
  \Omega(\lambda_n,z) {\rm d}\lambda_n
\end{align}
where
\begin{align}
    \Omega(\lambda_n,z)=\int\limits_{\mathcal{I}} \sum_{k=1}^{n-1}
   \frac{e^{(1-z)\beta \lambda_k}}{\displaystyle \prod_{\substack{i=1\\i\neq k}}^{n-1} \left(\lambda_k-\lambda_i\right)}
   \Delta_{n-1}^2(\boldsymbol{\lambda})
   \prod_{j=1}^{n-1} \lambda_j^\alpha e^{-\lambda_j}
   \left(\lambda_n-\lambda_j\right)^2 {\rm d}\lambda_j
    \end{align}
 in which we have used the decomposition $\Delta_n^2(\boldsymbol{\lambda})=\prod_{j=1}^{n-1} (\lambda_n-\lambda_j)^2 \Delta_{n-1}^2(\boldsymbol{\lambda})$ and $\mathcal{I}=\{0<\lambda_1<\lambda_2<\ldots<\lambda_n\}$. For convenience, we relabel the variables as $\lambda_n=x$ and $x_j=\lambda_{j}$, $j=1,2,\ldots,n-1$, to arrive at  
\begin{align}
\label{eq pdf znx}
  f^{(\alpha)}_{Z^{(n)}_n}(z)=\frac{C^\beta_{n,\alpha}}{\beta^{n-2}}  
  \int_{0}^\infty x^\alpha e^{-(1-\beta z)x}
  \Omega(x,z) {\rm d}x
\end{align}
where
\begin{align}
    \Omega(x,z)=\int\limits_{\mathcal{I}_x} \sum_{k=1}^{n-1}
   \frac{e^{(1-z)\beta x_k}}{\displaystyle \prod_{\substack{i=1\\i\neq k}}^{n-1} \left(x_k-x_i\right)}
   \Delta_{n-1}^2(\boldsymbol{x})
   \prod_{j=1}^{n-1} x_j^\alpha e^{-x_j}
   \left(x_j-x\right)^2 {\rm d}x_j
    \end{align}
with $\mathcal{I}_x=\{0<x_1<x_2<\ldots<x_{n-1}<x\}$. Since the integrand in the above $(n-1)$-fold integral is symmetric in $x_1,x_2,\ldots,x_{n-1}$, we may remove the ordered region of integration to obtain
\begin{align}
    \Omega(x,z)=\frac{1}{(n-1)!}\int\limits_{(0,x)^{n-1}} \sum_{k=1}^{n-1}
   \frac{e^{(1-z)\beta x_k}}{\displaystyle \prod_{\substack{i=1\\i\neq k}}^{n-1} \left(x_k-x_i\right)}
   \Delta_{n-1}^2(\boldsymbol{x})
   \prod_{j=1}^{n-1} x_j^\alpha e^{-x_j}
   \left(x_j-x\right)^2 {\rm d}x_j.
    \end{align}
Consequently, we can observe that the each term in the above summation evaluates to the same amount. Therefore, capitalizing on that observation, we may simplify the above $(n-1)$-fold integral to yield

    \begin{align}
    \Omega(x,z)=\frac{1}{(n-2)!}\int\limits_{(0,x)^{n-1}}
   \frac{e^{(1-z)\beta x_1}}{\displaystyle \prod_{{i=2}}^{n-1} \left(x_1-x_i\right)}
   \Delta_{n-1}^2(\boldsymbol{x})
   \prod_{j=1}^{n-1} x_j^\alpha e^{-x_j}
   \left(x_j-x\right)^2 {\rm d}x_j.
    \end{align}
Now it is convenient to introduce the variable transformations $xt_j= x_j$, $j=1,2,\ldots,n-1$, to  the above integral to arrive at
\begin{align}
    \Omega(x,z)=\frac{x^{(n-1)(n+\alpha)+1}}{(n-2)!}\int\limits_{(0,1)^{n-1}}
    \frac{e^{(1-z)\beta x t_1}}{\displaystyle \prod_{{i=2}}^{n-1} \left(t_1-t_i\right)}
    \Delta_{n-1}^2(\boldsymbol{t})
    \prod_{j=1}^{n-1} t_j^\alpha (1-t_j)^2 e^{-xt_j}
    {\rm d}t_j.
\end{align}
To facilitate further analysis, let us relabel the variables as, $t=t_1$, $y_{j-1}=t_j$, $j=2,3,\ldots,n-1$, and 
keep the integration with respect to $t$ last giving
\begin{align}
\label{Phi int formn}
    \Omega(x,z)=\frac{x^{(n-1)(n+\alpha)+1}}{(n-2)!}\int_0^1
    e^{-(1-(1-z)\beta)xt} (1-t)^2 t^\alpha
    \mathcal{P}_{n-2}(t,x) {\rm d}t
\end{align}
where
\begin{align}
\label{intP}
    \mathcal{P}_{n-2}(t,x)=
    \int\limits_{(0,1)^{n-2}}
    \Delta_{n-2}^2(\boldsymbol{y})
    \prod_{j=1}^{n-2} y_j^\alpha (1-y_j)^2 (t-y_j)e^{-xy_j}
    {\rm d}y_j
\end{align}
and we have used the decomposition $\Delta_{n-1}^2(\boldsymbol{t})=\prod_{j=1}^{n-1}(t-y_j)^2 \Delta_{n-2}^2(\boldsymbol{y})$.
The $(n-2)$-fold integral in (\ref{intP}) can be evaluated with the help of \cite[Corollary 1]{ref:chiani} and \cite[Eq. 6.5.1]{ref:erdelyi} to yield the Hankel determinant
\begin{align}
\label{pdetans}
   \mathcal{P}_{n-2}(t,x)=(n-2)! \det\left[t \mathcal{A}_{i,j}^{(\alpha)}(x)-\mathcal{A}_{i,j}^{(\alpha+1)}(x)\right]_{i,j=1,2,\ldots,(n-2)} 
\end{align}
where
\begin{align}
\label{pdetarg}
    \mathcal{A}^{(\alpha)}_{i,j}(x)=\mathcal{B}(3,i+j+\alpha-1) {}_1F_1\left(\alpha+i+j-1;\alpha+i+j+2;-x\right), 
\end{align}
$\mathcal{B}(p,q)=(p-1)!(q-1)!/(p+q-1)!$, and ${}_1F_1 (a;c;z)$ is the confluent hypergeometric function of the first kind\footnote{This function assumes a finite series expansion involving exponentials and powers of $x$, since $\alpha,i$, and $j$ take non-negative integer values in (\ref{pdetarg}).}. Finally, we use (\ref{pdetarg}) and (\ref{pdetans}) in (\ref{Phi int formn}) and substitute the resultant integral expression into (\ref{eq pdf znx}) with some algebraic manipulation to arrive at the exact p.d.f. of $Z^{(n)}_n$ as given in the following theorem.
\begin{theorem}
\label{pdf max vec}
    Let $\mathbf{W}\sim\mathcal{CW}_n\left(m,\mathbf{I}_n+\theta \mathbf{v}\mathbf{v}^\dagger\right)$ with $||\mathbf{v}||=1$ and $\theta>0 $. Let $\mathbf{u}_n$ be the eigenvector corresponding to the largest eigenvalue of $\mathbf{W}$. Then the p.d.f. of $Z^{(n)}_n=|\mathbf{v}^\dagger \mathbf{u}_n|^2\in(0,1)$ is given by
    \begin{align}
    \label{pdfmaxfinal}
        f_{Z^{(n)}_n}^{(\alpha)}(z)=\frac{C^\beta_{n,\alpha}}{\beta^{n-2}} \int_0^\infty
        x^{n^2+n\alpha-n+1} e^{-\left(1-\beta z\right)x}
        \mathcal{J}_{n}(x,z) {\rm d}x
    \end{align}
    where
    \begin{align}
    \label{pdfmaxfinalaux}
        \mathcal{J}_{n}(x,z)=\int_0^1
        e^{-(1-(1-z)\beta)xt} t^\alpha (1-t)^2 
    \det\left[t \mathcal{A}_{i,j}^{(\alpha)}(x)-\mathcal{A}_{i,j}^{(\alpha+1)}(x)\right]_{i,j=1,2,\ldots,(n-2)}  {\rm d}t
    \end{align}
    and $\beta=\theta/(1+\theta)$.
\end{theorem}
It is noteworthy that the above p.d.f. is valid only for $\theta>0$, since there exist a singularity of order $n-2$ for $\theta=0$. Nevertheless, the direct substitution of $\theta=0$ is possible in the case of $n=2$.
Although further simplification of (\ref{pdfmaxfinal}) and (\ref{pdfmaxfinalaux}) seems intractable for general matrix dimensions $m$ and $n$ (i.e., $n$ and $\alpha$), as shown in the following corollary, closed-form solutions are possible for the important configurations of $n=2,3$, and $n=4$ for arbitrary $\alpha$.
\begin{corollary}
\label{corpdfmaxf}
The exact p.d.f. of $Z^{(n)}_n\in(0,1)$ corresponding to $n=2,3$, and $n=4$ are given, respectively, by
\begin{align}
    f_{Z^{(2)}_2}^{(\alpha)}(z)&=\frac{2(2\alpha+3)! (1-\beta)^{\alpha+2}}{(\alpha+1)! (\alpha+3)! (2-\beta)^{2\alpha+4}}
    {}_2F_1\left(3,2\alpha+4;\alpha+4;\frac{1-\beta(1-z)}{2-\beta}\right)\\
    f_{Z^{(3)}_3}^{(\alpha)}(z)&=\frac{4 (3\alpha+7)!(1-\beta)^{\alpha+3}}{(\alpha+2)!(\alpha+3)! (\alpha+4)! \beta (3-\beta)^{3\alpha+8}}\left(\mathcal{F}^{(\alpha,\beta)}_2(4,5,z)-\mathcal{F}^{(\alpha,\beta)}_2(5,4,z)\right)\\
    f_{Z^{(4)}_4}^{(\alpha)}(z)&=\frac{(1-\beta)^{\alpha+4}}{2 \alpha ! (\alpha+1)! (\alpha+2)!(\alpha+4)!\beta^2}
    \sum_{k=0}^2
    \left[\mathcal{G}^{(\alpha,\beta)}_k(a_k,b_k,z)-\mathcal{G}^{(\alpha,\beta)}_k(c_k,d_k,z)\right.\nonumber\\
    & \hspace{7cm}-2\mathcal{G}^{(\alpha,\beta)}_{k+1}(a_k,b_k,z)+2\mathcal{G}^{(\alpha,\beta)}_{k+1}(c_k,d_k,z)\nonumber\\
    & \hspace{7.5cm}\left.+\mathcal{G}^{(\alpha,\beta)}_{k+2}(a_k,b_k,z)-\mathcal{G}^{(\alpha,\beta)}_{k+2}(c_k,d_k,z)\right]\label{pdfmaxn3}
\end{align}
where $(a_0\; a_1\; a_2)\equiv(5\; 5\; 4)$, $(b_0\; b_1\; b_2)\equiv(7\; 6\; 6)$, $(c_0\; c_1\; c_2)\equiv(6\; 4\; 5)$, $(d_0\; d_1\; d_2)\equiv(6\; 7\; 5)$,
$\mathcal{F}_2^{(\alpha,\beta)}(a,b,z)=
F_2\left(3\alpha+8,3,3,\alpha+a,\alpha+b;\frac{1}{3-\beta}, \frac{1-\beta(1-z)}{3-\beta}\right)$,
\begin{align}
    \mathcal{G}^{(\alpha,\beta)}_N(b,c,z)=
    \mathcal{B}(3,\alpha+b-3)&\mathcal{B}(3,\alpha+c-3)\frac{ (\alpha+N)!}{\left(1-\beta(1-z)\right)^{\alpha+N+1}}\nonumber\\
    & \quad\times \left[
    \frac{(3\alpha+12-N)!}{(3-\beta z)^{3\alpha+13-N}}
    \mathcal{\widetilde{F}}_2^{(\alpha, \beta)}
    \left(3\alpha+13-N,b,c,\frac{1}{3-\beta z}\right)\right.\nonumber\\
    & \qquad \quad  \Biggl.-\sum_{k=0}^{\alpha+N}
    \frac{(3\alpha+12+k-N)!}{k! (4-\beta)^{3\alpha+13-N+k}}
    (1-\beta(1-z))^k\nonumber\\
    & \qquad \qquad \qquad \times 
    \mathcal{\widetilde{F}}_2^{(\alpha, \beta)}
    \left(3\alpha+13-N+k,b,c,\frac{1}{4-\beta }\right)
    \Biggr],
\end{align}
and $\mathcal{\widetilde{F}}_2^{(\alpha, \beta)}
    \left(a,b,c,x\right)=F_2(a,3,3,\alpha+b,\alpha+c;x,x)$ with $F_2(a,b,c,d,g;x,y)$ denoting the hypergeometric function of two arguments \cite[Eq. 5.7.1.7]{ref:erdelyi}.
\end{corollary}

The above formulas follow by simplifying the determinant in (\ref{pdfmaxfinalaux}) 
and subsequent integration using \cite[Eqs. 7.621.4 and 7.622.3]{ref:gradshteyn} with some algebraic manipulation. However, for $n=4$ case, to keep the final answer tractable, after simplifying the determinant in (\ref{pdfmaxfinalaux}), we use integration by parts to evaluate the integral with respect to $t$. Subsequently, we employ \cite[Eq. 7.622.3]{ref:gradshteyn} to evaluate the integral in (\ref{pdfmaxfinal}) with some tedious algebraic manipulation to obtain (\ref{pdfmaxn3}).

{\color{blue}
It turns out that the analytical formula given in Theorem 4 can be used to show that, for $n=2$,  $f_{Z^{(2)}_2}^{(\alpha)}(z)$ is convex in $z$. To demonstrate this, following Theorem 4, let us rewrite the p.d.f. corresponding to $n=2$ as
\begin{align}
   f_{Z^{(2)}_2}^{(\alpha)}(z)=C^\beta_{2,\alpha}\int_0^\infty \int_0^1 e^{z(1-t)\beta x}
        \rho(x,t) {\rm d}t {\rm d}x
    \end{align}
    where
    \begin{align}
        \rho(x,t)= e^{-[1+(1-\beta)t] x} x^{2\alpha+1} t^\alpha (1-t)^2.
\end{align}
Since for each fixed $t,x$, $e^{z(1-t)\beta x}$ is convex in $z$ and $\rho(x,t)\geq 0$, we conclude that $f_{Z^{(n)}_n}^{(\alpha)}(z)$ is convex in $z$. This observation does not hold for $n\geq 3$, since the determinant term in the integrand changes it sign depending on the values of $x$ and $t$ . However, keeping in mind that
\begin{align}
    f_{Z^{(n)}_n}^{(\alpha)}(z)\leq \frac{C^\beta_{n,\alpha}}{\beta^{n-2}} \int_0^\infty
        x^{n^2+n\alpha-n+1} e^{-\left(1-\beta z\right)x}
        {J}_{n}(x,z) {\rm d}x
    \end{align}
    where
    \begin{align}
        {J}_{n}(x,z)=\int_0^1
        e^{-(1-(1-z)\beta)xt} t^\alpha (1-t)^2 
    \left|\det\left[t \mathcal{A}_{i,j}^{(\alpha)}(x)-\mathcal{A}_{i,j}^{(\alpha+1)}(x)\right]_{i,j=1,2,\ldots,(n-2)}\right|  {\rm d}t,
\end{align}
we may follow similar arguments as before to establish the fact that there exists a {\it convex function}  $h_n(z)$ such that
\begin{align}
    f_{Z^{(n)}_n}^{(\alpha)}(z)\leq h_n(z),\; n=2,3,\ldots,
\end{align}
for which the equality holds for $n=2$.}

\begin{remark}
  It is noteworthy that although the hypergeometric functions ${}_2F_1 (\cdot)$ and $F_2(\cdot)$ assume infinite series expansions, for the parameters of our interest here, they can be expressed as finite sums of rational functions of $z$ with the help of \cite[Eq. 2.1.4.23]{ref:erdelyi} and  \cite[Eq. 5.8.1.2]{ref:erdelyi}. However, those formulas are not presented here to keep the representations clear and concise.
\end{remark}
{\color{blue}A careful inspection of (\ref{pdfmaxfinal}) and (\ref{pdfmaxfinalaux}) reveals that the determinant of square matrix, whose size depends on $n$, prevents us from applying the same asymptotic framework that we have developed to characterize the asymptotic behavior of  $Z^{(n)}_1$ (i.e., as $m,n\to\infty$ such that $m-n$ is fixed). This stems from the fact that, as we are well aware of, in this asymptotic regime, $Z_n^{(n)}$ undergoes a phase transition at $\theta=1$ (i.e., above it is $O(1)$, whereas below it is $o(1)$).  To circumvent this difficulty, it is natural to consider other asymptotic regimes\footnote{The most prevalent asymptotic domain assumes $m,n\to \infty$ such that $n/m\to \gamma\in(0,1)$.} in view of obtaining relatively simple expressions for the p.d.f. of $Z^{(n)}_n$. 
In this respect, various stochastic convergence results for $Z^{(n)}_n$ have been established in \cite{ref:paul,ref:bGeorges,ref:bloemendal,ref:haokai,ref:bao,ref:wWang, ref:sandeep,ref:bai,ref:baoChi,ref:oliver,ref:wolf}.} Nevertheless, a finite dimensional analysis has not been available in the literature to date.  
 
\begin{figure}
    \centering
    \includegraphics[width=0.9\textwidth]{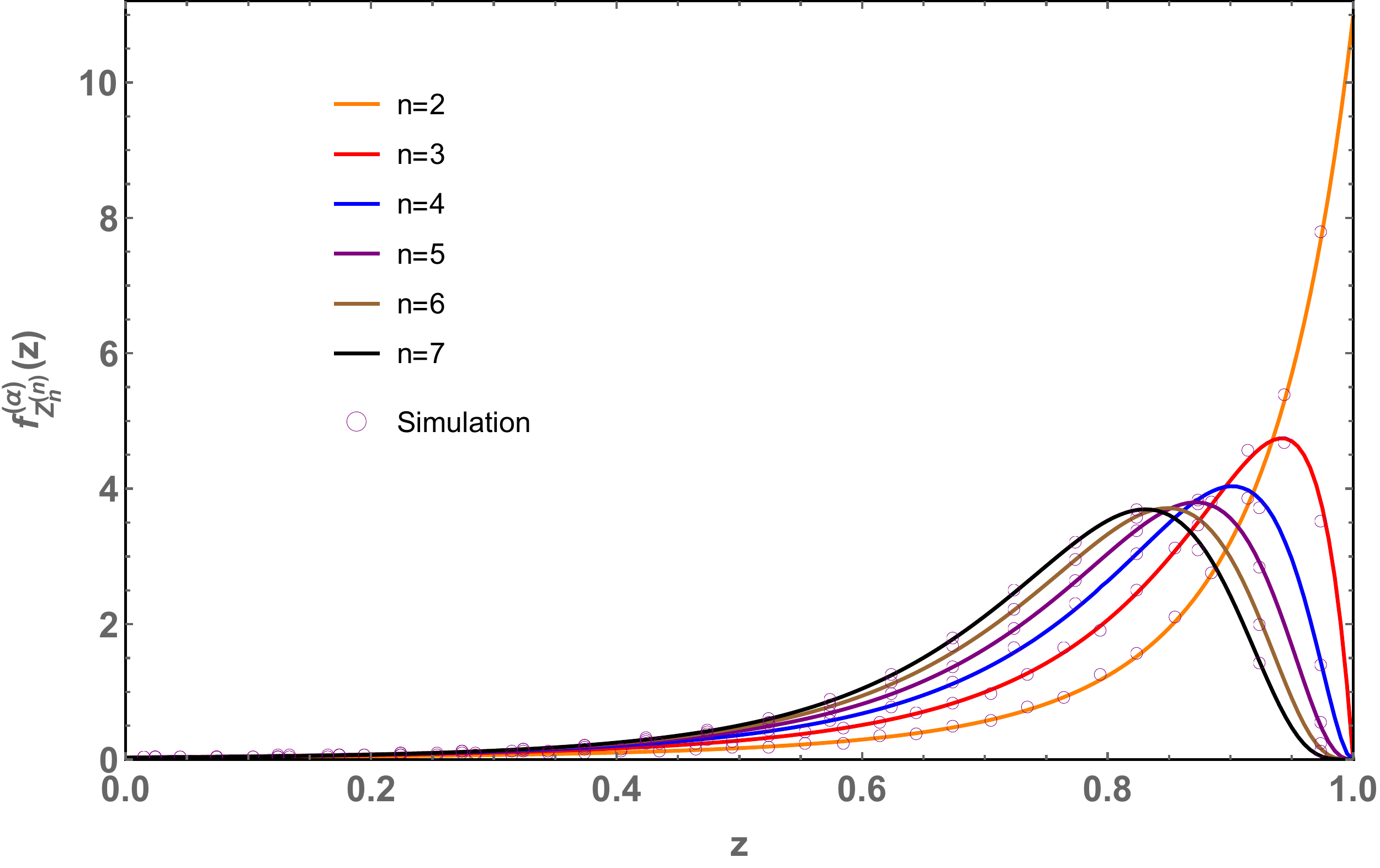}
    \caption{Comparison of simulated data points and the analytical p.d.f. $f^{(\alpha)}_{Z^{(n)}_n}(z)$ for different values of $n$ with $\alpha=2$ and $\theta=3$.}
    \label{fig_maxeigvec_different_n}
\end{figure}

\begin{figure}
    \centering
    \begin{subfigure}[b]{.45\textwidth}
    \includegraphics[width=\textwidth]{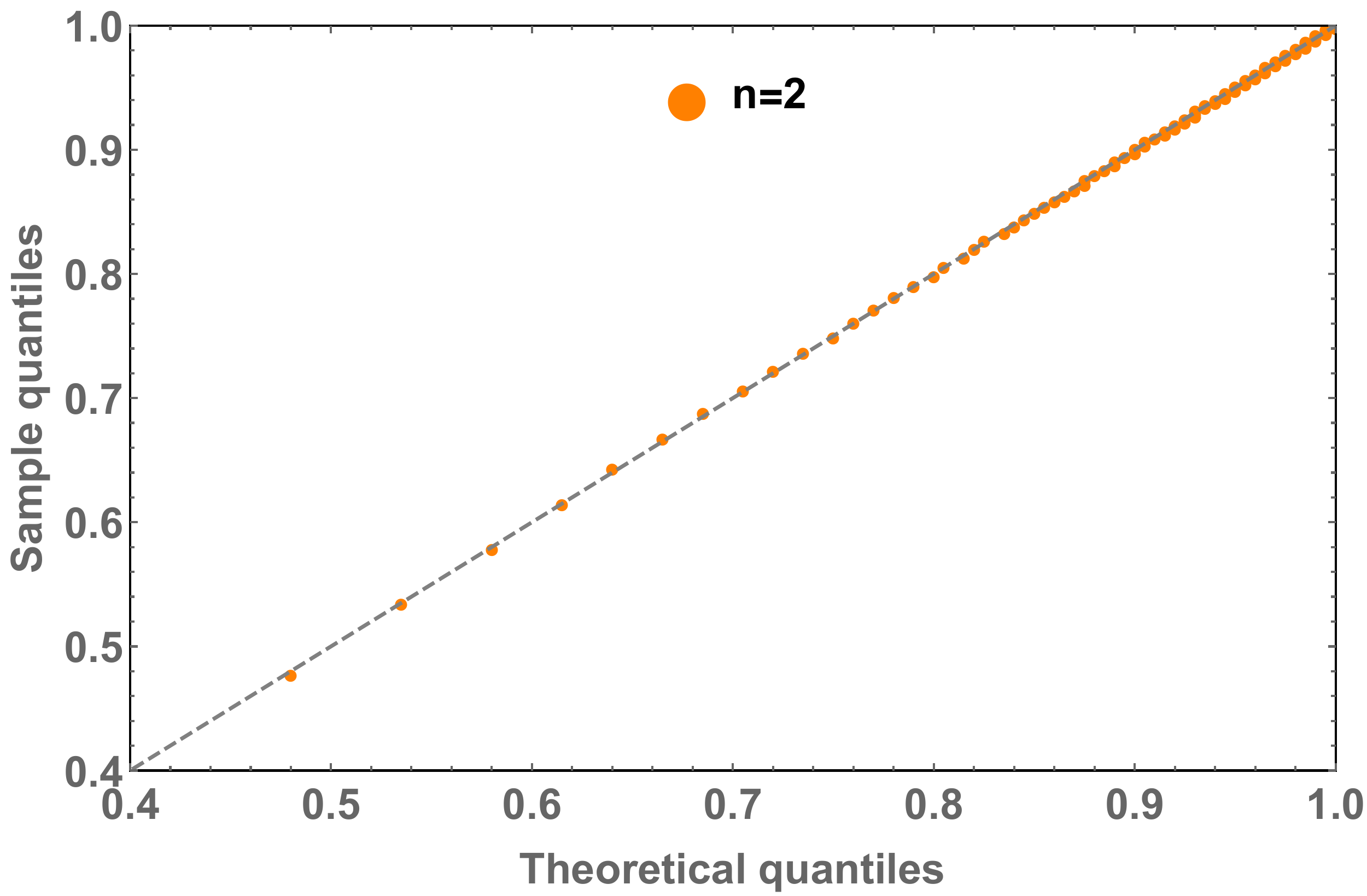}
    \end{subfigure}
    \begin{subfigure}[b]{.45\textwidth}
    \includegraphics[width=\textwidth]{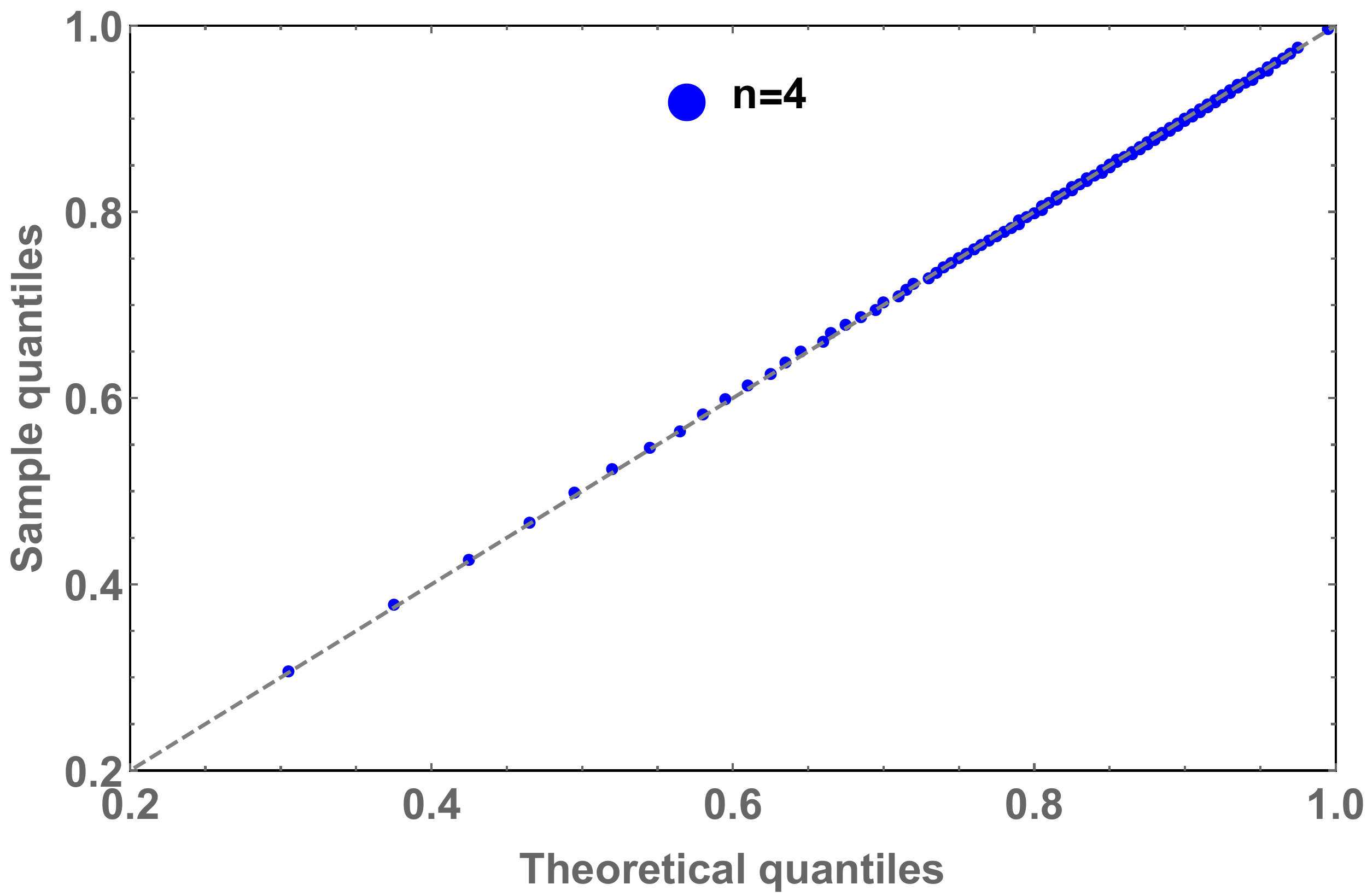}
    \end{subfigure}
    
    \begin{subfigure}[b]{.45\textwidth}
    \includegraphics[width=\textwidth]{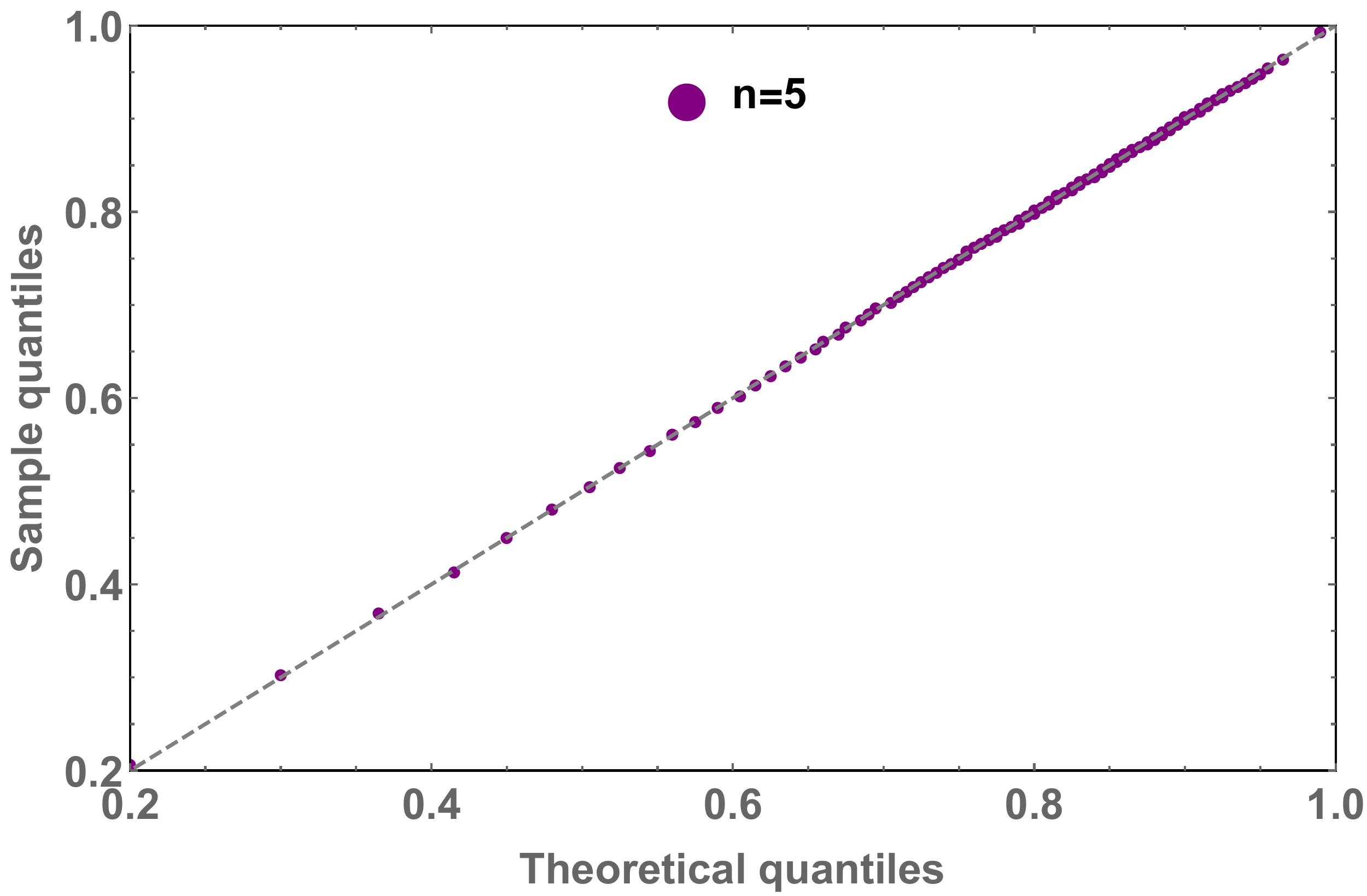}
    \end{subfigure}
    \begin{subfigure}[b]{.45\textwidth}
    \includegraphics[width=\textwidth]{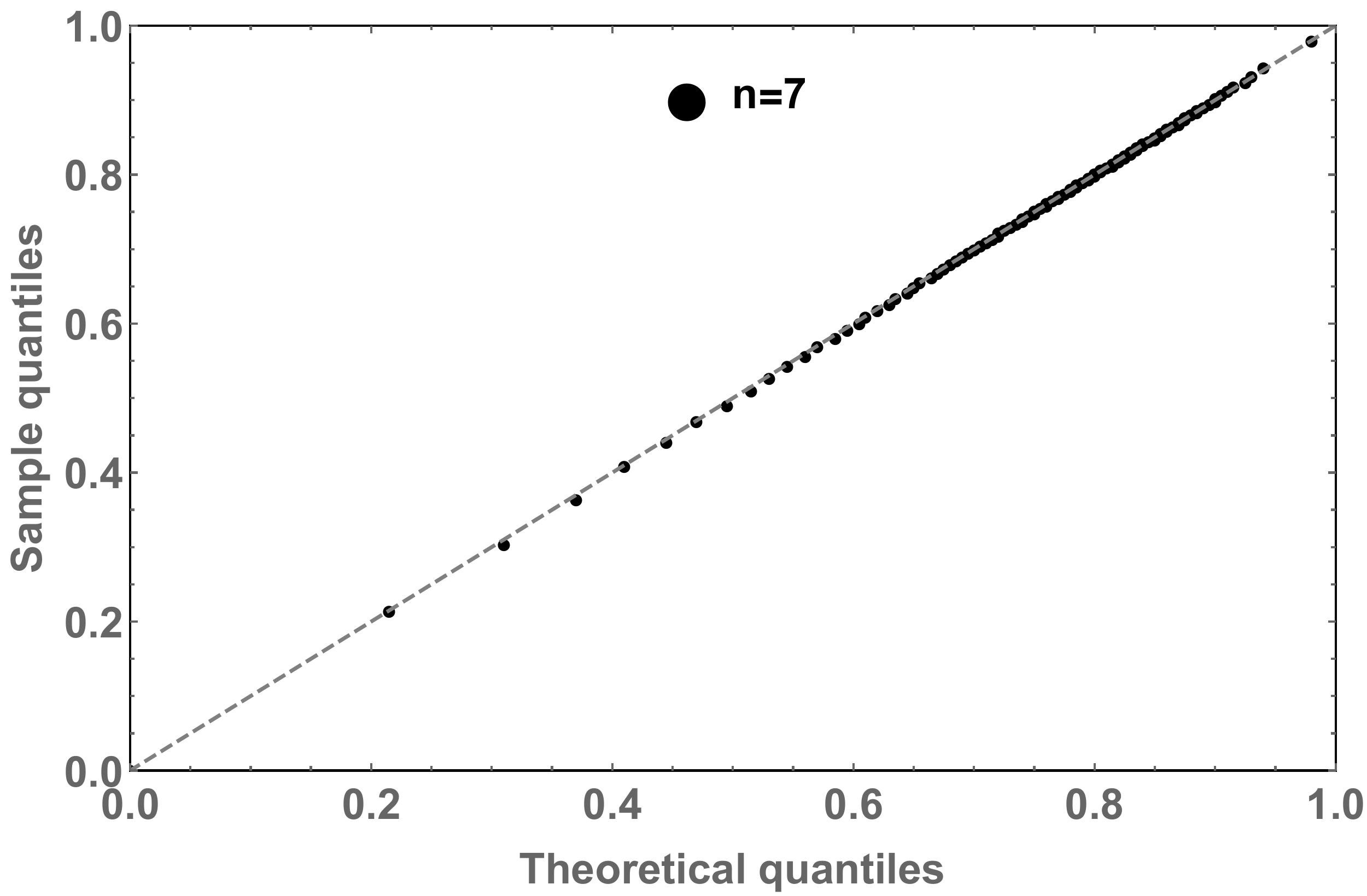}
    \end{subfigure}
    \caption{Quantile-Quantile plots of simulated data of $Z^{(n)}_n$ drawn from  $f^{(\alpha)}_{Z^{(n)}_n}(z)$ for different values of $n$ with $\alpha=2$ and $\theta=3$.}
    \label{fig_qq_maxeigvec_different_n}
\end{figure}

\begin{figure}
    \centering
    \includegraphics[width=0.9\textwidth]{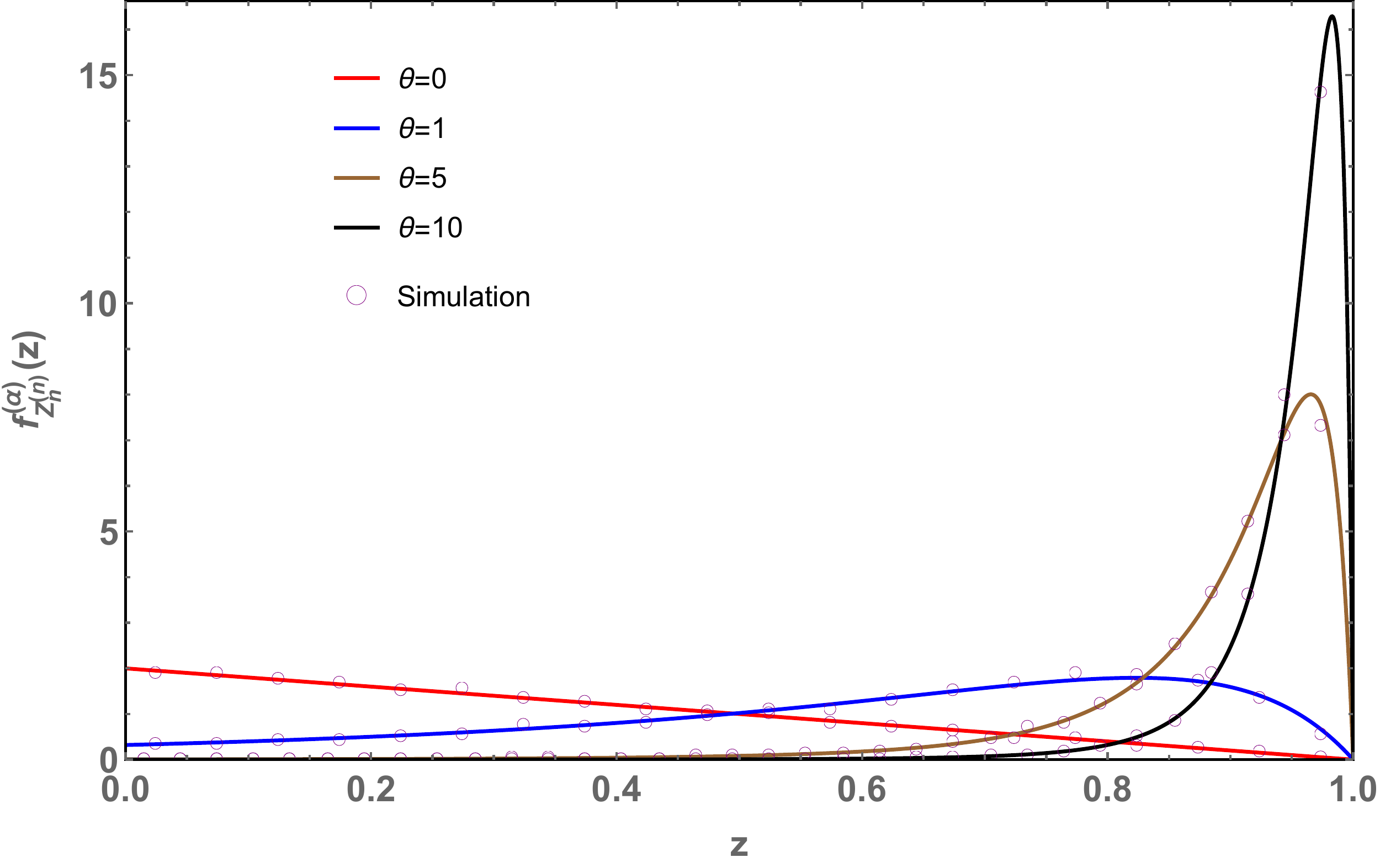}
    \caption{Comparison of simulated data points and the analytical p.d.f. $f^{(\alpha)}_{Z^{(n)}_n}(z)$ for different values of $\theta$ with $\alpha=2$ and $n=3$.}
    \label{fig_maxeigvec_different_theta}
\end{figure}

\begin{figure}
    \centering
    \begin{subfigure}[b]{.45\textwidth}
    \includegraphics[width=\textwidth]{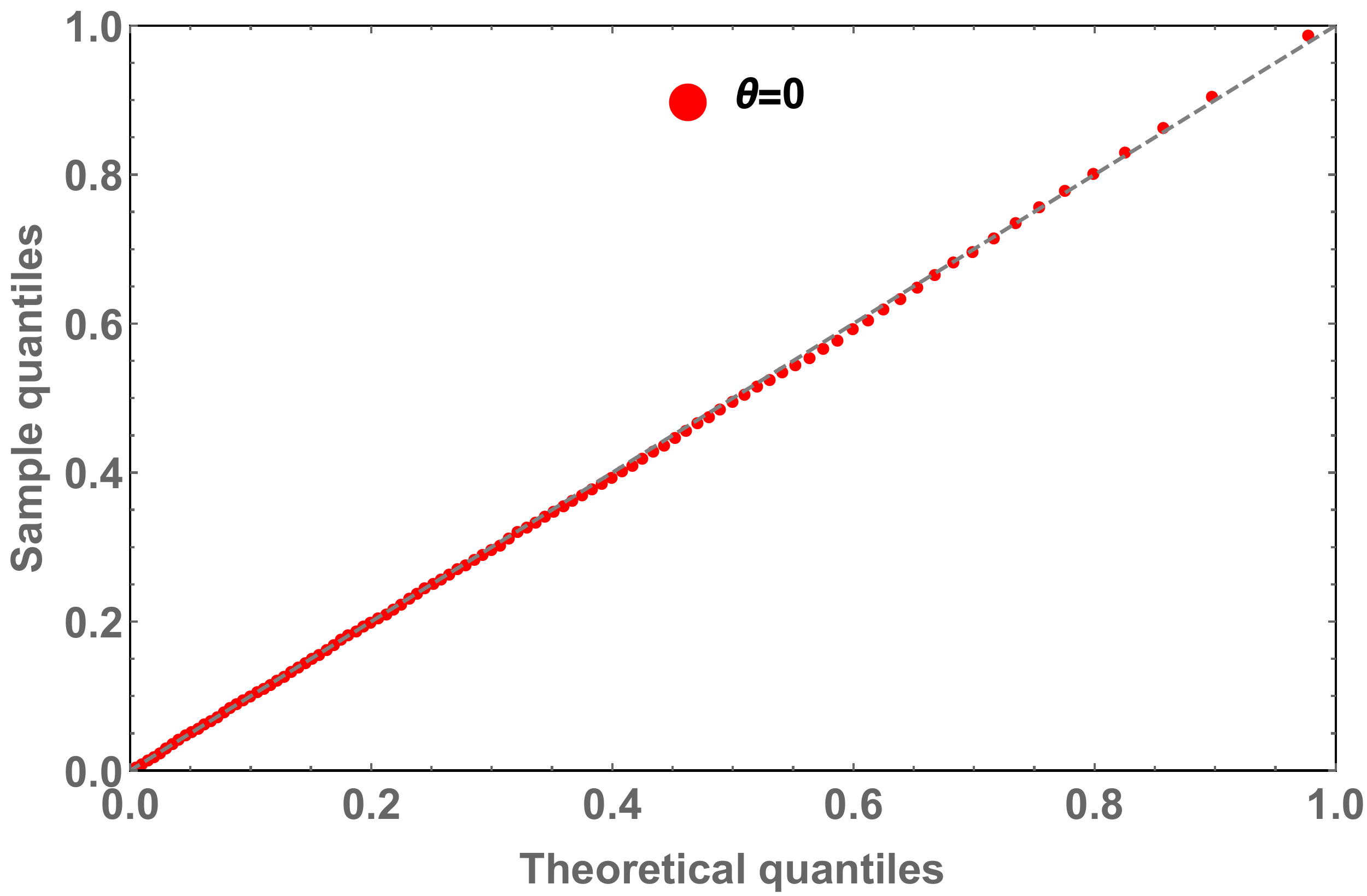}
    \end{subfigure}
    \begin{subfigure}[b]{.45\textwidth}
    \includegraphics[width=\textwidth]{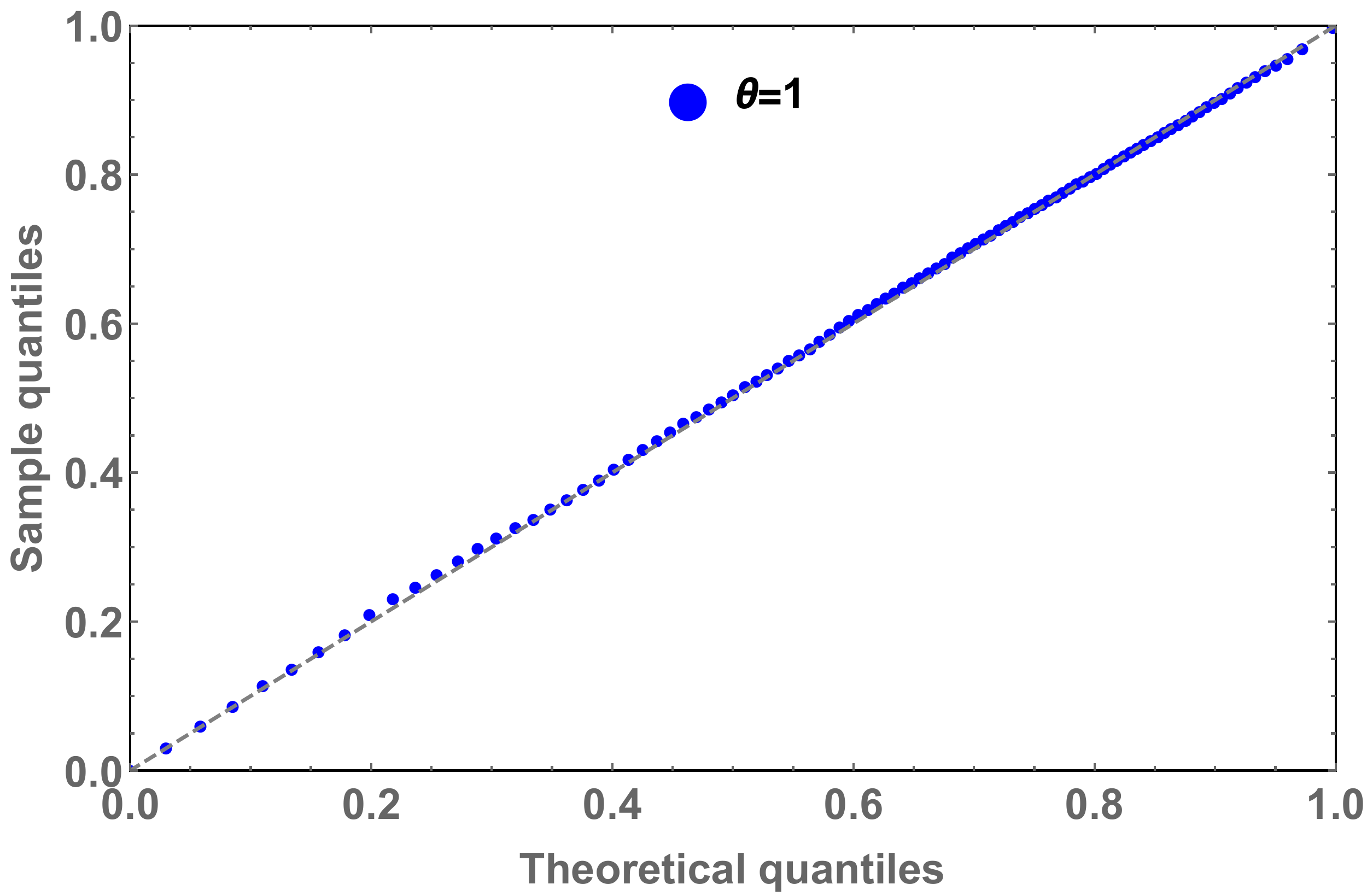}
    \end{subfigure}
    
    \begin{subfigure}[b]{.45\textwidth}
    \includegraphics[width=\textwidth]{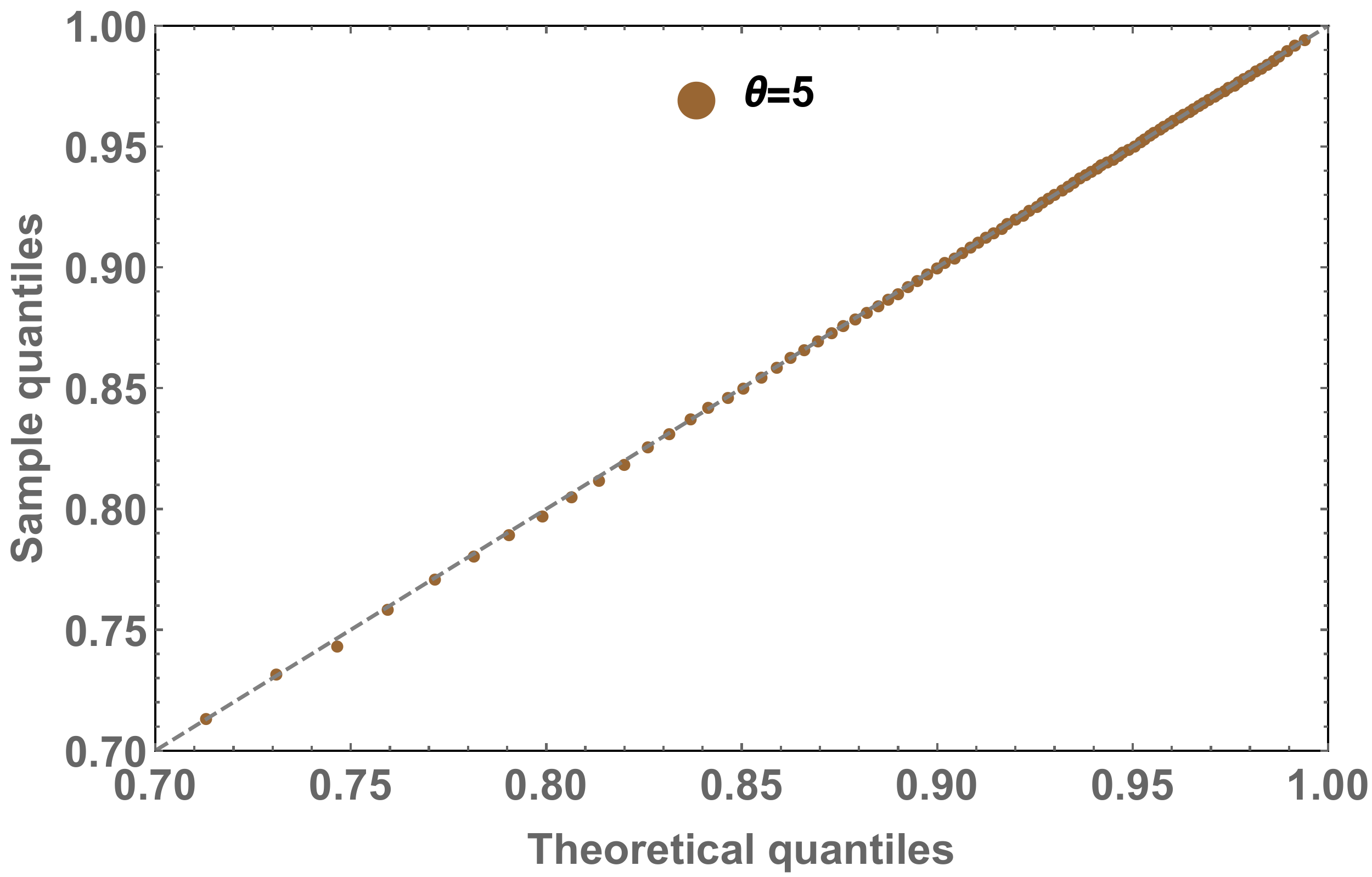}
    \end{subfigure}
    \begin{subfigure}[b]{.45\textwidth}
    \includegraphics[width=\textwidth]{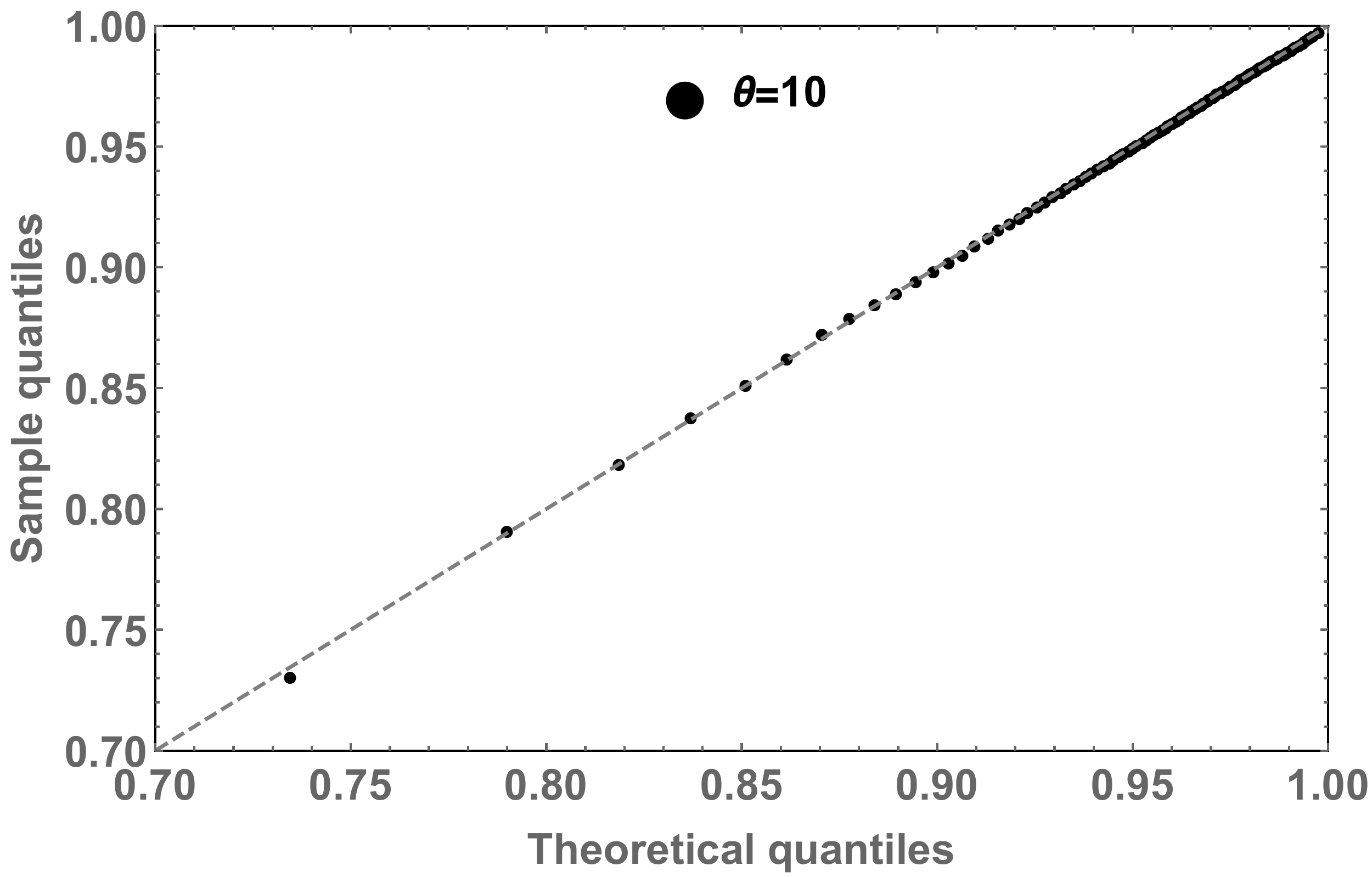}
    \end{subfigure}
    \caption{Quantile-Quantile plots of simulated data of $Z^{(n)}_n$ drawn from  $f^{(\alpha)}_{Z^{(n)}_n}(z)$ for different values of $\theta$ with $\alpha=2$ and $n=3$.}
    \label{fig_qq_maxeigvec_different_theta}
\end{figure}

 Figures \ref{fig_maxeigvec_different_n}, \ref{fig_qq_maxeigvec_different_n},  \ref{fig_maxeigvec_different_theta}, and \ref{fig_qq_maxeigvec_different_theta} verify the accuracy of our formulation. In particular, Fig. \ref{fig_maxeigvec_different_n} shows the effect of $n$ on the p.d.f. of $Z^{(n)}_n$ for $\alpha=2$ and $\theta=3$, whereas Fig. \ref{fig_maxeigvec_different_theta} demonstrates the effect of $\theta$ for $n=3$ and $\alpha=2$. The corresponding Q-Q plots are shown in \ref{fig_qq_maxeigvec_different_n} and \ref{fig_qq_maxeigvec_different_theta}, respectively.  The theoretical curves corresponding to $n=5,6$ and $n=7$ in Fig. \ref{fig_maxeigvec_different_n} have been generated by numerically evaluating the integrals in (\ref{pdfmaxfinal}) and (\ref{pdfmaxfinalaux}). As can be seen from the figures, our theoretical formulations are corroborated by the simulation results. 

 \subsection{The P.D.F. of $Z_\ell^{(n)},\; \ell\neq 1,n$}

The m.g.f. based machinery, which we have developed, can be readily extended to determine the p.d.f. of $Z_\ell^{(n)}$, for $\ell\neq 1,n$, however, at the expense of increased algebraic complexity. To further strengthen this claim, in what follows, we present a summary of technical arguments which lead to the exact p.d.f. of $Z_2^{(n)}$.

{\color{blue} 
For $\ell=2$, (\ref{eq unitary comp}) specializes to  
\begin{align}
    \Psi_2(\boldsymbol{\lambda},s)=\frac{(n-1)!}{2\pi \mathrm{i}} 
        \oint\limits_{\mathcal{C}}
        \frac{e^\omega}{\displaystyle \left(\omega-\beta\lambda_1\right) \left(s+\omega-\beta\lambda_2\right)\prod_{j=3}^n \left(\omega-\beta \lambda_j\right)} {\rm d}\omega
\end{align}
which in turn enables us to rewrite (\ref{eq mgf cont}) as
\begin{align}
    \mathcal{M}_{Z^{(n)}_2}(s)=C^\beta_{n,\alpha}
   \int\limits_{\mathcal{R}}
   \Delta_n^2(\boldsymbol{\lambda}) \prod_{j=1}^n \lambda_j^\alpha e^{-\lambda_j} \frac{1}{2\pi \mathrm{i}} 
    \oint\limits_{\mathcal{C}}
    \frac{e^\omega}{\displaystyle \left(\omega-\beta\lambda_1\right) \left(s+\omega-\beta\lambda_2\right)\prod_{j=3}^n \left(\omega-\beta \lambda_j\right)} {\rm d}\omega\; {\rm d}\boldsymbol{\Lambda}.
\end{align}
Consequently, the inverse Laplace transform yields
\begin{align}
   f^{(\alpha)}_{Z^{(n)}_2}(z)=\frac{C^\beta_{n,\alpha}}{\beta^{n-2}}
   \int\limits_{\mathcal{R}} 
   e^{\beta \lambda_2 z}
   \Delta_n^2(\boldsymbol{\lambda}) \prod_{j=1}^n \lambda_j^\alpha e^{-\lambda_j}
   \sum_{\substack{k=1\\k\neq 2}}^{n}
   \frac{e^{(1-z)\beta \lambda_k}}{\displaystyle \prod_{\substack{i=1\\i\neq 2,k}}^{n} \left(\lambda_k-\lambda_i\right)}\;
  {\rm d}\boldsymbol{\Lambda}.
\end{align}
To facilitate further analysis, we split the integral into two parts as
\begin{align}
    \label{pdfU2 split}
    f^{(\alpha)}_{Z^{(n)}_2}(z) = \frac{C^\beta_{n,\alpha}}{\beta^{n-2}} \left( \hat{\mathcal{A}}^{(\alpha)}_n\left(z\right) + \hat{\mathcal{B}}^{(\alpha)}_n\left(z\right) \right)
\end{align}
where
\begin{align}
    \hat{\mathcal{A}}^{(\alpha)}_n\left(z\right) =  \int\limits_{\mathcal{R}} 
   e^{\beta \lambda_2 z}
   \Delta_n^2(\boldsymbol{\lambda}) \prod_{j=1}^n \lambda_j^\alpha e^{-\lambda_j}
   \frac{e^{(1-z)\beta \lambda_1}}{\displaystyle \prod_{i=3}^{n} \left(\lambda_1-\lambda_i\right)}\;
  {\rm d}\boldsymbol{\Lambda}
\end{align}
and
\begin{align}
   \hat{\mathcal{B}}^{(\alpha)}_n\left(z\right) = \int\limits_{\mathcal{R}} 
   e^{\beta \lambda_2 z}
   \Delta_n^2(\boldsymbol{\lambda}) \prod_{j=1}^n \lambda_j^\alpha e^{-\lambda_j}
   \sum_{k=3}^{n}
   \frac{e^{(1-z)\beta \lambda_k}}{\left(\lambda_k - \lambda_1\right)\displaystyle \prod_{\substack{i=3\\i\neq k}}^{n} \left(\lambda_k-\lambda_i\right)}\;
  {\rm d}\boldsymbol{\Lambda}.
\end{align}

Let us first consider $\hat{\mathcal{B}}^{(\alpha)}_n\left(z\right)$. Employing the decomposition $\Delta_n^2\left(\boldsymbol{\lambda}\right)=(\lambda_2-\lambda_1)^2\prod_{i=3}^n (\lambda_i-\lambda_1)^2(\lambda_i-\lambda_2)^2 \Delta_{n-2}^2(\boldsymbol{\lambda})$ along with relabelling the variables as $t=\lambda_1$, $x=\lambda_2$ and $t_{j-2}=\lambda_j$,  $j=3,\dots,n$, gives us
\begin{align}
    \label{Bhat with upsilon}
    \hat{\mathcal{B}}^{(\alpha)}_n\left(z\right) &= \int_0^\infty x^\alpha e^{-x\left(1-\beta z\right)} \int_0^{x} \left(x-t\right)^2 t^\alpha e^{-t} \Upsilon\left(t,x,z\right) {\rm d}t {\rm d}x
\end{align}
where
\begin{align}
\label{upsilon definition}
    \Upsilon\left(t,x,z\right) = \int\limits_{\mathcal{D}} \sum_{k=1}^{n-2}
   \frac{e^{(1-z)\beta t_k}}{\left(t_k - t\right)\displaystyle \prod_{\substack{i=1\\i\neq k}}^{n-2} \left(t_k-t_i\right)} \prod_{j=1}^{n-2} \left(t_j-t\right)^2 \left(t_j-x\right)^2 t_j^\alpha e^{-t_j} \Delta^2_{n-2}\left(\boldsymbol{t}\right)  {\rm d}t_j
\end{align}
and $\mathcal{D}=\{x<t_1<\dots<t_{n-2}\}$. An important observation that facilitates further simplification of $\Upsilon(t,x,z)$ can be made at this point. Since the integrand is symmetric with respect to $t_j$'s, the ordered region of integration can be removed. Consequently, all the terms in the summation will contribute the same value to the total. Capitalizing  on this observation, we simplify $\Upsilon(t,x,z)$ to obtain
\begin{align}
    \Upsilon(t, &x, z) \nonumber\\
    &= \frac{1}{(n-3)!} \int_{[x,\infty)^{n-2}} \frac{e^{(1-z)\beta t_1}}{\left(t_1 - t\right)\displaystyle \prod_{i=2}^{n-2} \left(t_1-t_i\right)} \prod_{j=1}^{n-2} \left(t_j-t\right)^2 \left(t_j-x\right)^2 t_j^\alpha e^{-t_j} \Delta^2_{n-2}\left(\boldsymbol{t}\right)  {\rm d}t_j.
\end{align}
Now we introduce the variable transformations $w=t_1-x$ and $w_{j-1}=t_j-x$, $j=2,\dots,n-2$, followed by some algebraic manipulation to yield
\begin{align}
\label{upsilon in w}
    \Upsilon(t, &x,z) \nonumber \\
    &= \frac{e^{-x\left(n-2-\beta[1-z]\right)}}{(n-3)!} \int_0^{\infty} e^{-w\left(1-\beta[1-z]\right)} w^2 (w+x-t) (w+x)^{\alpha} \Tilde{\mathcal{Q}}_{n-3}(w,t,x) {\rm d}w
\end{align}
where
\begin{align}
      \Tilde{\mathcal{Q}}_{n-3}(w,t,x) = \int_{[0,\infty)^{n-3}} \prod_{j=1}^{n-3} (w-w_j)(t-x-w_j)^2 (w_j+x)^\alpha w_j^2 e^{-w_j} \Delta^2_{n-3}(\boldsymbol{w}) {\rm d}w_j. 
\end{align}
Following similar steps as in Appendix \ref{appendix_mehta21}, $\Tilde{\mathcal{Q}}_{n-3}(w,t,x)$ can be evaluated to obtain
\begin{multline}
    \Tilde{\mathcal{Q}}_{n-3}(w,t,x) = \frac{(-1)^{n - 1} (n+\alpha-1)!(n+\alpha-2)! \hat{K}_{n-3,\alpha}}{t^{2\alpha}(w+x-t)^2(w+x)^\alpha } \\
    \times \det \left[L_{n+i-4}^{(2)}\left(w\right) \;\; L_{n+i-2-j}^{(j)}\left(t-x\right) \;\; L_{n+i-k}^{(k-2)}\left(-x\right)\right]_{\substack{i=1,\dots,\alpha+3 \\ j=2,3 \\ k=4,\dots,\alpha+3}}.
\end{multline}
Now, by observing that only the first column of the above determinant depends on $w$, we can re-write (\ref{upsilon in w}) to get
\begin{multline}
    \label{upsilon final}
    \Upsilon(t, x, z) = (-1)^{n-1}\hat{K}_{n-3,\alpha} \frac{(n+\alpha-1)! (n+\alpha-2)!}{(n-3)!} \frac{e^{-x\left(n-2-\beta[1-z]\right)}}{t^{2\alpha}} \\
    \times \det \left[\varphi_i^{(\beta)}(x,t,z) \;\; L_{n+i-2-j}^{(j)}\left(t-x\right) \;\; L_{n+i-k}^{(k-2)}\left(-x\right)\right]_{\substack{i=1,\dots,\alpha+3 \\ j=2,3 \\ k=4,\dots,\alpha+3}}
\end{multline}
with
\begin{align}
    \varphi_i^{(\beta)}(x,t,z) &= \int_0^\infty \frac{e^{-w\left(1-\beta[1-z]\right)}w^2}{w+x-t} L_{n+i-4}^{(2)}\left(w\right) {\rm d}w.
\end{align}
With the intention of simplifying $\varphi_i^{(\beta)}(x,t,z)$ further, we employ the variable transformation $u = w / \vartheta(z)$, where $\vartheta(z) = 1 / \left(1-\beta + \beta z\right)$. Consequently, $\varphi_i^{(\beta)}(x,t,z)$ becomes
\begin{align}
    \label{phi define}
    \varphi_i^{(\beta)}(x,t,z) &= \vartheta^2(z) \int_0^\infty \frac{u^2 e^{-u} L_{n+i-4}^{(2)}\left(\vartheta(z)u\right)}{u + (x-t)/\vartheta(z)} {\rm d}u.
\end{align}
Following \cite{ref:szego} and \cite[Eq. 6.15.1.1]{ref:erdelyi} we obtain
\begin{align}
    L_{n+i-4}^{(2)}\left(\vartheta(z)u\right) &= \left(3\right)_{n+i-4} \sum_{\ell=0}^{n+i-4} \frac{\left(1-\vartheta(z)\right)^\ell u^\ell}{(n+i-4-\ell)! \ell! (3)_\ell} {}_1F_1\left(\ell - (n+i-4); 3+\ell; u\right).
\end{align}
Here, $(a)_b = \Gamma(a + b)/\Gamma(a)$ is the Pochhammer symbol. Substituting the above result back into (\ref{phi define}) gives us
\begin{multline}
    \varphi_i^{(\beta)}(x,t,z) = (3)_{n+i-4} \vartheta^2(z) \sum_{\ell=0}^{n+i-4} \frac{\left(1-\vartheta(z)\right)^\ell}{(n+i-4-\ell)! \ell! (3)_\ell} \\ 
    \times\int_0^\infty \frac{u^{2+\ell} e^{-u}}{u + (x-t)/\vartheta(z)} {}_1F_1\left(\ell - (n+i-4); 3+ \ell; u \right) {\rm d}u.
\end{multline}
The integral in the above expressions can be evaluated with the help of \cite[Eq. 6.15.2.16]{ref:erdelyi} to get
\begin{multline}
    \int_0^\infty \frac{u^{2+\ell} e^{-u}}{u + (x-t)/\vartheta(z)} {}_1F_1\left(\ell - (n+i-4); 3+ \ell; u \right) {\rm d}u \\
    = \Gamma(3+\ell) \Gamma(n+i-3-\ell) \left(\frac{x-t}{\vartheta(z)}\right)^{2+\ell} U\left(n+i-1; 3+\ell; \frac{x-t}{\vartheta(z)}\right),
\end{multline}
which along  with some algebraic manipulation yields
\begin{multline}
    \varphi_i^{(\beta)}(x,t,z) = (n+i-2)! (x-t)^2 \sum_{\ell=0}^{n+i-4} \frac{\left(-\beta (1-z)(x-t)\right)^\ell }{\ell!}  \\
    \times U\left(n+i-1; 3+\ell; (x-t)(1-\beta + \beta z)\right)
\end{multline}
where $U(a;c;x)$ being the confluent hypergeometric function of the second kind\cite{ref:erdelyi}. Having simplified $\varphi_i^{(\beta)}(x,t,z)$, we can now substitute (\ref{upsilon final}) into (\ref{Bhat with upsilon}) to obtain
\begin{align}
    \label{Bhat final}
    \hat{\mathcal{B}}^{(\alpha)}_n\left(z\right) = & (-1)^{n-1} \hat{K}_{n-3,\alpha}\frac{(n+\alpha-1)! (n+\alpha-2)!}{(n-3)!} \int_0^\infty x^\alpha e^{-x\left(n-1-\beta\right)} \nonumber  \\
    & \times \int_0^{x} \left(x-t\right)^2 t^{-\alpha} e^{-t} \det \left[\varphi_i^{(\beta)}(x,t,z) \;\; L_{n+i-2-j}^{(j)}\left(t-x\right) \;\; L_{n+i-k}^{(k-2)}\left(-x\right)\right]_{\substack{i=1,\dots,\alpha+3 \\ j=2,3 \\ k=4,\dots,\alpha+3}} 
    {\rm d}t {\rm d}x.
\end{align}
The evaluation of $\hat{\mathcal{A}}_n^{(\alpha)}(z)$ follows along similar line of arguments as before and hence is omitted. As such, we obtain
\begin{multline}
    \label{Ahat final}
    \hat{\mathcal{A}}_n^{(\alpha)}(z) = \frac{(-1)^n \hat{K}_{n-2,\alpha}}{(n-2)!} \int_0^\infty x^\alpha e^{-x(n-1-\beta z)} \int_0^x (x-t)^2 e^{-t(1-\beta + \beta z)} \\
    \times \det \left[L_{n+i-3}^{(2)}\left(t-x\right) \; L_{n+i-1-j}^{(j)}\left(-x\right)\right]_{\substack{i=1,\dots,\alpha+1 \\ j=2,\dots,\alpha+1}} {\rm d}t {\rm d}x.
\end{multline}
By substituting (\ref{Bhat final}) and (\ref{Ahat final}) into (\ref{pdfU2 split}) along with the variable transformation $y=1 - \frac{t}{x}$, we arrive at the final result which is stated in  Theorem \ref{thm u2}.
} 

\begin{theorem}\label{thm u2}
 Let $\mathbf{W}\sim\mathcal{CW}_n\left(m,\mathbf{I}_n+\theta \mathbf{v}\mathbf{v}^\dagger\right)$ with $||\mathbf{v}||=1$ and $\theta>0 $. Let $\mathbf{u}_2$ be the eigenvector corresponding to the {\it second smallest} eigenvalue of $\mathbf{W}$. Then the p.d.f. of $Z^{(n)}_2=|\mathbf{v}^\dagger \mathbf{u}_2|^2\in(0,1)$ is given by
    \begin{align}
    \label{eq pdf second}
       f_{Z^{(n)}_2}^{(\alpha)}(z) &= (-1)^n \beta^{\alpha+2} \left(\frac{1}{\beta}-1\right)^{n+\alpha} \int_0^{\infty} x^5 e^{-x(n-\beta)} \int_0^{1} y^2 e^{xy}
       \mathcal{Z}_n^{(\alpha,\beta)}(x,y,z)
       {\rm d}y {\rm d}x
    \end{align}
    where 
    \begin{align}
        \mathcal{Z}_n^{(\alpha,\beta)}(x,y,z)=\frac{(n-1)!}{(n+\alpha-1)!}  e^{-\beta x y(1-z)} x^{\alpha-2} \mathcal{V}_{n}^{(\alpha)}(x,y) -  \frac{y^2}{ \left(1-y\right)^{\alpha}} \mathcal{U}_{n}^{(\alpha,\beta)}(x,y,z),
    \end{align}
    \begin{align}
        \mathcal{V}_{n}^{(\alpha)}(x,y) &= \det\left[L^{(2)}_{n+i-3}(-xy) \;\;\; L^{(j)}_{n+i-j-1}(-x)\right]_{\substack{i=1,\dots,\alpha+1 \\ j=2,\dots,\alpha+1}},
    \end{align}
    \begin{align}
        \mathcal{U}_{n}^{(\alpha,\beta)}(x,y,z) &= \det\left[\varphi_{i}^{(\beta)}(x,y,z) \;\; L_{n+i-j-2}^{(j)}(-xy) \;\;L_{n+i-k}^{(k-2)}(-x)\right]_{\substack{i=1,\dots,\alpha+3\\j=2,3\\k=4,\dots,\alpha+3}},
    \end{align}
    \begin{align}
        \varphi_{i}^{(\beta)}&(x,y,z) &= (n+i-2)! \sum_{\ell=0}^{n+i-4} \frac{\left(-\beta xy(1-z)\right)^{\ell}}{\ell!} U\left(n+i-1; 3+\ell; xy(1-\beta+\beta z)\right),
    \end{align}
    and $U\left(a;c;z\right)$ is the confluent hypergeometric function of the second kind \cite{ref:erdelyi}.
\end{theorem}
 One of the most important features in the above formula is that the dimensions of the square matrices, whose determinants are in the integrand, depend only on the relative difference $m-n$. Consequently, it is plausible  that this further facilities an asymptotic analysis akin to what we have demonstrated in Corollary \ref{cor min asy}; nevertheless, we do not pursue it here.

\begin{figure}
    \includegraphics[width=0.9\textwidth]{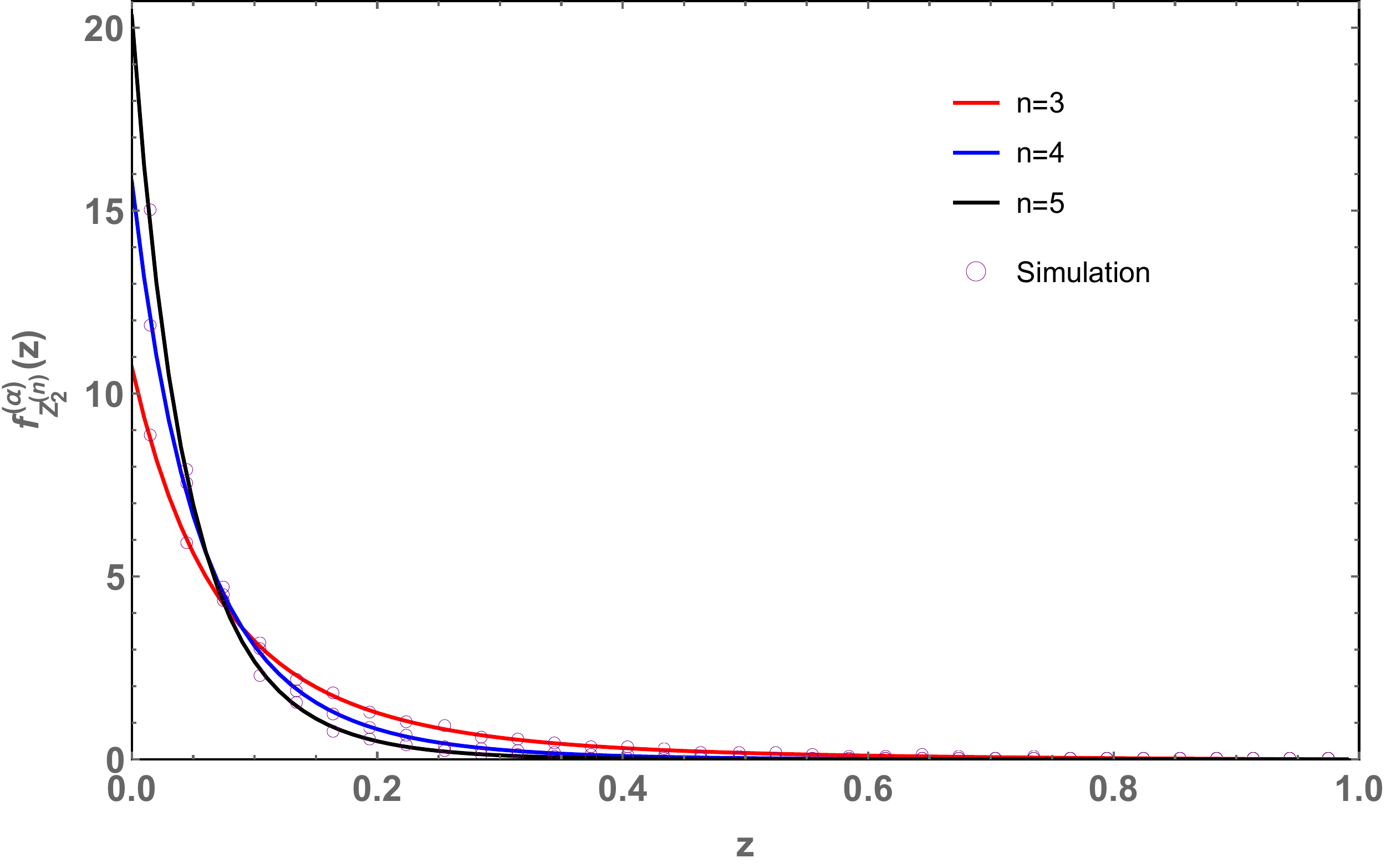}
    \caption{Comparison of simulated data points and the analytical p.d.f. $f^{(\alpha)}_{Z^{(n)}_2}(z)$ for different values of $n$ with $\alpha=1$ and $\theta=3$.}
    \label{fig_seceigvec_different_n}
\end{figure}

\begin{figure}
    \centering
    \begin{subfigure}[b]{.45\textwidth}
    \includegraphics[width=\textwidth]{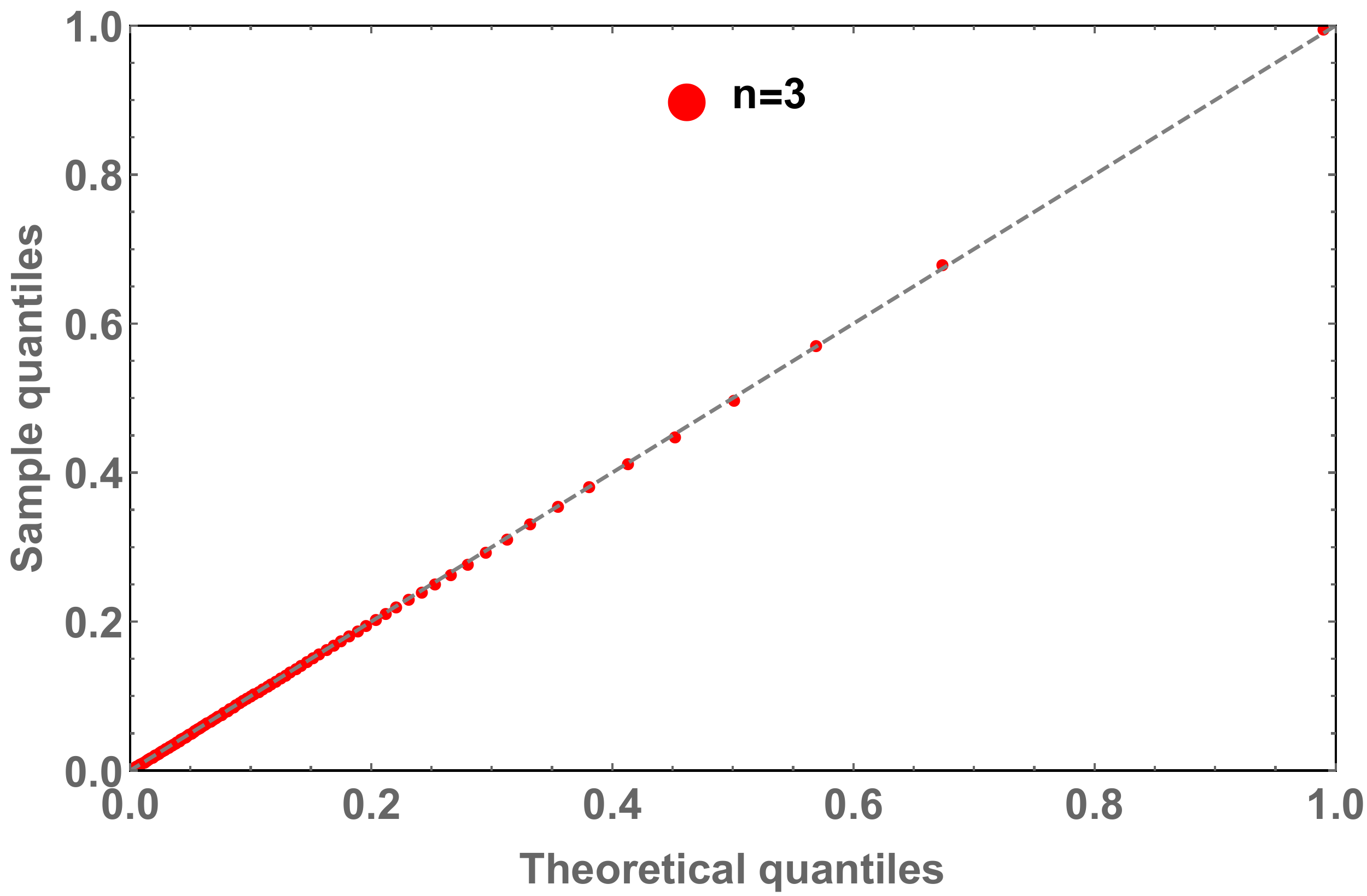}
    \end{subfigure}
    \begin{subfigure}[b]{.45\textwidth}
    \includegraphics[width=\textwidth]{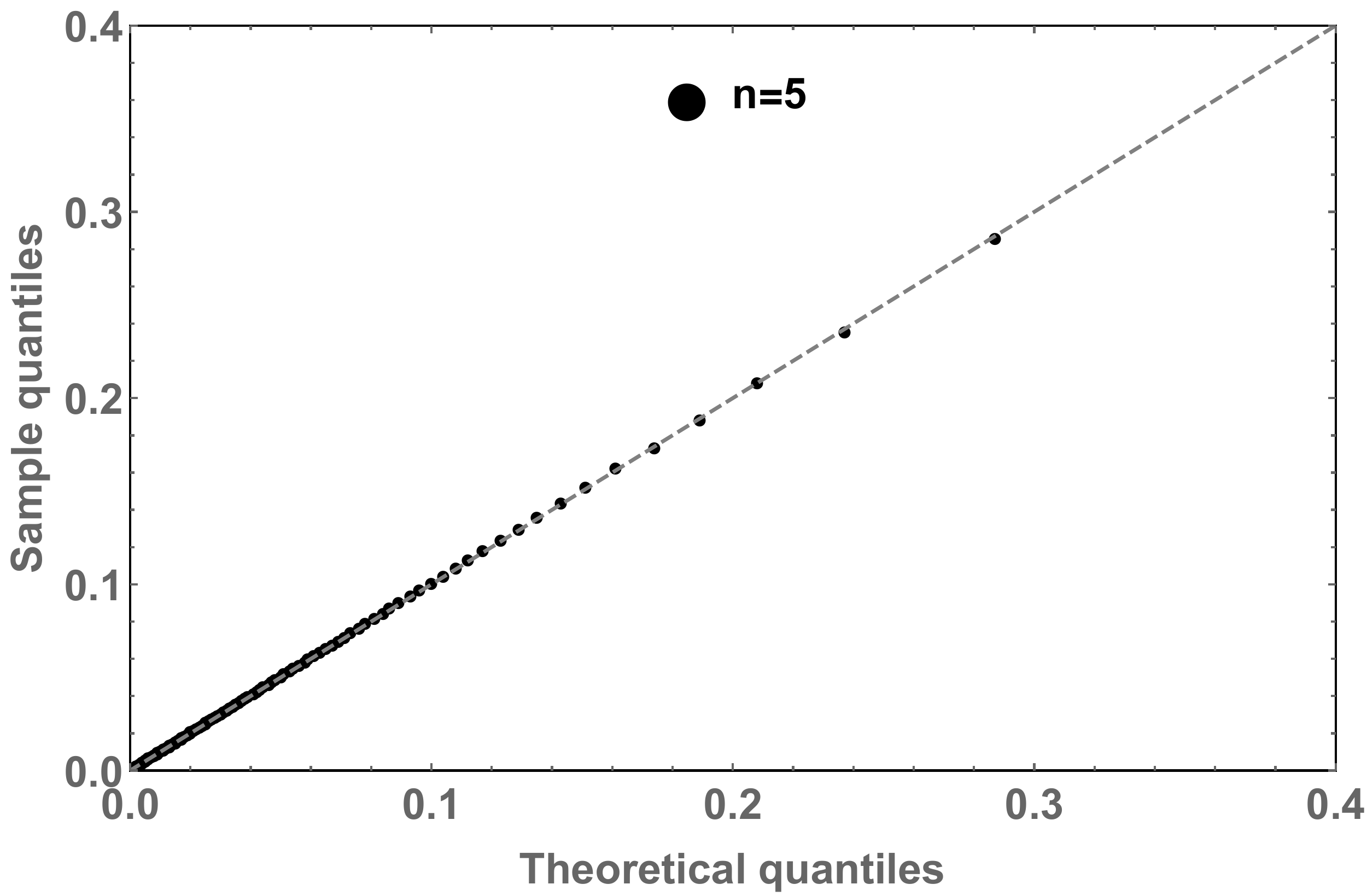}
    \end{subfigure}
    \caption{Quantile-Quantile plots of simulated data of $Z^{(n)}_2$ drawn from  $f^{(\alpha)}_{Z^{(n)}_2}(z)$ for different values of $n$ with $\alpha=1$ and $\theta=3$.}
     \label{fig_qq_seceigvec_different_n}
\end{figure}

To further strengthen our claim in the above theorem, we compare the simulated data points and the analytical p.d.f. $f^{(\alpha)}_{Z^{(n)}_2}(z)$ in Fig. \ref{fig_seceigvec_different_n} for different values of $n$ with $\theta=3$ and $\alpha=1$. Here the analytical p.d.f. is obtained by numerically evaluating the double integral in (\ref{eq pdf second}). Q-Q plots in this respect are depicted in Fig. \ref{fig_qq_seceigvec_different_n}.The agreement between the analytical and simulation results is  clearly evident from the figures; moreover, this verifies the accuracy of Theorem \ref{thm u2}.

\begin{figure}
    \centering
    \includegraphics[width=0.9\textwidth]{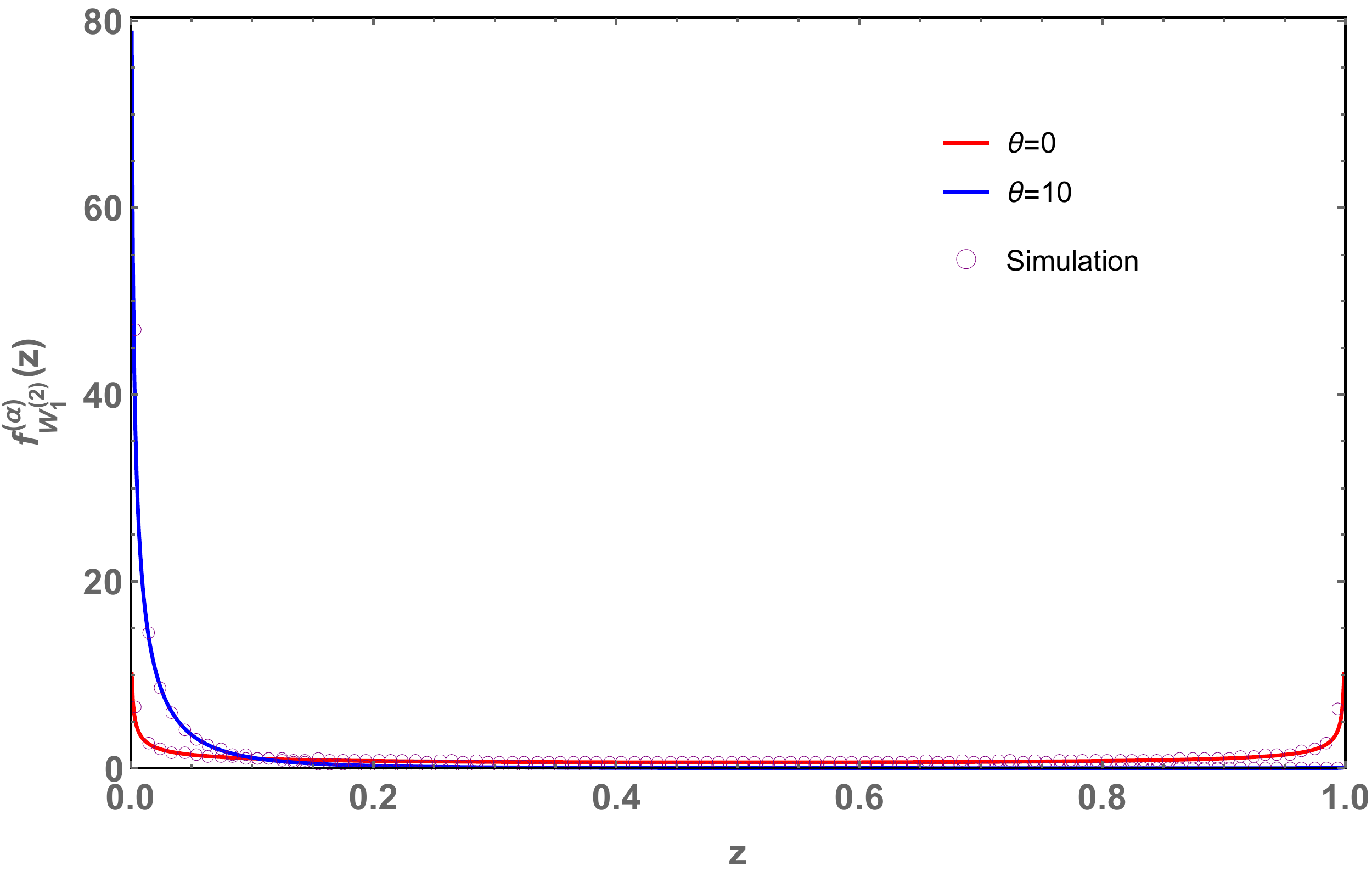}
    \caption{Comparison of simulated data points and the analytical p.d.f. $f^{(\alpha)}_{W^{(2)}_1}(z)$ for different values of $\theta$ with $n=2$ and $\alpha=3$.}
    \label{fig_real_mineigvec_different_theta}
\end{figure}

\begin{figure}
    \centering
    \begin{subfigure}[b]{.45\textwidth}
    \includegraphics[width=\textwidth]{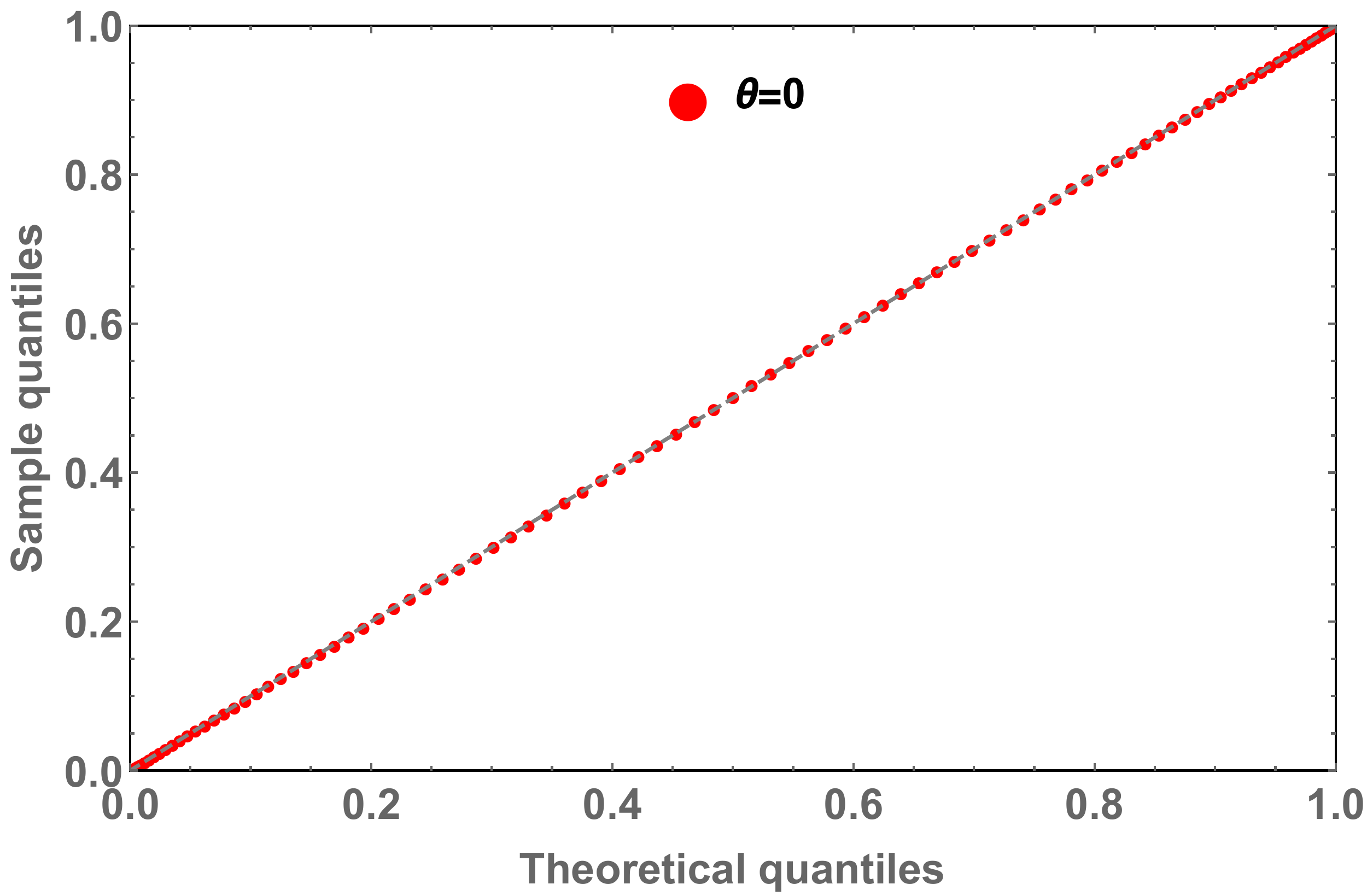}
    \end{subfigure}
    \begin{subfigure}[b]{.45\textwidth}
    \includegraphics[width=\textwidth]{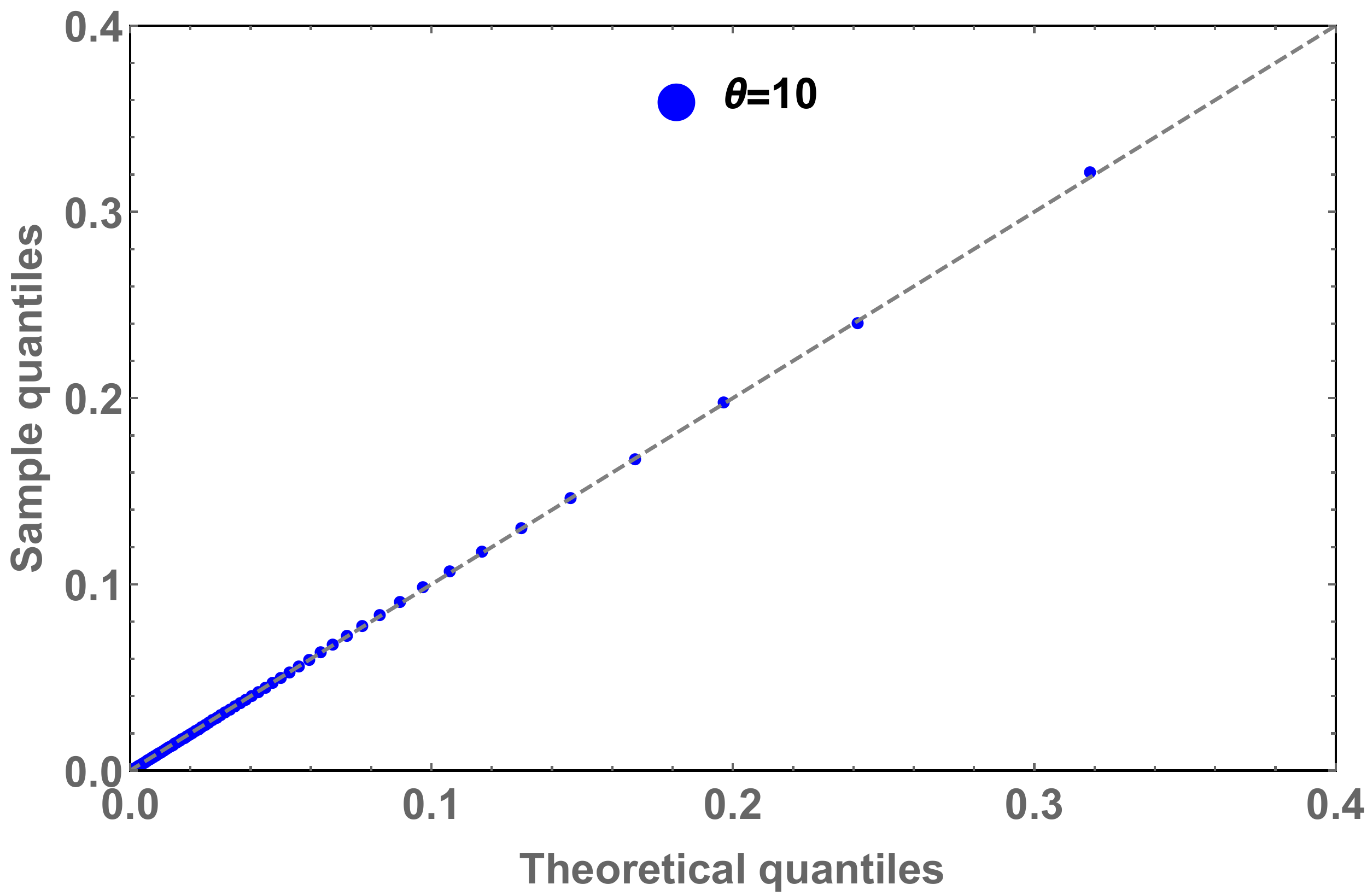}
    \end{subfigure}
    \caption{Quantile-Quantile plots of simulated data of $W^{(2)}_1$ drawn from  $f^{(\alpha)}_{W^{(2)}_1}(z)$ for different values of $\theta$ with $n=2$ and $\alpha=3$.}
     \label{fig_qq_real_mineigvec_different_theta}
\end{figure}

{\color{blue}\section{Real and Singular Wishart Extensions}
Here we provide a detailed analysis on how to extend the results derived in the former sections to real and singular (i.e., $m<n$) Wishart scenarios. As shown in the sequel, the main technical challenge in this respect is to evaluate the corresponding contour integrals in closed-form which are amenable to further analysis.
Although, in general, no closed-form solutions are possible, we  demonstrate a few important special cases for which closed-form solutions lead to exact p.d.f.s.


\subsection{Real Wishart Case}
 To shed some light on this scenario,
  following \cite{ref:james}, \cite{ref:wang}, \cite{ref:dharma}, \cite{ref:onatskiSphe} let us rewrite the real analogy of Theorem \ref{hcizint}  as
\begin{align}
 \label{hcizreal}
 {}_0{\mathsf{F}}_0(\mathbf{vv}^T,\mathbf{T})=
 \int\limits_{\mathcal{O}_n} \text{etr}\left(\mathbf{vv}^T\mathbf{O}\mathbf{TO}^T\right) {\rm d}\mathbf{O}
 = \frac{\Gamma(n/2)}{2\pi \rm i}\oint\limits_{\mathcal{K}}
 \frac{e^{\omega}}{\displaystyle \prod_{j=1}^n\left(\omega-\tau_j\right)^{\frac{1}{2}}} {\rm d}\omega
 \end{align}
where ${\rm d}\mathbf{O}$ is the invariant measure on the orthogonal group $\mathcal{O}_n$ normalized to make the total measure unity, ${}_0{\mathsf{F}}_0(\cdot, \cdot)$ is the real hypergeometric function of two matrix arguments \cite{ref:james,ref:muirhead},  the contour $\mathcal{K}$ starts from $-\infty$ and encircles the eigenvalues of $\mathbf{T}\in \mathbb{R}^{n\times n}$ given by $\tau_1,\tau_2, \cdots,\tau_n$ in the positive direction (i.e., counter-clockwise) and goes back to $-\infty$, $\mathbf{v}\in\mathbb{R}^n$ with $||\mathbf{v}||=1$, and $(\cdot)^T$ represents the transpose operator. Now, for $\mathbf{W}\sim \mathcal{W}_n\left(m,\mathbf{I}_n+\theta \mathbf{vv}^T\right)$, we have from \cite{ref:james}
\begin{align}
    f(\mathbf{W}){\rm d}\mathbf{W}=K^\theta_{m,n} \prod_{k=1}^n
    \lambda_k^{\frac{1}{2}(\alpha-1)} e^{-\frac{\lambda_k}{2}}{\Delta}_n(\boldsymbol{\lambda})
    \text{etr}\left(\frac{\beta}{2} \mathbf{vv}^T \mathbf{O}\boldsymbol{\Lambda}\mathbf{O}^T\right) {\rm d}\boldsymbol{\Lambda}{\rm d}\mathbf{O}
\end{align}
where
\begin{align}
    K^\theta_{m,n}=\frac{\pi^{\frac{n}{2}} 2^{-\frac{mn}{2}} (1+\theta)^{-\frac{m}{2}}}{\prod_{k=1}^n\Gamma\left(\frac{m-k+1}{2}\right)\Gamma\left(\frac{n-k+1}{2}\right) }.
\end{align}
Let us denote the $\ell$th column of $\mathbf{O}$ as $\mathbf{o}_\ell$. Consequently, the m.g.f. of $W_\ell ^{(n)}=|\mathbf{v}^T\mathbf{o}_\ell|^2$ can be written as
\begin{align}
    \mathcal{M}_{W_\ell ^{(n)}}(s)=K^\theta_{m,n} \int_{\mathcal{R}}  \prod_{k=1}^n
    \lambda_k^{\frac{1}{2}(\alpha-1)} e^{-\frac{\lambda_k}{2}}{\Delta}_n(\boldsymbol{\lambda})
    \int_{\mathcal{O}_n} \text{etr}\left\{\mathbf{vv}^T \mathbf{O} \left(\frac{\beta}{2} \boldsymbol{\Lambda}-s\mathbf{e}_\ell \mathbf{e}^T_\ell\right)\mathbf{O}^T\right\} {\rm d}\mathbf{O} {\rm d}\boldsymbol{\Lambda}
\end{align}
from which we obtain, in view of (\ref{hcizreal}), the m.g.f. of $W_1^{(n)}$ as
\begin{align}
    \mathcal{M}_{W_1 ^{(n)}}(s)=K^\theta_{m,n} \Gamma\left(\frac{n}{2}\right) 
    \int_{\mathcal{R}}  \prod_{k=1}^n
    \lambda_k^{\frac{1}{2}(\alpha-1)} e^{-\frac{\lambda_k}{2}}{\Delta}_n(\boldsymbol{\lambda})
    \Psi (\boldsymbol{\lambda}, s) {\rm d} \boldsymbol{\Lambda}
\end{align}
where
\begin{align}
    \Psi (\boldsymbol{\lambda}, s) =
    \frac{1}{2\pi \rm i}\oint\limits_{\mathcal{K}}\frac{e^\omega}{\left(s+\omega-\frac{\beta}{2}\lambda_1\right)^{\frac{1}{2}}\displaystyle \prod_{k=2}^n \left(\omega-\frac{\beta}{2}\lambda_k\right)^{\frac{1}{2}}}
    {\rm d}\omega.
\end{align}
To facilitate further analysis, following \cite[Eq. 5.3.21]{ref:batemanTransforms}, let us take the inverse Laplace transform of the above m.g.f. to yield
\begin{align}
   f^{(\alpha)}_{W_1 ^{(n)}}(z)= \frac{K^\theta_{m,n} \Gamma\left(\frac{n}{2}\right)}{\sqrt{\pi}} z^{-\frac{1}{2}}
   \int_{\mathcal{R}}  \prod_{k=1}^n
    \lambda_k^{\frac{1}{2}(\alpha-1)} e^{-\frac{\lambda_k}{2}}{\Delta}_n(\boldsymbol{\lambda}) e^{\frac{\beta}{2}\lambda_1 z}
    \frac{1}{2\pi \rm i}\oint\limits_{\mathcal{K}}\frac{e^{\omega(1-z)}}{\displaystyle \prod_{k=2}^n \left(\omega-\frac{\beta}{2}\lambda_k\right)^{\frac{1}{2}}}
    {\rm d}\omega {\rm d}\boldsymbol{\Lambda},
\end{align}
from which we obtain in view of \cite[Eq. 5.4.9]{ref:batemanTransforms}
\begin{align}
\label{remin}
    f^{(\alpha)}_{W_1 ^{(n)}}(z)& = \frac{K^\theta_{m,n} \Gamma\left(\frac{n}{2}\right)}{\sqrt{\pi} \Gamma\left(\frac{n-1}{2}\right)}  z^{-\frac{1}{2}}(1-z)^{\frac{n-3}{2}}\nonumber\\
  & \qquad \times  \int_{\mathcal{R}}  \prod_{k=1}^n
    \lambda_k^{\frac{1}{2}(\alpha-1)} e^{-\frac{\lambda_k}{2}}{\Delta}_n(\boldsymbol{\lambda}) e^{\frac{\beta}{2}\lambda_1 z}\nonumber\\
    & \qquad \qquad \times 
    \Phi_2^{(n-1)}\left(\frac{1}{2},\ldots,\frac{1}{2};\frac{n-1}{2}; \frac{\beta}{2}\lambda_2(1-z),\ldots,\frac{\beta}{2}\lambda_n(1-z)\right) {\rm d}\boldsymbol{\Lambda}
\end{align}
where $\Phi_2^{(n)}(\cdot)$ is the confluent form of the generalized Lauricella series\footnote{To be precise, it assumes the infinite series expansion given by \cite{ref:erdelyiGer}
\begin{align*}
    \Phi_2^{(n)}(b_1,b_2,\cdots,b_n;c;x_1,x_2,\cdots,x_n)=\sum_{k_1,k_2,\ldots,k_n=0}^\infty \frac{(b_1)_{k_1} (b_2)_{k_2}\ldots (b_n)_{k_n}}{(c)_{k_1+k_2+\cdots+k_n}} \frac{x_1^{k_1}}{k_1 !} \frac{x_2^{k_2}}{k_2!} \ldots \frac{x_n^{k_n}}{k_n!}.
\end{align*}} defined in \cite[Eq. 7.2]{ref:erdelyiGer}.  Similarly, we can write the p.d.f. of $W_n^{(n)}$ as
\begin{align}
\label{remax}
    f^{(\alpha)}_{W_n ^{(n)}}(z)& = \frac{K^\theta_{m,n} \Gamma\left(\frac{n}{2}\right)}{\sqrt{\pi} \Gamma\left(\frac{n-1}{2}\right)}  z^{-\frac{1}{2}}(1-z)^{\frac{n-3}{2}}\nonumber\\
  & \qquad \times  \int_{\mathcal{R}}  \prod_{k=1}^{n}
    \lambda_k^{\frac{1}{2}(\alpha-1)} e^{-\frac{\lambda_k}{2}}{\Delta}_n(\boldsymbol{\lambda}) e^{\frac{\beta}{2}\lambda_n z}\nonumber\\
    & \qquad \qquad \times 
    \Phi_2^{(n-1)}\left(\frac{1}{2},\ldots,\frac{1}{2};\frac{n-1}{2}; \frac{\beta}{2}\lambda_1(1-z),\ldots,\frac{\beta}{2}\lambda_{n-1}(1-z)\right) {\rm d}\boldsymbol{\Lambda}.
\end{align}

Further manipulation of the above multiple integrals is highly challenging due to the presence of $\Phi_2^{(n-1)}$ in each of the integrands. Therefore, unlike in the complex case, a closed-form p.d.f. for $W_1^{(n)}$ seems intractable for general values of $m$ and $n$\footnote{The p.d.f. of $W_1^{(n)}$, for $\theta=0$, can be obtained by setting $\beta=0$ in (\ref{remin}) as
$
f^{(\alpha)}_{W_1 ^{(n)}}(z) = \frac{ \Gamma\left(\frac{n}{2}\right)}{\sqrt{\pi} \Gamma\left(\frac{n-1}{2}\right)}  z^{-\frac{1}{2}}(1-z)^{\frac{n-3}{2}}
$. Moreover, the same p.d.f. result remains valid for the p.d.f. of $W_\ell^{(n)},\; \ell=2,3,\ldots,n$.}. Nevertheless, it is noteworthy that in the important case of $n=2$, $\Phi_2^{(1)}$ degenerates to an exponential function. Capitalizing on this observation, we may write the p.d.f. of 
$W_1^{(2)}$ as
\begin{align}
    f^{(\alpha)}_{W_1 ^{(2)}}(z) = \frac{K^\theta_{m,2}}{\pi}  z^{-\frac{1}{2}}(1-z)^{-\frac{1}{2}}  \iint\limits_{0<\lambda_1<\lambda_2<\infty} \prod_{k=1}^2
    \lambda_k^{\frac{1}{2}(m-3)} e^{-\frac{\lambda_k}{2}}{\Delta}_2(\boldsymbol{\lambda}) e^{\frac{\beta}{2}\left(\lambda_1 z+\lambda_2(1-z)\right)} {\rm d}\lambda_1 {\rm d}\lambda_2,
\end{align}
from which we obtain, after some algebraic manipulation,
\begin{align}
    f^{(\alpha)}_{W_1 ^{(2)}}(z) & = \frac{K^\theta_{m,2}}{\pi}  z^{-\frac{1}{2}}(1-z)^{-\frac{1}{2}}\nonumber\\
    & \qquad \times \iint\limits_{0<\lambda_1<\lambda_2<\infty}
    \left(\lambda_1^{\frac{m-3}{2}} \lambda_2^{\frac{m-1}{2}}-\lambda_1^{\frac{m-1}{2}} \lambda_2^{\frac{m-3}{2}}\right)
    e^{-\frac{\lambda_1}{2}\left( 1-\beta z\right)}
   e^{-\frac{\lambda_2}{2}\left[ 1-\beta(1- z)\right]} {\rm d}\lambda_1 {\rm d}\lambda_2.
\end{align}
Now we may evaluate the above double integral with the help of \cite[Eq. 3.194.1]{ref:gradshteyn} to yield
\begin{align}
   f^{(\alpha)}_{W_1 ^{(2)}}(z) &= \frac{2^{m-1}(m-1)}{\pi (1+\theta)^{\frac{m}{2}}} z^{-\frac{1}{2}}(1-z)^{-\frac{1}{2}} \left[1-\beta(1-z)\right]^{-m}\nonumber\\
   & \qquad \times \left[
   \frac{1}{m-1} {}_2F_1\left(m,\frac{m-1}{2};\frac{m+1}{2};-\frac{1-\beta z}{1-\beta(1-z)}\right)\right.\nonumber\\
   & \qquad \qquad \qquad \qquad -\left.\frac{1}{m+1} {}_2F_1\left(m,\frac{m+1}{2};\frac{m+3}{2};-\frac{1-\beta z}{1-\beta(1-z)}\right) \right]. 
\end{align}
Consequently, noting that $W_2^{(2)}+W_1^{(2)}=1$, we obtain the p.d.f. of 
   $W_2^{(2)}$ as $f^{(\alpha)}_{W_2 ^{(2)}}(z)=f^{(\alpha)}_{W_1 ^{(2)}}(1-z)$.}
Finally, Fig. \ref{fig_real_mineigvec_different_theta} and the accompanying Q-Q plot in Fig. \ref{fig_qq_real_mineigvec_different_theta} verify the accuracy of our result for different values of $\theta$ with $n=2$ and $\alpha=3$.

{\color{blue}\subsection{Singular Complex Wishart Case}
Another scenario of technical interest is when the number of observations $m$ is less than the dimension $n$ of the random vectors (i.e., $m<n$). In this situation, the matrix $\mathbf{W}$ degenerates to the so called singular Wishart matrix \cite{ref:ratnarajahSing,ref:uhlig} whose rank is $m$ almost surely. As such, the density of $\mathbf{W}$ is defined on the space of $m\times m$ Hermitian positive {\it semi-definite} matrices of rank $m$ \cite{ref:ratnarajahSing}. Consequently, $\mathbf{W}$ assumes the eigen-decomposition $\mathbf{W}=\mathbf{U}_1\boldsymbol{\Lambda}\mathbf{U}_1^\dagger$ with $\mathbf{U}_1^\dagger \mathbf{U}_1=\mathbf{I}_m$  and $\boldsymbol{\Lambda}=\text{diag}\left(\lambda_1, \lambda_2,\ldots,\lambda_m\right)$ are the non-zero eigenvalues of $\mathbf{W}$ ordered such that $0<\lambda_1<\lambda_2<\ldots<\lambda_m<\infty$. It is noteworthy that the set of all $n\times m$ complex matrices $\mathbf{U}_1$ such that $\mathbf{U}_1^\dagger \mathbf{U}_1=\mathbf{I}_m$ (i.e., with orthonormal columns), denoted by $\mathcal{V}_{m,n}$, is known as the complex {\it Stiefel manifold}.

Following \cite[Eqs. 2, 22]{ref:ratnarajahSing}, the joint density of the eigenvalues and eigenvectors of $\mathbf{W}\sim\mathcal{CW}_n\left(m, \mathbf{I}_n+\theta \mathbf{vv}^\dagger\right)$, for $m<n$, can be written as
\begin{align}
    f(\boldsymbol{\Lambda}, \mathbf{U_1})=\frac{\pi^{m(m-n-1)}2^{-m}}{\tilde{\Gamma}_m(n) (1+\theta)^{m}}
    \prod_{k=1}^m \lambda_k^{n-m} e^{-\lambda_k} \Delta^2_m(\boldsymbol{\lambda})
    \text{etr}\left(\beta \mathbf{vv}^\dagger \mathbf{U}_1\boldsymbol{\Lambda}\mathbf{U}_1^\dagger \right).
\end{align}  
Therefore, the m.g.f. of $Y_{\ell}^{(n)}=\left|\mathbf{v}^\dagger \mathbf{u}_\ell\right|,\; \ell=1,2,\ldots, m$, can be written as
\begin{align}
\label{mgfsingint}
    \mathcal{M}_{Y_{\ell}^{(n)}}(s)=\frac{\pi^{m(m-n-1)}2^{-m}}{\tilde{\Gamma}_m(n) (1+\theta)^{m}}
    \int_{\mathcal{R}}
    \prod_{k=1}^m \lambda_k^{n-m} e^{-\lambda_k} \Delta^2_m(\boldsymbol{\lambda})
    \Psi_{\ell}(\boldsymbol{\lambda},s) {\rm d}\boldsymbol{\Lambda}
\end{align}
where
\begin{align}
\label{singHCIZ}
    \Psi_{\ell}(\boldsymbol{\lambda},s)=\int_{\mathcal{V}_{m,n}} 
    \text{etr}\left\{\mathbf{vv}^\dagger \mathbf{U}_1\left(\beta\boldsymbol{\Lambda}-s\mathbf{g}_\ell \mathbf{g}_\ell^T \right)\mathbf{U}_1^\dagger \right\} \left(\mathbf{U}_1^\dagger {\rm d}\mathbf{U}_1\right)
\end{align}
with $\left(\mathbf{U}_1^\dagger {\rm d}\mathbf{U}_1\right)$ denoting the exterior differential form representing the uniform measure on the complex Stiefel manifold \cite{ref:ratnarajahSing,ref:AOSing} and $\mathbf{g}_\ell$ denotes the $\ell$th column of the $m\times m$ identity matrix. The next technical challenge is to evaluate the matrix integral in (\ref{singHCIZ}). To this end, for convenience, let us consider an equivalent matrix integral given by
\begin{align}
    J(\mathbf{T})=\int_{\mathcal{V}_{m,n}} 
    \text{etr}\left\{\mathbf{vv}^\dagger \mathbf{U}_1\mathbf{T}\mathbf{U}_1^\dagger \right\} \left(\mathbf{U}_1^\dagger {\rm d}\mathbf{U}_1\right)
\end{align}
where $\mathbf{T}=\text{diag}\left(\tau_1,\tau_2,\ldots,\tau_m
\right)\in\mathbb{R}^{m\times m}$ is a diagonal matrix. Since further manipulation of this matrix integral in the current form is an arduous task, following \cite[Appendix A]{ref:AOSing}, we may rewrite it as a matrix integral over the unitary manifold as
\begin{align}
    J(\mathbf{T})=\frac{2^m \pi ^{mn}}{\tilde\Gamma_m(m)}\int_{\mathcal{U}_{n}}
    \text{etr}\left\{\mathbf{vv}^\dagger \mathbf{U}\widetilde{\mathbf{T}}\mathbf{U}^\dagger \right\} {\rm d}\mathbf{U}
\end{align}
where $\mathbf{U}=\left(\mathbf{U}_1\; \mathbf{U}_2\right)\in\mathcal{U}_n$ and $\widetilde{\mathbf{T}}=\text{diag}\left(\tau_1,\tau_2,\ldots,\tau_m, 0, 0, \ldots,0\right)\in\mathbb{R}^{n\times n}$. Now, the above matrix integral yields \cite{ref:james}
\begin{align}
    J(\mathbf{T})=\frac{2^m \pi ^{mn}}{\tilde\Gamma_m(m)}
    {}_0\widetilde{\mathsf F}_0 \left(\mathbf{vv}^\dagger, \widetilde{\mathbf{T}}\right).
\end{align}
Keeping in mind that $\widetilde{\mathbf{T}}$ is rank deficient, we apply Theorem \ref{hcizint} to arrive at
\begin{align}
\label{eqsingmatint}
    \int_{\mathcal{V}_{m,n}} 
    \text{etr}\left\{\mathbf{vv}^\dagger \mathbf{U}_1\mathbf{T}\mathbf{U}_1^\dagger \right\} \left(\mathbf{U}_1^\dagger {\rm d}\mathbf{U}_1\right)=\frac{\mathrm{K}_{m,n}}{2\pi \mathrm{i}} 
        \oint\limits_{\mathcal{C}}
        \frac{\omega ^{m-n}e^\omega}{\displaystyle \prod_{j=1}^m \left(\omega-\tau_j\right)} {\rm d}\omega
\end{align}
where $\mathrm{K}_{m,n}=\frac{2^m \pi ^{mn}\Gamma(n)}{\tilde\Gamma_m(m)}$
Therefore, in view of (\ref{eqsingmatint}), (\ref{singHCIZ}) specializes to
\begin{align}
\Psi_\ell (\boldsymbol{\lambda}, s)=\left\{
    \begin{array}{cc}\displaystyle 
    \frac{\mathrm{K}_{m,n}}{2\pi \mathrm{i}} 
        \oint\limits_{\mathcal{C}}
        \frac{\omega ^{m-n}e^\omega}{(s+\omega-\beta \lambda_1)\displaystyle \prod_{j=2}^m \left(\omega-\beta \lambda_j\right)} {\rm d}\omega & \text{for $\ell=1$}\\
        \displaystyle  
    \frac{\mathrm{K}_{m,n}}{2\pi \mathrm{i}} 
        \oint\limits_{\mathcal{C}}
        \frac{\omega ^{m-n}e^\omega}{(s+\omega-\beta \lambda_m)\displaystyle \prod_{j=1}^{m-1} \left(\omega-\beta \lambda_j\right)} {\rm d}\omega & \text{for $\ell=m$},
\end{array}\right.
\end{align}
from which we obtain after taking the inverse Laplace transform
\begin{align}
\label{singinvlap1}
    \mathcal{L}^{-1}\left\{\Psi_\ell (\boldsymbol{\lambda}, s)\right\}=\left\{\begin{array}{cc}\displaystyle  e^{\beta \lambda_1 z}
    \frac{\mathrm{K}_{m,n}}{2\pi \mathrm{i}} 
        \oint\limits_{\mathcal{C}}
        \frac{\omega ^{m-n}e^{(1-z)\omega}}{\displaystyle \prod_{j=2}^m \left(\omega-\beta \lambda_j\right)} {\rm d}\omega & \text{for $\ell=1$}\\
        \displaystyle  e^{\beta \lambda_m z}
    \frac{\mathrm{K}_{m,n}}{2\pi \mathrm{i}} 
        \oint\limits_{\mathcal{C}}
        \frac{\omega ^{m-n}e^{\omega(1-z)}}{\displaystyle \prod_{j=1}^{m-1} \left(\omega-\beta \lambda_j\right)} {\rm d}\omega & \text{for $\ell=m$}
\end{array}\right.
\end{align}
where $\mathcal{L}^{-1}\{\cdot\}$ denotes the inverse Laplace operator. Now we evaluate the contour integrals, for $m\geq 2$, to yield
\begin{align}
\label{singinteval}
    \mathcal{L}^{-1}\left\{\Psi_\ell (\boldsymbol{\lambda}, s)\right\}=\left\{\begin{array}{cc}\displaystyle \mathrm{K}_{m,n} e^{\beta \lambda_1 z}\displaystyle
  \left(\frac{1}{\beta^{n-2}}\sum_{k=2}^m
        \frac{\lambda_k ^{m-n}e^{(1-z)\beta\lambda_k}}{\displaystyle \prod_{\substack{i=2\\i\neq k}}^m \left(\lambda_k- \lambda_i\right)}+\mathrm{J}^{(n-m)}_1(z,\boldsymbol{\lambda})\right) & \text{for $\ell=1$}\\
        \displaystyle \mathrm{K}_{m,n} e^{\beta \lambda_m z}
    \displaystyle
  \left(\frac{1}{\beta^{n-2}}\sum_{k=1}^{m-1}
        \frac{\lambda_k ^{m-n}e^{(1-z)\beta\lambda_k}}{\displaystyle \prod_{\substack{i=1\\i\neq k}}^{m-1} \left(\lambda_k- \lambda_i\right)}+\mathrm{J}^{(n-m)}_m(z,\boldsymbol{\lambda})\right) & \text{for $\ell=m$}
\end{array}\right.
\end{align}
where
\begin{align}
\label{contintsing}
    \mathrm{J}^{(n-m)}_\ell(z,\boldsymbol{\lambda})=\left\{\begin{array}{cc}\displaystyle \frac{1}{2\pi \mathrm{i}}\oint\limits_{0}
        \frac{\omega ^{m-n}e^{\omega(1- z)}}{\displaystyle \prod_{j=2}^{m} \left(\omega-\beta \lambda_j\right)} {\rm d}\omega & {\text{for $\ell=1$}}\\
        \displaystyle \frac{1}{2\pi \mathrm{i}}\oint\limits_{0}
        \frac{\omega ^{m-n}e^{\omega(1-z)}}{\displaystyle \prod_{j=1}^{m-1} \left(\omega-\beta \lambda_j\right)} {\rm d}\omega & \text{for $\ell=m$}\end{array}\right.
\end{align}
in which the contour $0$ is taken to be a small circle around the origin such that all $\lambda_j,\; j=1,2,\ldots,m$ are in the exterior of the contour. A careful inspection of (\ref{contintsing}) reveals that the pole of order $n-m$ at the origin does not facilitate the direct use of Cauchy integral formula to evaluate these contour integrals\footnote{To be precise, for $\ell=1$, in light of Cauchy integral theorem, we obtain $\mathrm{J}^{(n-m)}_\ell(z,\boldsymbol{\lambda})=g(0,z, \boldsymbol{\lambda})$, where $
  g(\omega,z, \boldsymbol{\lambda})=  \frac{1}{(m-n-1)!} \frac{{\rm d}^{m-n-1}}{{\rm d}\omega ^{m-n-1}}\left[\frac{e^{(1-z)\omega}}{ \prod_{j=2}^m \left(\omega-\beta \lambda_j\right)}
  \right]$. Here it is worth noting that obtaining a general expression for the $(m-n-1)$th derivative of $e^{(1-z)\omega}/ \prod_{j=2}^m \left(\omega-\beta \lambda_j\right)$ is intractable. Similar arguments apply to the case corresponding to $\ell=m$ as well. } for general values of $m$ and $n$ such that $n>m$. To circumvent this difficulty, noting the power series expansion \cite[Eq. 247]{ref:wang}
  \begin{align}
      \frac{1}{\displaystyle \prod_{j=1}^{N} \left(\omega-\beta \lambda_j\right)}=\frac{1}{(-\beta)^N \displaystyle \prod_{j=1}^N \lambda_j}
      \sum_{k=0}^\infty \frac{\omega^k}{(-\beta)^k} C_k\left(\text{diag}\left(\lambda_1^{-1},\ldots, \lambda_N^{-1}\right)\right),
  \end{align}
  where $k\equiv(k,0,0,\ldots,0)$ denotes the length $m$ partition of $k$ into not more than one part and $C_k(\cdot)$ is the complex zonal polynomial,
we may rewrite (\ref{contintsing}) as
\begin{align}
\label{contintsingcauchy}
    \mathrm{J}^{(n-m)}_\ell(z,\boldsymbol{\lambda})=\left\{\begin{array}{cc}\displaystyle \frac{1}{(-\beta)^{m-1} \displaystyle \prod_{j=2}^m \lambda_j}
      \sum_{k=0}^\infty
      \frac{C_k\left(\boldsymbol{\Lambda}_1^{-1}\right)}{(-\beta)^k}
      \frac{1}{2\pi \mathrm{i}}
      \oint\limits_{0}
        \frac{e^{\omega(1- z)}}{\omega^{n-m-k}} {\rm d}\omega & {\text{for $\ell=1$}}\\
        \displaystyle \frac{1}{(-\beta)^{m-1} \displaystyle \prod_{j=1}^{m-1} \lambda_j}
      \sum_{k=0}^\infty
      \frac{C_k\left(\boldsymbol{\Lambda}_m^{-1}\right)}{(-\beta)^k}
      \frac{1}{2\pi \mathrm{i}}
      \oint\limits_{0}
        \frac{e^{\omega(1- z)}}{\omega^{n-m-k}} {\rm d}\omega & \text{for $\ell=m$}\end{array}\right.
\end{align}
where $\boldsymbol{\Lambda}_1=\text{diag}\left(\lambda_2,\ldots, \lambda_m\right)$ and $\boldsymbol{\Lambda}_m=\text{diag}\left(\lambda_1,\ldots, \lambda_{m-1}\right)$. Now, keeping in mind that $ \oint\limits_{0}\frac{e^{\omega(1- z)}}{\omega^{n-m-k}} {\rm d}\omega=0$ for $k\geq n-m$, we may evaluate the above contour integrals to yield
\begin{align}
\label{contintsingcauchy1}
    \mathrm{J}^{(n-m)}_\ell(z,\boldsymbol{\lambda})=\left\{\begin{array}{cc}\displaystyle \frac{1}{(-\beta)^{m-1} \displaystyle \prod_{j=2}^m \lambda_j}
      \sum_{k=0}^{n-m-1}
      \frac{C_k\left(\boldsymbol{\Lambda}_1^{-1}\right)(1-z)^{n-m-1-k}}{(-\beta)^k(n-m-1-k)!}& {\text{for $\ell=1$}}\\
        \displaystyle \frac{1}{(-\beta)^{m-1} \displaystyle \prod_{j=1}^{m-1} \lambda_j}
      \sum_{k=0}^{n-m-1}
      \frac{C_k\left(\boldsymbol{\Lambda}_m^{-1}\right)(1-z)^{n-m-1-k}}{(-\beta)^k(n-m-1-k)!} & \text{for $\ell=m$}.
      \end{array}\right.
\end{align}
The above general forms are not amenable to further analysis due to the availability of the zonal polynomials. However, for $n-m=1,2$, the above formulas simplify to 
\begin{align}
\label{contintsing1}
    \mathrm{J}^{(1)}_\ell(z,\boldsymbol{\lambda})=\left\{\begin{array}{cc}\displaystyle \frac{\lambda_1}{(-\beta)^{m-1} \displaystyle \prod_{j=1}^m \lambda_j} & {\text{for $\ell=1$}}\\
        \displaystyle \frac{\lambda_m}{(-\beta)^{m-1} \displaystyle \prod_{j=1}^m \lambda_j} & \text{for $\ell=m$}\end{array}\right.
\end{align}
and 
\begin{align}
\label{contintsing2}
    \mathrm{J}^{(2)}_\ell(z,\boldsymbol{\lambda})=\left\{\begin{array}{cc}\displaystyle \frac{(1-z)\lambda_1}{(-\beta)^{m-1} \displaystyle \prod_{j=1}^m \lambda_j}-\frac{\lambda_1}{(-\beta)^{m} \displaystyle \prod_{j=1}^m \lambda_j}\sum_{k=2}^m \frac{1}{\lambda_k} & {\text{for $\ell=1$}}\\
        \displaystyle \frac{(1-z)\lambda_m}{(-\beta)^{m-1} \displaystyle \prod_{j=1}^m \lambda_j}-\frac{\lambda_m}{(-\beta)^{m} \displaystyle \prod_{j=1}^m \lambda_j}\sum_{k=1}^{m-1} \frac{1}{\lambda_k} & \text{for $\ell=m$}.
        \end{array}\right.
\end{align}

Since we are primarily interested in the p.d.f. of $Y_\ell^{(n)}$, for $\ell=1,m$, we use (\ref{singinteval}) and (\ref{contintsing}) to take the inverse Laplace transform of (\ref{mgfsingint}) to obtain
\begin{align}
\label{singinteval1}
    f_{Y_\ell^{(n)}}(z)=\left\{\begin{array}{cc}\displaystyle  \frac{\pi^{m(m-n-1)}2^{-m}}{\tilde{\Gamma}_m(n) (1+\theta)^{m}}\int_{\mathcal{R}}
    \prod_{k=1}^m \lambda_k^{n-m} e^{-\lambda_k} \Delta^2_m(\boldsymbol{\lambda})
    \displaystyle  \mathcal{L}^{-1}\left\{\Psi_1 (\boldsymbol{\lambda}, s)\right\}& \text{for $\ell=1$}\\
        \displaystyle\frac{\pi^{m(m-n-1)}2^{-m}}{\tilde{\Gamma}_m(n) (1+\theta)^{m}} \int_{\mathcal{R}}
    \prod_{k=1}^m \lambda_k^{n-m} e^{-\lambda_k} \Delta^2_m(\boldsymbol{\lambda})
    \displaystyle
   \mathcal{L}^{-1}\left\{\Psi_m (\boldsymbol{\lambda}, s)\right\}& \text{for $\ell=m$}
\end{array}\right.
\end{align}
Consequently, the p.d.f.s corresponding to $n>m\geq 2$ can be determined, in principle, by capitalizing on the framework developed in the previous sections along with (\ref{singinteval}), (\ref{contintsing}) and (\ref{contintsingcauchy1}). However, the ensuing algebraic complexity prevents us from obtaining general closed-form solutions. Nevertheless, for the trivial case of $m=1$, it can be shown that
\begin{align}
    \mathcal{L}^{-1}\left\{\Psi_1(\boldsymbol{\lambda},s\right\}=\frac{\mathrm{K}_{1,n}}{\Gamma(n-1)(1+\theta)} e^{\beta \lambda_1 z} (1-z)^{n-2},
\end{align}
which upon substituting into (\ref{singinteval1}) with some algebraic manipulation gives
\begin{align}
    f_{Y_1^{(n)}}(z)=\frac{(1-z)^{n-2}}{\Gamma(n-1)}
    \int_0^\infty  \lambda_1^{n-1} e^{-\lambda_1(1-\beta z)} {\rm d}\lambda_1=\frac{(n-1)(1-z)^{n-2}}{(1+\theta)(1-\beta z)^n},\; z\in(0,1).
\end{align}
Moreover, for $n-m=1$, the p.d.f. of $Y_1^{(n)}$ in (\ref{singinteval1}) can be simplified in view of 
(\ref{singinteval}), (\ref{contintsing}) and (\ref{contintsing1}) as
\begin{align}
\label{singminpdf}
    f_{Y_1^{(n)}}(z)=\frac{\pi^{m(m-1)}}{\tilde{\Gamma}_m(m)\tilde{\Gamma}_m(m+1)\left(1+\theta\right)^m\beta^{m-1}}
    \left(\mathsf{A}_m(z)+\mathsf{B}_m(z)\right)
\end{align}
where 
\begin{align}
\label{singadef}
    \mathsf{A}_m(z)=\int\limits_{\mathcal{R}} 
   e^{\beta \lambda_1 z}
   \Delta_m^2(\boldsymbol{\lambda}) \prod_{j=1}^m \lambda_j e^{-\lambda_j}
   \sum_{k=2}^{m}
   \frac{e^{(1-z)\beta \lambda_k}}{\displaystyle \prod_{\substack{i=2\\i\neq k}}^{n} \left(\lambda_k-\lambda_i\right)}\;
  {\rm d}\boldsymbol{\Lambda}
\end{align}
and
\begin{align}
\label{singbdef}
    \mathsf{B}_m(z)=(-1)^{m-1}\int\limits_{\mathcal{R}} \lambda_1
   e^{\beta \lambda_1 z}
   \Delta_m^2(\boldsymbol{\lambda}) \prod_{j=1}^m e^{-\lambda_j}
  {\rm d}\boldsymbol{\Lambda}.
\end{align}
Let us first focus on $\mathsf{A}_m(z)$. To this end, 
following the similar arguments as in Subsection II.A, we may decompose, after some algebraic manipulation, the multiple integral in (\ref{singadef}) as
\begin{align}
\label{singadecom}
    \mathsf{A}_m(z)=\frac{1}{(m-2)!} \int_0^\infty 
    x e^{-(m-\beta)x} \int_0^\infty y^2 e^{-\left[1-(1-z)\beta\right]y} 
    \mathcal{Q}_{m-2}(y,x) {\rm d}y {\rm d}x
\end{align}
where \begin{align}
    \mathcal{Q}_{m-2}(y,x)=(-1)^m \hat{K}_{m-2,1}
    (x+y)^{-1}\det\left[L^{(2)}_{m+i-3}(y)\;\;\; L^{(j)}_{m+i-j-1}(-x)\right]_{\substack{i=1,2\\
    j=2}}.
\end{align}
To facilitate further analysis, we may further simplify $ \mathcal{Q}_{m-2}(y,x)$ using the  Christoffel–Darboux formula \cite[Eq. 22.12.1]{ref:milton} as
\begin{align}
    \mathcal{Q}_{m-2}(y,x)=(-1)^m m\hat{K}_{m-2,1}
    \sum_{\ell=0}^{m-2}
    \frac{(m-2-\ell)!}{(m-\ell)!} L^{(2)}_{m-2-\ell}(-x)L^{(2)}_{m-2-\ell}(y).
\end{align}
Now we substitute $\mathcal{Q}_{m-2}(y,x)$ into (\ref{singadecom}) and solve the resultant double integral with the help of \cite[Eq. 7.414.8]{ref:gradshteyn} and the definition 3 with some algebraic manipulation to arrive at
\begin{align}
\label{singaans}
     \mathsf{A}_m(z)=\frac{m\hat{K}_{m-2,1}}{(m-2)!}
     \sum_{\ell=0}^{m-2}
     \sum_{k=0}^{m-2-\ell}
     a_{\ell,k}(m,\beta)
     \frac{(1-z)^{m-2-\ell}}{[1-(1-z)\beta]^{m-\ell+1}}
\end{align}
where 
\begin{align}
    a_{\ell,k}(m,\beta)=\frac{(-1)^\ell (m-\ell)! \beta^{m-2-\ell}}{(k+2) k! (m-2-\ell-k)! (m-\beta)^{k+2}}
\end{align}
Having evaluated $\mathsf{A}_m(z)$, we now focus $\mathsf{B}_m(z)$. As such, we may decompose the corresponding multiple integral, after some algebraic manipulation, to yield
\begin{align}
   \mathsf{B}_m(z)=\frac{(-1)^{m-1}}{(m-1)!}
   \int_0^\infty
   x e^{(m-\beta z)x} {\rm d}x
   \int_{(0,\infty)^{m-1}} \Delta^2_{m-1}(\boldsymbol{\lambda})
   \prod_{k=2}^m \lambda_k^2 e^{-\lambda_k} {\rm d}\lambda_k,
\end{align}
from which we obtain
\begin{align}
\label{singbans}
    \mathsf{B}_m(z)=\frac{(-1)^{m-1} \prod_{k=0}^{m-2} (k+1)! (k+2)!}{(m-1)!(m-\beta z)^2}.
\end{align}
Finally, we may substitute (\ref{singaans}) and (\ref{singbans}) into (\ref{singminpdf}) with some algebraic manipulation to arrive at the p.d.f. of $Y_1^{(n)}$ corresponding to $n-m=1$ as
\begin{align}
    f_{Y_1^{(n)}}(z)=\frac{m}{\beta^{m-1}} (1-\beta)^m \left(
     \sum_{\ell=0}^{m-2}
     \sum_{k=0}^{m-2-\ell}
     a_{\ell,k}(m,\beta)
     \frac{(1-z)^{m-2-\ell}}{[1-(1-z)\beta]^{m-\ell+1}}+\frac{(-1)^{m-1}}{(m-\beta z)^2}\right).
\end{align}
Similarly, following the algebraic machinery shown  in Subsection II.B, we can obtain the p.d.f. of $Y_n^{(n)}$, for $n-m=1$, as
\begin{align}
    f_{Y_n^{(n)}}(z)& =\frac{(1-\beta)^m\beta^{1-m}}{\prod_{j=1}^m (m-j)! (m-j)!} \nonumber\\
    & \qquad \times \int_0^\infty
    x^{m^2} e^{-(1-\beta z)x}
    \left(
    \mathcal{J}^{s}_m(x,z)+(-1)^{m-1}
    \det\left[ \mathcal{A}_{i,j}^{(0)}(x)\right]_{i,j=1,2,\ldots,(m-1)}\right)  {\rm d}x
\end{align}
where
\begin{align}
    \mathcal{J}^{s}_m(x,z)=\int_0^1 (1-t)^2 e^{-[1-(1-z)\beta]xt}
     \det\left[t \mathcal{A}_{i,j}^{(1)}(x)-\mathcal{A}_{i,j}^{(2)}(x)\right]_{i,j=1,2,\ldots,(m-2)}  {\rm d}t.
\end{align}
Here, an empty determinant (i.e., when $m=2$) is interpreted as unity. It is noteworthy that, although further simplification of the above integrals seems an arduous task, they can be evaluated numerically for not so large values of $m$.

Figures \ref{fig_sing_mineigvec_different_n} and \ref{fig_qq_sing_mineigvec_different_n} compare the theoretical and simulated data of $Y_1^{(n)}$ with respect to the p.d.f. and Q-Q plots for different values of $n$ with $n-m=1$ and $\theta=0.3$. A close match between the results verifies the accuracy of our formulation. Moreover, Figs. \ref{fig_sing_maxeigvec_different_n} and \ref{fig_qq_sing_maxeigvec_different_n} depict the respective plots corresponding to $Y_n^{(n)}$. Again, a close match between the analytical and simulation results validates the accuracy of our expressions.
}

\begin{figure}
    \centering
    \includegraphics[width=0.9\textwidth]{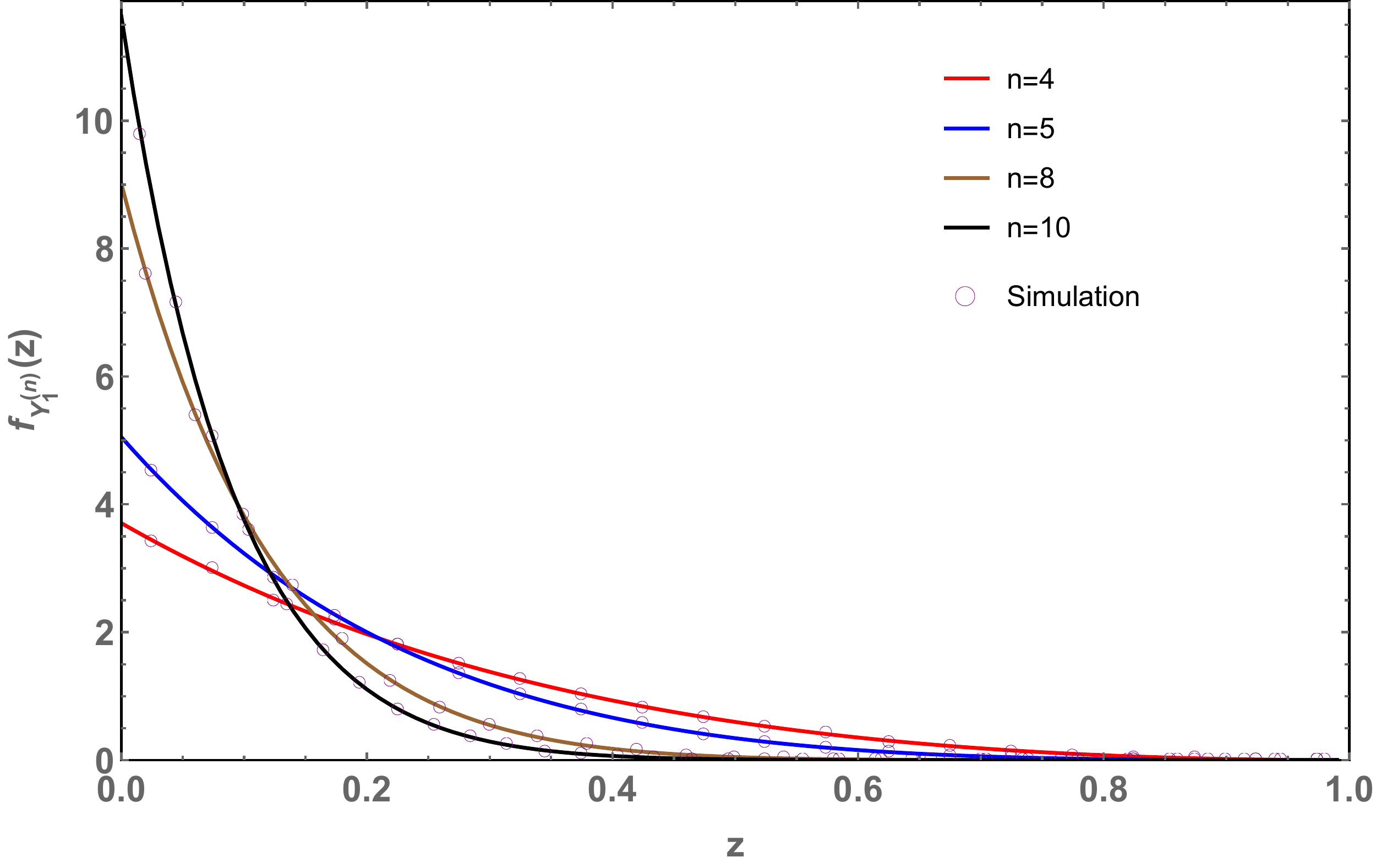}
    \caption{Comparison of simulated data points and the analytical p.d.f. $f_{Y^{(n)}_1}(z)$ for different values of $n$ with $n-m=1$ and $\theta=0.3$.}
    \label{fig_sing_mineigvec_different_n}
\end{figure}

\begin{figure}
    \centering
    \begin{subfigure}[b]{.45\textwidth}
    \includegraphics[width=\textwidth]{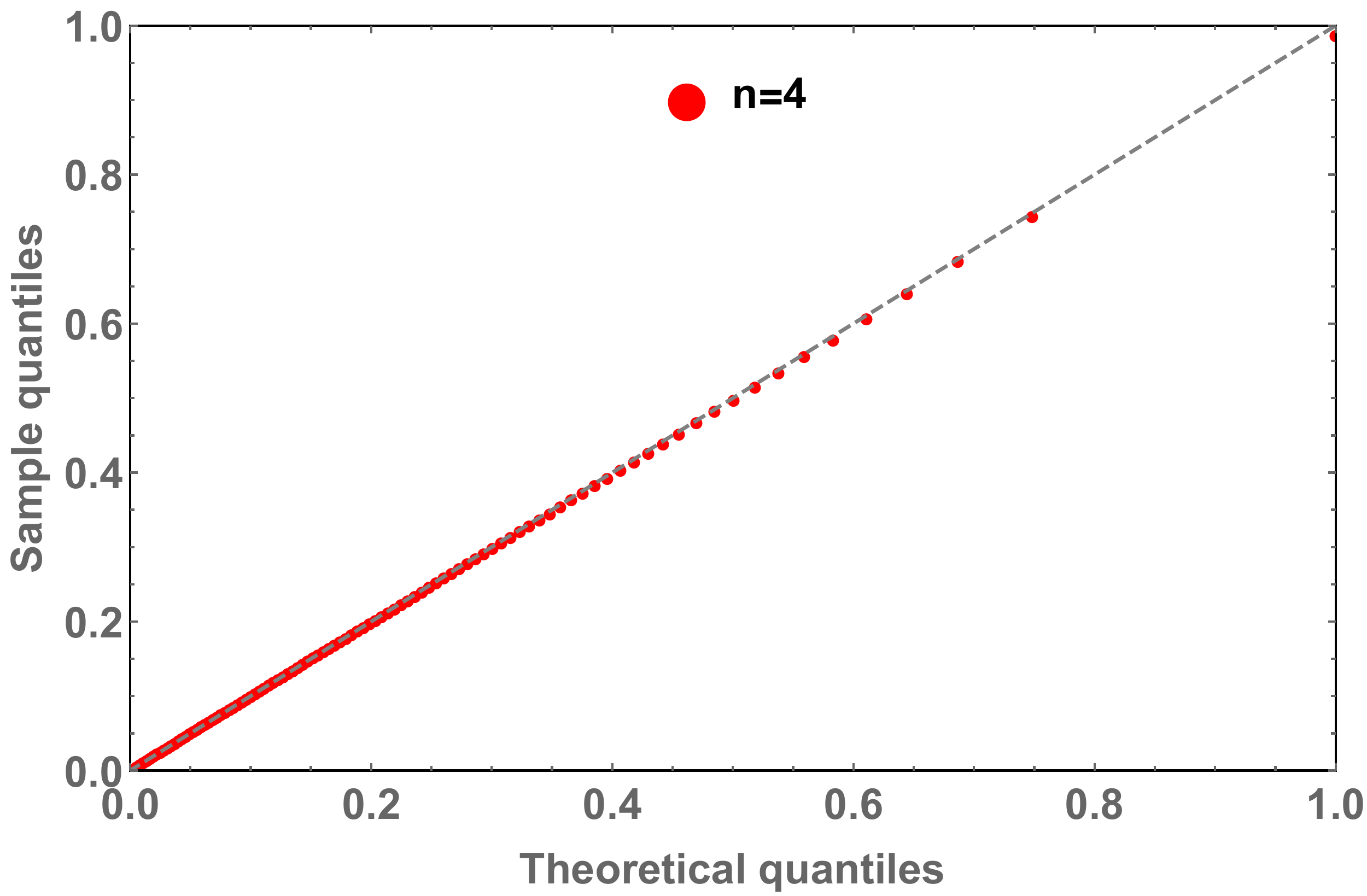}
    \end{subfigure}
    \begin{subfigure}[b]{.45\textwidth}
    \includegraphics[width=\textwidth]{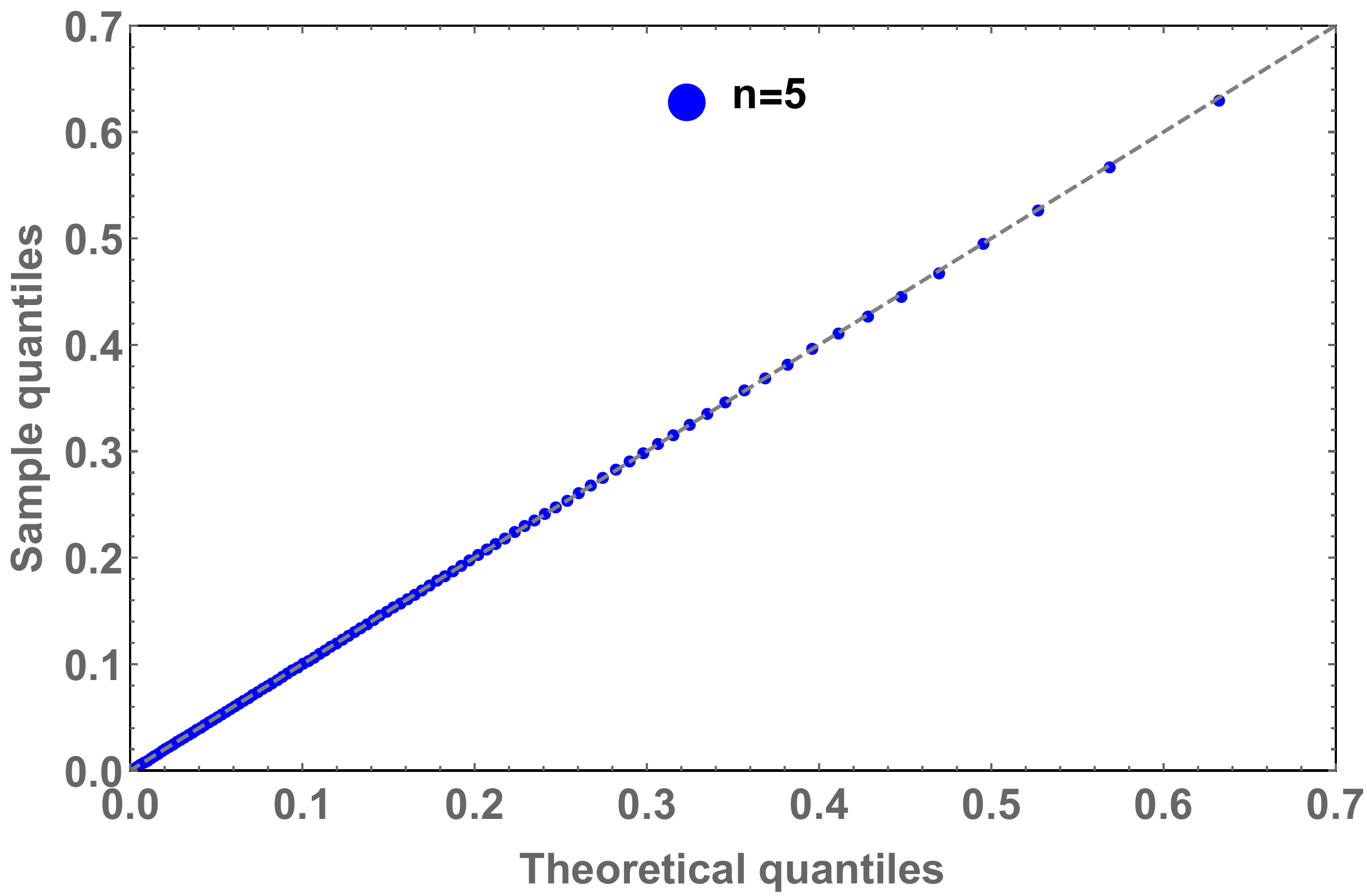}
    \end{subfigure}
    
    \begin{subfigure}[b]{.45\textwidth}
    \includegraphics[width=\textwidth]{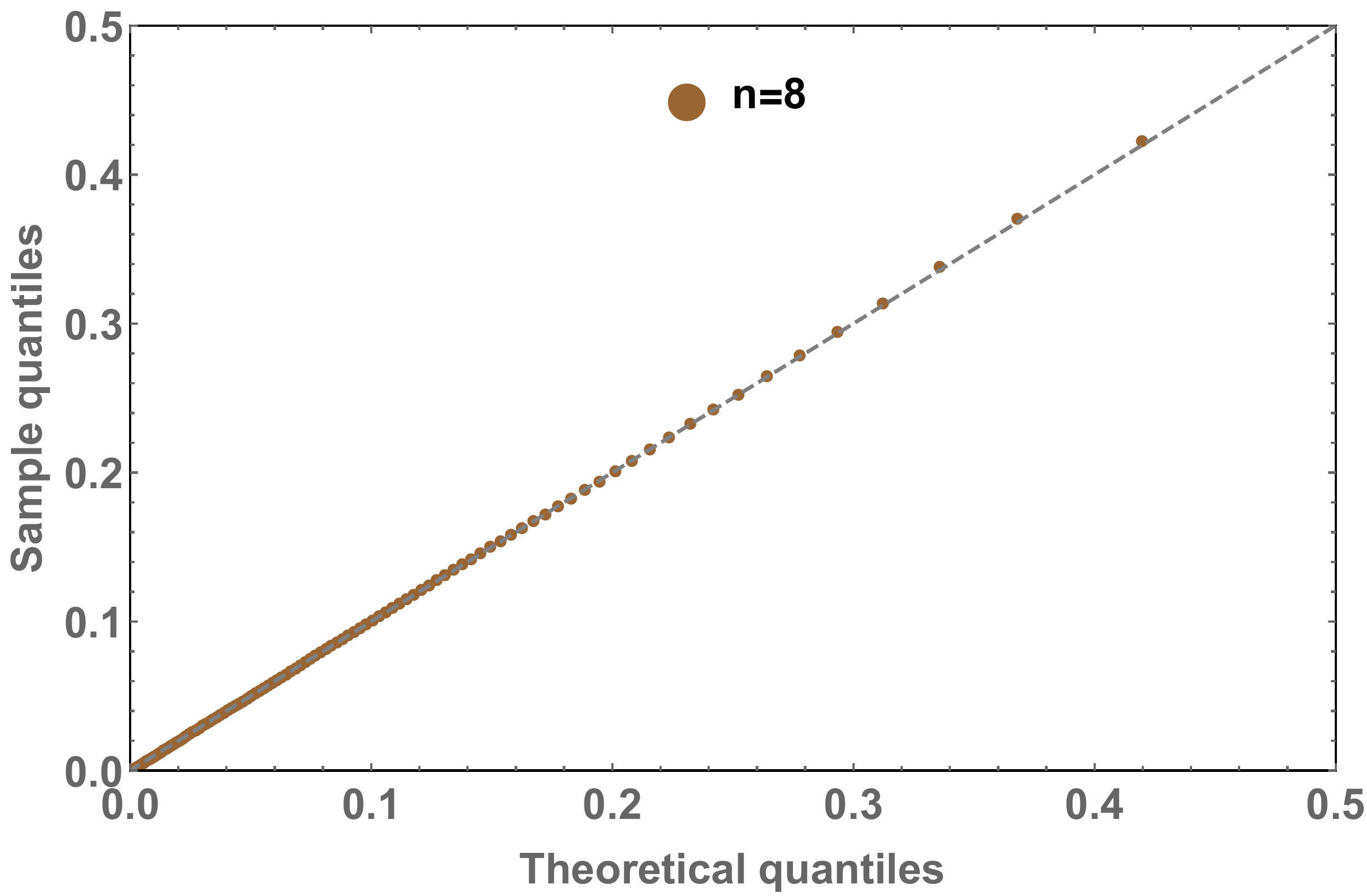}
    \end{subfigure}
    \begin{subfigure}[b]{.45\textwidth}
    \includegraphics[width=\textwidth]{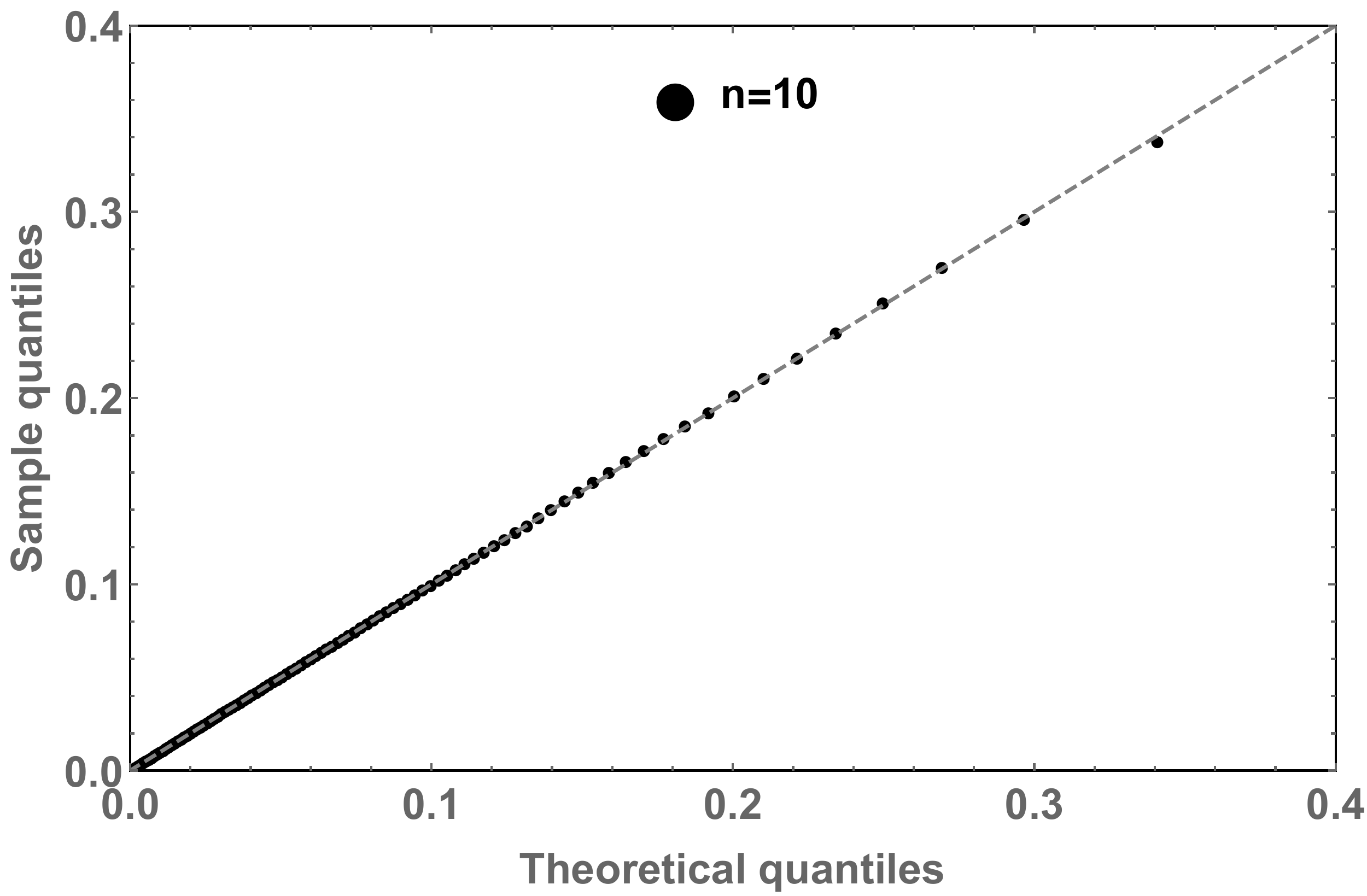}
    \end{subfigure}
    \caption{Quantile-Quantile plots of simulated data of $Y^{(n)}_1$ drawn from  $f_{Y^{(n)}_1}(z)$ for different values of $n$ with $n-m=1$ and $\theta=0.3$.}
    \label{fig_qq_sing_mineigvec_different_n}
\end{figure}

\begin{figure}
    \centering
    \includegraphics[width=0.9\textwidth]{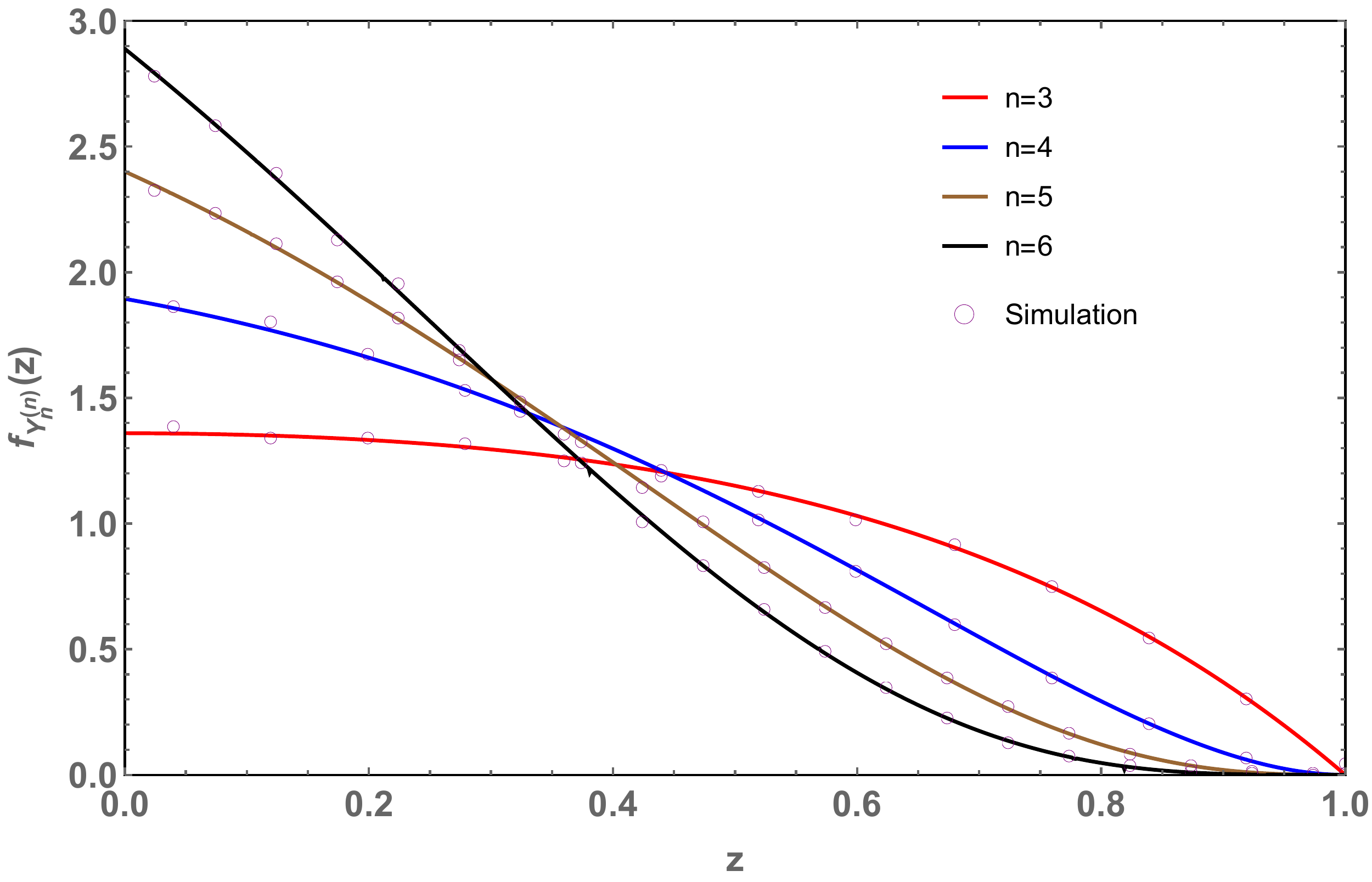}
    \caption{Comparison of simulated data points and the analytical p.d.f. $f_{Y^{(n)}_n}(z)$ for different values of $n$ with $n-m=1$ and $\theta=0.3$.}
    \label{fig_sing_maxeigvec_different_n}
\end{figure}

\begin{figure}
    \centering
    \begin{subfigure}[b]{.45\textwidth}
    \includegraphics[width=\textwidth]{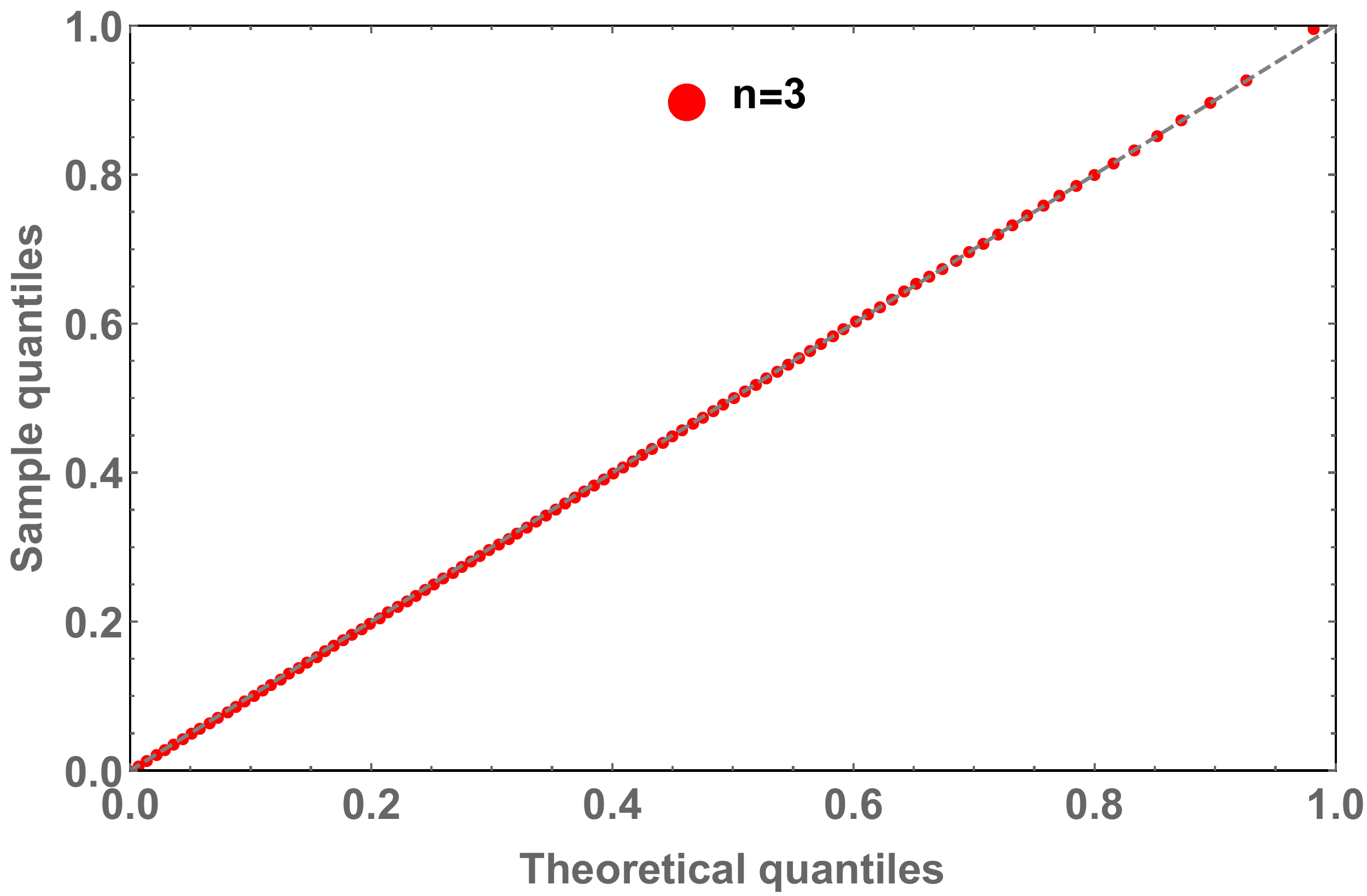}
    \end{subfigure}
    \begin{subfigure}[b]{.45\textwidth}
    \includegraphics[width=\textwidth]{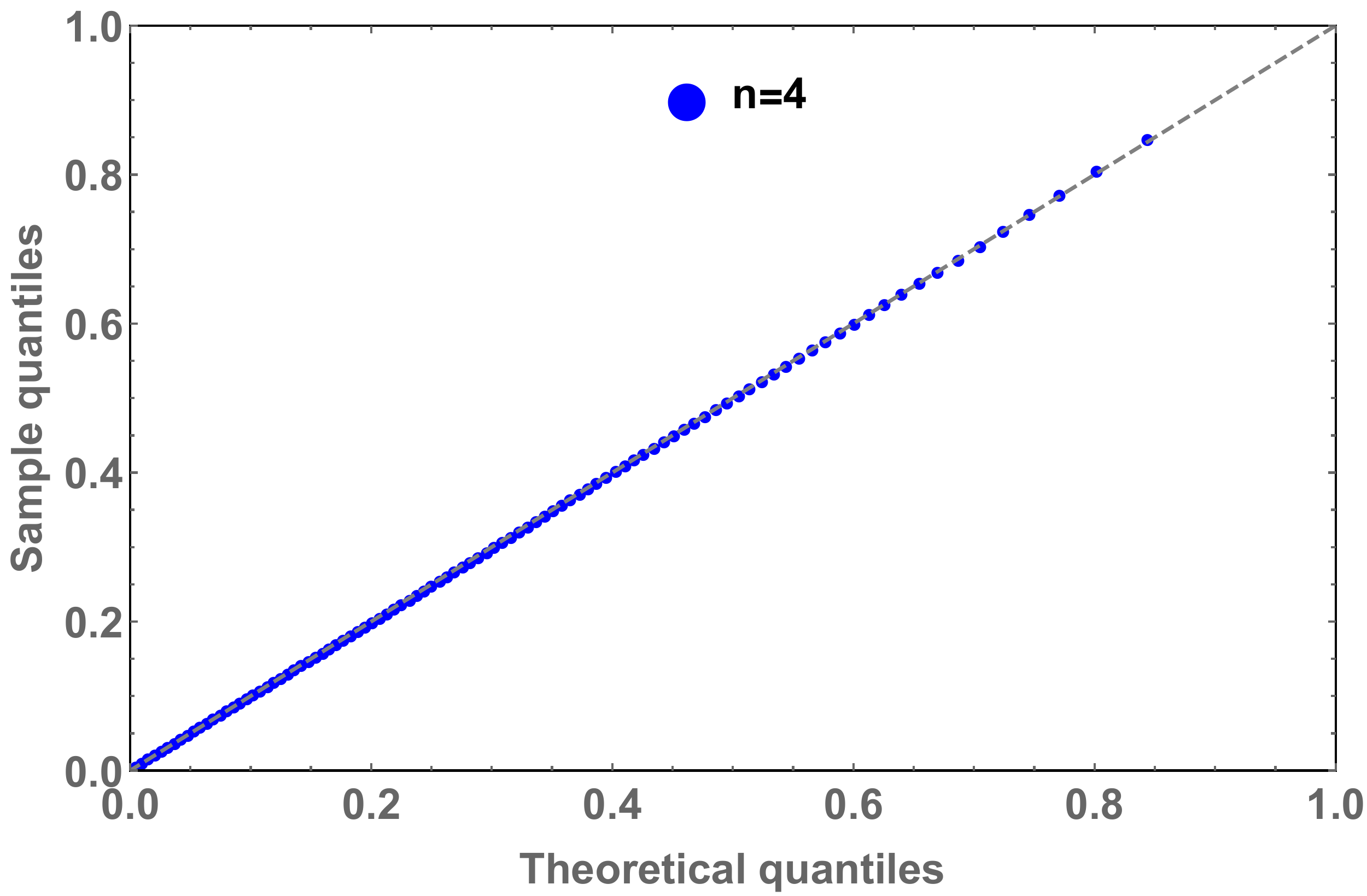}
    \end{subfigure}
    
    \begin{subfigure}[b]{.45\textwidth}
    \includegraphics[width=\textwidth]{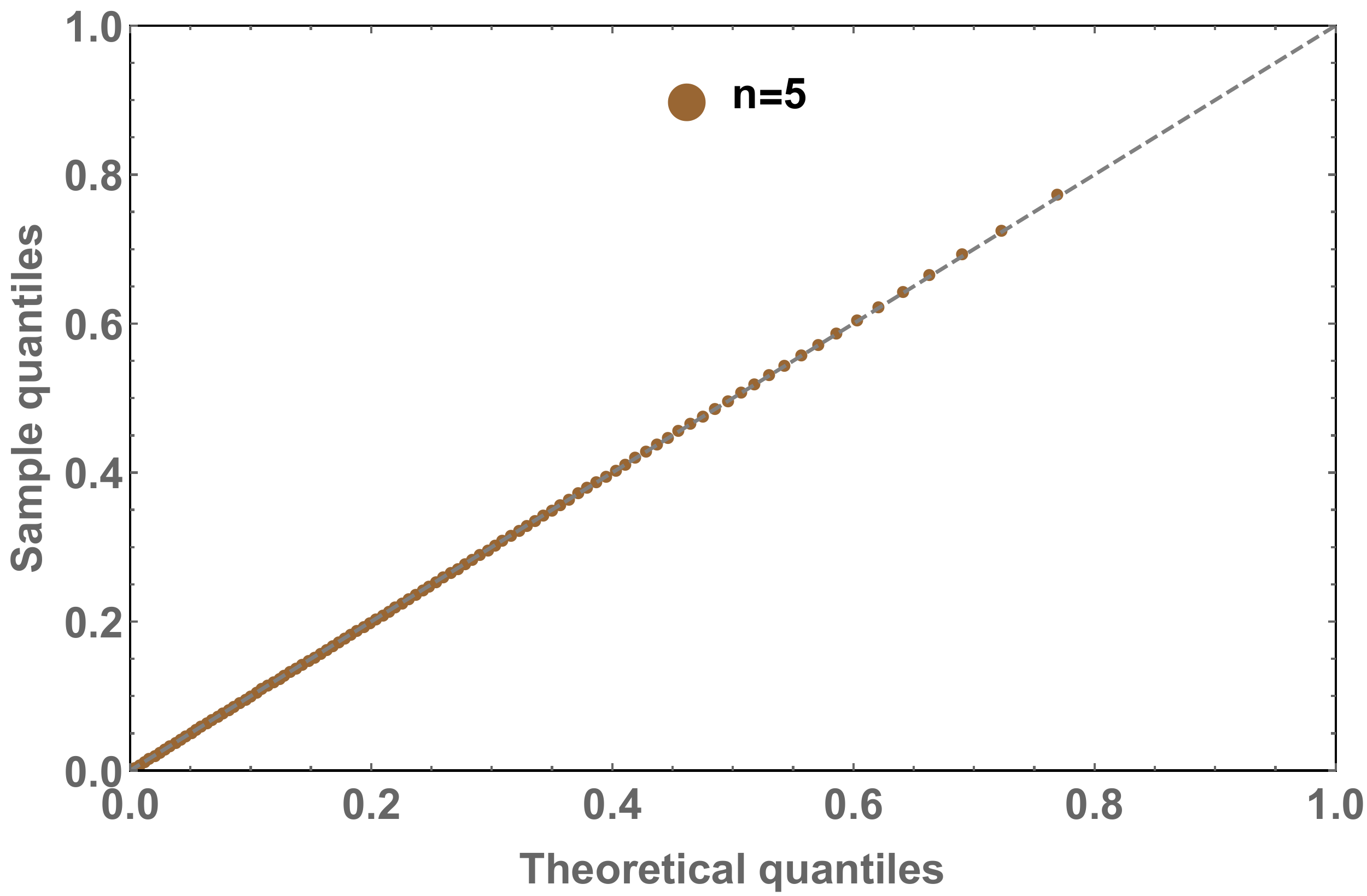}
    \end{subfigure}
    \begin{subfigure}[b]{.45\textwidth}
    \includegraphics[width=\textwidth]{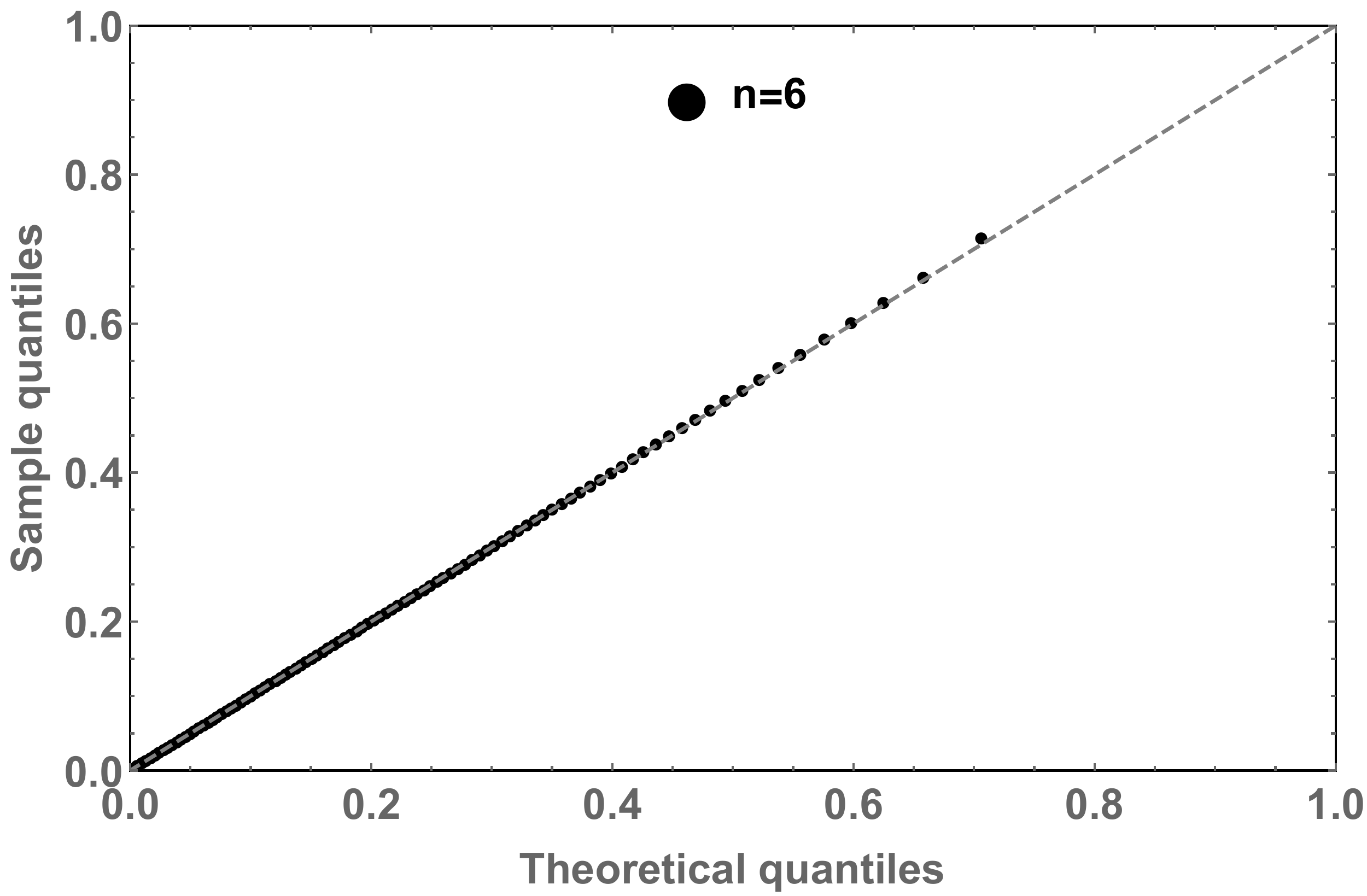}
    \end{subfigure}
    \caption{Quantile-Quantile plots of simulated data of $Y^{(n)}_n$ drawn from  $f_{Y^{(n)}_n}(z)$ for different values of $n$ with $n-m=1$ and $\theta=0.3$.}
    \label{fig_qq_sing_maxeigvec_different_n}
\end{figure}


\section{Conclusions}
This paper investigates the finite dimensional distributions of the eigenvectors corresponding to the extreme eigenvalues (i.e., the minimum and the maximum) of $\mathbf{W}\sim\mathcal{CW}_n\left(m,\mathbf{I}_n+\theta \mathbf{vv}^\dagger\right)$. In this respect, the exact p.d.f.s of $|\mathbf{v}^\dagger \mathbf{u}_1|^2$ and $|\mathbf{v}^\dagger \mathbf{u}_n|^2$ have been derived. In particular, the p.d.f. of $|\mathbf{v}^\dagger \mathbf{u}_1|^2$ assumes a closed-form involving the determinant of a square matrix of dimension $m-n$, whereas the p.d.f. of $|\mathbf{v}^\dagger \mathbf{u}_n|^2$ is expressed as a double integral with its integrand containing the determinant of a square matrix of dimension $n-2$. Our analytical p.d.f.s reveal that, in the finite dimensional setting, the least  eigenvector of $\mathbf{W}$ contains information about the dominant eigenvector (i.e., the spike $\mathbf{v}$) of the population covariance matrix. In this respect, a somewhat interesting stochastic convergence result derived here shows that, as $m,n\to\infty$ such that $m-n$ is fixed, $n |\mathbf{v}^\dagger \mathbf{u}_1|^2$ converges in distribution to $\chi^2_2/2(1+\theta)$. This further highlights the discrimination power of the least eigenvector of $\mathbf{W}$ in certain asymptotic domains. Moreover, our numerical results show that, although derived in the asymptotic setting, this analytical relationship remains valid for small values of $m$ and $n$ as well. On the other hand, the double integral form of p.d.f. derived for $|\mathbf{v}^\dagger \mathbf{u}_n|^2$ does not seem to admit a simple form except in the case of $n=2,3$, and $4$. Nevertheless, this double integral can be evaluated numerically as shown in our numerical results. 

The m.g.f. framework developed here can be easily  adapted  to derive the p.d.f.s pertaining to the other eigenmodes, however with an extra algebraic complexity, as shown in Theorem \ref{thm u2}. {\color{blue}Moreover, the same framework has been extended to derive the corresponding p.d.f.s for real and singular Wishart scenarios; however, with closed-form solutions limited to a few special configurations of $m$ and $n$. This stems from the fact that the analogous contour integrals, in general, do  not admit tractable forms, which are amenable to further analysis, except in those special cases. Be that as it may, one of the other major technical challenges is to extend the current finite dimensional analysis to rank-$r (\leq n)$ spiked Wishart matrices (i.e., $\boldsymbol{\Sigma}=\mathbf{I}_n+\sum_{k=1}^r \mathbf{v}_k\mathbf{v}_k^\dagger$ with $\left(\mathbf{v}_1\;\mathbf{v}_1\; \ldots\; \mathbf{v}_r\right)\in\mathcal{V}_{r,n}$), which remains as an open problem.}

\appendices
\section{The Evaluation of $\mathcal{Q}_{n-2}(y_1,x)$ in (\ref{Qintdef})}
\label{appendix_mehta21}
For clarity, let us consider the following integral
\begin{align}
\label{eq def R}
    \mathcal{T}_{n}(y,x)=
    \int\limits_{(0,\infty)^{n}}
    \Delta_{n}^2(\boldsymbol{z}) 
    \prod_{j=1}^{n} (y-z_j)(x-z_j)^\alpha z_j^2 e^{-z_j}
    {\rm d}z_j
\end{align}
which in turn gives
\begin{align}
\label{RQ}
    \mathcal{Q}_{n-2}(y_1,x)=(-1)^{n\alpha}\mathcal{T}_{n-2}(y_1,-x).
\end{align} Therefore, in the sequel, we evaluate $\mathcal{T}_n(y,x)$. To this end, following \cite[Eqs. 22.4.2, 22.4.11]{ref:mehta}, we focus on the related integral
\begin{align}
\label{qbeg}
& \int_{[0,\infty)^n}\prod_{j=1}^ne^{-z_j}\prod_{i=1}^{\alpha+1}(r_i-z_j)\Delta_n^2(\boldsymbol{z}) {\rm d}z_1\cdots{\rm d}z_n 
=
\prod_{i=0}^{n-1}(i+1)!(i+2)!\;\frac{\det\left[\mathfrak{P}_{n+i-1}(r_j)\right]_{i,j=1,\ldots,\alpha+1}}{\Delta_{\alpha+1}(\boldsymbol{r})},
\end{align}
where $\mathfrak{P}_{k}(x)$'s are monic polynomials orthogonal with respect to $x^2e^{-x}$, over $0\leq x<\infty$. Consequently, we choose $\mathfrak{P}_k(x)=(-1)^kk!L_k^{(2)}(x)$, which upon substituting into (\ref{qbeg}) gives
\begin{align}
\label{eq appa lag}
\int_{[0,\infty)^n}\prod_{j=1}^n z_j^2e^{-z_j}\prod_{i=1}^{\alpha+1}(r_i-z_j)&\Delta_n^2(\boldsymbol{z}) {\rm d}z_1\cdots{\rm d}z_n= \tilde{K}_{n,\alpha}\frac{\det\left[L^{(2)}_{n+i-1}(r_j)\right]_{i,j=1,\ldots,\alpha+1}}{\Delta_{\alpha+1}(\boldsymbol{r})}
\end{align}
where 
\begin{align*}
\tilde{K}_{n,\alpha}=(-1)^{(n-1)(\alpha+1)}\prod_{i=0}^{n-1}(i+1)!(i+2)!\prod_{i=1}^{\alpha+1}(-1)^i(n+i-1)!.
\end{align*}
Now we need to choose the $r_i$'s in the above relation such that the multiple integral on the left coincides with (\ref{eq def R}). In this respect, we 
 select $r_i$'s as follows:
\begin{equation*}
r_i=\left\{\begin{array}{ll}
y & \text{if $i=1$}\\
x & \text{if $i=2,\ldots,\alpha+1$.}
\end{array}\right.
\end{equation*}
However, since the $r_i$'s in (\ref{eq appa lag}) are distinct in general, the above choice of parameters drives the right side to  $0/0$ indeterminate form. To alleviate this technical difficulty, capitalizing on an approach given in \cite{ref:khatri}, instead of direct substitution, we use the limiting argument
\begin{align}
\label{qdef}
\mathcal{T
}_n(y,x)&=\tilde{K}_{n,\alpha}\lim_{r_2,\ldots,r_{\alpha+1}\to x}\frac{\det\left[L^{(2)}_{n+i-1}(y)\;\;\;L^{(2)}_{n+i-1}(r_j)\right]_{\substack{i=1,\ldots,\alpha+1\\
j=2,\ldots,\alpha+1}}}{\det[y^{i-1}\;\;\; r_j^{i-1}]_{\substack{i=1,\ldots,\alpha+1\\
j=2,\ldots,\alpha+1}}}\nonumber\\
&=\tilde{K}_{n,\alpha} \frac{\det\left[L^{(2)}_{n+i-1}(y)\;\;\;\displaystyle \frac{{\rm d}^{j-2}}{{\rm d}x^{j-2}}L^{(2)}_{n+i-1}(x)\right]_{\substack{i=1,\ldots,\alpha+1\\
j=2,\ldots,\alpha+1}}}{\det\left[y^{i-1}\;\;\; \displaystyle \frac{{\rm d}^{j-2}}{{\rm d}x^{j-2}}x^{i-1}\right]_{\substack{i=1,\ldots,\alpha+1\\
j=2,\ldots,\alpha+1}}}.
\end{align}
The denominator of (\ref{qdef}) gives
\begin{equation}
\label{denom}
\det\left[y^{i-1}\;\;\; \displaystyle \frac{{\rm d}^{j-2}}{{\rm d}x^{j-2}}x^{i-1}\right]_{\substack{i=1,\ldots,\alpha+1\\
j=2,\ldots,\alpha+1}}=\prod_{i=1}^{\alpha-1} i!\;(x-y)^\alpha.
\end{equation}
The numerator can be simplified using (\ref{lagderi}) to yield
\begin{align}
\label{num}
&\det\left[L^{(2)}_{n+i-1}(y)\;\;\;\displaystyle \frac{{\rm d}^{j-2}}{{\rm d}x^{j-2}}L^{(2)}_{n+i-1}(x)\right]_{\substack{i=1,\ldots,\alpha+1\\
j=2,\ldots,\alpha+1}}\nonumber\\
&\hspace{3cm} = 
(-1)^{\frac{1}{2}\alpha(\alpha-1)}\det\left[L^{(2)}_{n+i-1}(y)\;\;\;L^{(j)}_{n+i+1-j}(x)\right]_{\substack{i=1,\ldots,\alpha+1\\
j=2,\ldots,\alpha+1}}.
\end{align}
Therefore, we substitute (\ref{denom}) and (\ref{num}) into  (\ref{qdef}) with some algebraic manipulation to arrive at
\begin{align}
\label{Rans}
    \mathcal{T}_n(y,x)=(-1)^{n+\alpha(n+\alpha)}\frac{\hat{K}_{n,\alpha}}{(x-y)^\alpha}\det\left[L^{(2)}_{n+i-1}(y)\;\;\;L^{(j)}_{n+i+1-j}(x)\right]_{\substack{i=1,\ldots,\alpha+1\\
j=2,\ldots,\alpha+1}}
\end{align}. 
where 
\begin{align*}
  \hat{K}_{n,\alpha}= \prod_{j=1}^{\alpha+1}(n+j-1)! \prod_{j=0}^{n-1}(j+1)!(j+2)!/\prod_{j=0}^{\alpha-1} j! .
\end{align*}
Finally, we use (\ref{Rans}) in (\ref{RQ}) with some simple algebraic manipulation to obtain (\ref{eq Q solu}).


\section{Proof of Corollary \ref{cor min asy}}
\label{appendix_corr2}
 {\color{blue}Our strategy is to first show that the p.d.f. of $nZ^{(n)}_1$ converges almost everywhere to the p.d.f. of a certain scaled Chi-squared random variable. Subsequent application of the Scheff\'e's theorem \cite[Corollary 2.3]{ref:vaart} then establishes the convergence in distribution result for $nZ^{(n)}_1$.}  
 
 Let us evaluate the limiting p.d.f. of $nZ^{(n)}_1$. To this end, we write
 \begin{align}
    f_{nZ^{(n)}_1}^{(\alpha)}(z) &=
   \frac{1}{n} f_{Z^{(n)}_1}^{(\alpha)}\left(\frac{ z}{n}\right),
 \end{align}
 from which we obtain, in view of (\ref{thm main}) and by employing elementary limiting arguments,
\begin{align}
\label{asy min vec pdf}
 \lim_{n\rightarrow\infty}   f_{nZ^{(n)}_1}^{(\alpha)}(z)
    &= (1-\beta)^{\alpha-1} e^{-\frac{ z}{1-\beta}} \sum_{k_1=0}^{\infty}\sum_{k_2=0}^{\infty}\ldots \sum_{k_{\alpha}=0}^{\infty}
    \frac{\left(\alpha+\sum_{j=1}^{\alpha}k_j\right)!}{
    \prod_{j=1}^{\alpha}(j+k_j+1)! k_j!} \lim_{n\rightarrow\infty}\left\{\xi^{(\alpha)}_n(z,k_j,\beta)\right\}
\end{align}
where
\begin{equation}
    \xi^{(\alpha)}_n(z,k_j,\beta) = \frac{\prod_{j=1}^\alpha (n+\alpha-j-1)!}{(n-\beta)^{\alpha+\sum_{j=1}^\alpha k_j + 1}} \det\left[\frac{(n+\alpha)!(-\beta)^i(1-\frac{ z}{n})^i}{n(n+i-2)![1-\beta(1-\frac{ z}{n})]^i}\;\;\;\;a_{i,j}(k_j)
    \right]_{\substack{i=0,\ldots,\alpha\\j=1,\dots,\alpha}}.
\end{equation}
Now multiplying and dividing $a_{i,j}(\ell)$ by $(n+\alpha-j-\ell-1)!$ along with some algebraic manipulation gives us
\begin{align}
  \xi^{(\alpha)}_n(z,k_j,\beta)&= 
    \prod_{j=1}^\alpha \frac{(n+\alpha-j-1)!}{(n-\beta)^{k_j}(n+\alpha-j-k_j-1)!} \nonumber\\
    & \qquad \times    \det\left[\frac{(n+\alpha)!(-\beta)^i(1-\frac{ z}{n})^i}{n(n-\beta)^{\alpha+1}(n+i-2)![1-\beta(1-\frac{ z}{n})]^i}\;\;\;\;\prod_{k=0}^{\alpha-i-1}(\tilde{c}_j-k)
    \right]_{\substack{i=0,\ldots,\alpha\\j=1,\dots,\alpha}}
\end{align}
where $\tilde{c}_j=n+\alpha-j-k_j-1$ and  an empty product is interpreted as $1$. To further simplify the determinant, we apply the following row operations 
\begin{equation}
    (i+1)\text{\textsuperscript{th} row} \leftarrow (i+1)\text{\textsuperscript{th} row} +  \sum_{k=0}^{\alpha-i-1} S_{\alpha-i}^{(k)} \times (\alpha-k+1)\text{\textsuperscript{th} row }(\text{for } i=0,1,\dots,\alpha-1),
\end{equation}
where $S_n^{(m)}$ is the Stirling number of the second kind \cite{ref:erdelyi}, to yield
\begin{equation}
    \label{asymp n dependents}
    \xi^{(\alpha)}_n(z,k_j,\beta)=\prod_{j=1}^\alpha \frac{(n+\alpha-j-1)!}{(n-\beta)^{k_j}(n+\alpha-j-k_j-1)!} \det\left[\mu_n(z,i,\beta)\;\;\;\;\tilde{c}_j^{\alpha-i}
    \right]_{\substack{i=0,\ldots,\alpha\\j=1,\dots,\alpha}}
\end{equation}
in which
\begin{multline}
    \mu_n \left(z,i,\beta\right) = \frac{(n+\alpha)!(-\beta)^i(1-\frac{ z}{n})^i}{n(n-\beta)^{\alpha+1}(n+i-2)![1-\beta(1-\frac{z}{n})]^i} \\
    + \sum_{k=0}^{\alpha-i-1}S_{\alpha-i}^{(k)} \frac{(n+\alpha)!(-\beta)^{\alpha-k}(1-\frac{ z}{n})^{\alpha-k}}{n(n-\beta)^{\alpha+1}(n+\alpha-k-2)![1-\beta(1-\frac{ z}{n})]^{\alpha-k}}.
\end{multline}
To facilitate further analysis, we expand the determinant in (\ref{asymp n dependents}) with its first column to obtain
\begin{align}
\label{eq det simp vand}
    \det\left[\mu_n(z,i,\beta)\;\;\;\;\tilde{c}_j^{\alpha-i}
    \right]_{\substack{i=0,\ldots,\alpha\\j=1,\dots,\alpha}} &= \sum_{i=0}^\alpha (-1)^i \mu_n(z,i,\beta) M_i(n)
\end{align}
where $M_{i}(n)$ is the minor corresponding to the $(i+1)^{\text{th}}$ element of the first column. Since it is clear that $\mu_n(z,i,\beta)=1/n^i+o(1/n^i)$ as $n$ grows large, we need to determine the highest power of $n$ in $M_i(n)$ and corresponding coefficient to determine the limit of (\ref{eq det simp vand}) for large $n$. To this end, following the definition of Schur-polynomials $s_{\nu_i}$ \cite{ref:takemura}, for $i>0$, we obtain
\begin{align}
    s_{\nu_i}\left(\tilde{c}_1,\dots,\tilde{c}_\alpha\right) &= \frac{M_{i}(n)}{\Delta_\alpha\left(\tilde{\mathbf{c}}\right)}
\end{align}
where $\nu_i=(\underbrace{1,\dots,1}_{i \text{ terms}},\underbrace{0,\dots,0}_{\alpha-i \text{ terms}})$ and $\Delta_\alpha\left(\tilde{\mathbf{c}}\right)$ is the Vandermonde determinant in terms of $\tilde{\mathbf{c}}=\{\tilde{c}_1,\dots,\tilde{c}_\alpha\}$. Moreover, for $i=0$, it is easy to show that $M_0(n)=\Delta_\alpha (\tilde{\mathbf{c}})$. Consequently, noting that $\Delta_\alpha\left(\tilde{\mathbf{c}}\right) = \Delta_\alpha\left(\mathbf{c}\right)$ for $\mathbf{c}=\{c_1(k_1),\dots,c_\alpha(k_\alpha)\}$ with $c_j(\ell)=j+\ell$, we rewrite $M_i(n)$ as
\begin{equation}
    M_i(n) = \left\{\begin{matrix}
        \Delta_\alpha\left(\mathbf{c}\right) & i=0 \\
        s_{\nu_i}\left(\tilde{c}_1,\dots,\tilde{c}_\alpha\right) \Delta_\alpha\left(\mathbf{c}\right) & i>0
    \end{matrix}\right..
\end{equation}
Therefore, the problem boils down to determining the highest power of $n$ in the expansion of $ s_{\nu_i}\left(\tilde{c}_1,\dotsc,\tilde{c}_\alpha\right)$ which is a symmetric homogeneous polynomial of degree $i$. To this end, following 
 \cite[Eq. 5.3.5]{ref:takemura}, we may write
\begin{equation}
    s_{\nu_i}\left(\tilde{c}_1,\dotsc,\tilde{c}_\alpha\right) = \sum_{1\leq\ell_1<\dotsc<\ell_i\leq\alpha} \tilde{c}_{\ell_1}\tilde{c}_{\ell_2}\dotsc\tilde{c}_{\ell_i}
\end{equation}
where the sum includes $\binom{\alpha}{i}$ different terms of degree $i$, where $\binom{\alpha}{i}=\alpha !/i! (\alpha-i)!$ denotes the binomial coefficient. Since we are interested in the highest power of $n$, we may write $s_{\nu_i}(\tilde{c}_1,\dots,\tilde{c}_\alpha)$ as a polynomial of $n$ as
\begin{equation}
\label{eq schur n}
    h_i(n) = \binom{\alpha}{i}n^i + \text{ other lower order terms}.
\end{equation}
As such, (\ref{asymp n dependents}) assumes
\begin{equation}
   \xi^{(\alpha)}_n(z,k_j,\beta)=\prod_{j=1}^\alpha \frac{(n+\alpha-j-1)!}{(n-\beta)^{k_j}(n+\alpha-j-k_j-1)!} \sum_{i=0}^\alpha (-1)^i \mu_n(z,i,\beta) h_i(n) \Delta_\alpha(\mathbf{c}).
\end{equation}
Noting that 
\begin{equation}
    \frac{(n+\alpha-j-1)!}{(n-\beta)^{k_j}(n+\alpha-j-k_j-1)!} = 1+o\left(1\right)
\end{equation}
and $\mu_n\left(z,i,\beta\right) = \frac{1}{n^i} +  o\left(\frac{1}{n^i}\right)$, we take the limit as $n\to\infty$ with the help of (\ref{eq schur n}) to yield
\begin{align}
    \lim_{n\rightarrow\infty}\xi^{(\alpha)}_n(z,k_j,\beta) &= \sum_{i=0}^\alpha (-1)^i \frac{(-\beta)^i}{(1-\beta)^i}\binom{\alpha}{i} \Delta_\alpha\left(\mathbf{c}\right) = \frac{\Delta_\alpha(\mathbf{c})}{(1-\beta)^\alpha}.
\end{align}
Therefore, noting that $\beta=\theta/(1+\theta)$, the limiting p.d.f. (\ref{asy min vec pdf}) specializes to
\begin{align}
\label{asy pdf min}
  \lim_{n\rightarrow\infty}  f_{nZ^{(n)}_1}^{(\alpha)}(z)
    &= \mathcal{K}_\alpha (1+\theta) e^{-(1+\theta)z}
\end{align}
where
\begin{align}
    \mathcal{K}_\alpha= \sum_{k_1=0}^{\infty}\sum_{k_2=0}^{\infty}\ldots \sum_{k_{\alpha}=0}^{\infty}
    \frac{\left(\alpha+\sum_{j=1}^{\alpha}k_j\right)!}{
    \prod_{j=1}^{\alpha}(j+k_j+1)! k_j!} \Delta_{\alpha}(\mathbf{c})
\end{align}
is a constant.
The direct evaluation of  $\mathcal{K}_\alpha$ seems to be an arduous task due to the presence of the term $\left(\alpha+\sum_{j=1}^{\alpha}k_j\right)!$. To circumvent this difficulty, we replace it with an equivalent integral to yield
\begin{align}
    \mathcal{K}_\alpha= \sum_{k_1=0}^{\infty}\sum_{k_2=0}^{\infty}\ldots \sum_{k_{\alpha}=0}^{\infty}
    \frac{\int_0^\infty x^{\alpha+\sum_{j=1}^\alpha k_j}e^{-x} {\rm d}x}{
    \prod_{j=1}^{\alpha}(j+k_j+1)! k_j!} \Delta_{\alpha}(\mathbf{c}),
\end{align}
from which we obtain
\begin{align}
    \mathcal{K}_\alpha=\int_0^\infty
    x^\alpha e^{-x}
    \det\left[\sum_{k_j=0}^\infty\frac{(j+k_j)^{i-1}x^{k_j}}{(j+k_j+1)! k_j!}\right]_{i,j=1,2,\ldots,\alpha}
    {\rm d}x.
\end{align}
Now we adopt the similar procedure as in \cite[Eqs. B.21-B.23]{ref:prathapJMVA} to obtain
\begin{align}
    \mathcal{K}_\alpha= \int_0^\infty x^\alpha e^{-x} \det\left[\frac{1}{(1)_{j+\rho-i}}\sum_{k_j=0}^\infty \frac{(j+\rho)_{k_j}}{(j+\rho+1-i)_{k_j}(j+2)_{k_j}}\frac{x^{k_j}}{k_j!}\right]_{i,j=1,\dots,\alpha} {\rm d}x
\end{align}
where $\rho>0$ is an arbitrary number. Consequently, we choose $\rho=2$ to further simplify the above integral as
\begin{align}
   \mathcal{K}_\alpha&= \int_0^\infty x^\alpha e^{-x} \det\left[\frac{1}{(1)_{j+2-i}}\sum_{k_j=0}^\infty \frac{x^{k_j}}{(j+3-i)_{k_j}k_j!}\right]_{i,j=1,\dots,\alpha} {\rm d}x \nonumber\\
    &= \int_0^\infty e^{-x} \det\left[I_{j-i+2}\left(2\sqrt{x}\right)\right]_{i,j=1,\dots,\alpha} {\rm d}x
\end{align}
where $I_p(z)$ is the modified Bessel function of the first kind and order $p$. A careful inspection of the above integral and 
\cite[Eq. 74]{ref:SCNArxiv} reveals that the function $e^{-x} \det\left[I_{j-i+2}\left(2\sqrt{x}\right)\right]_{i,j=1,\dots,\alpha}$ is another p.d.f, thereby $\mathcal{K}_\alpha=1$. Therefore, (\ref{asy pdf min}) simplifies to 
\begin{equation}
    \lim_{n\rightarrow\infty}  f_{nZ^{(n)}_1}^{(\alpha)}(z)
    =(1+\theta) e^{-(1+\theta)z}
\end{equation}
which is the p.d.f. of a $\displaystyle \frac{\chi_2^2}{2(1+\theta)}$ random variable. Finally, we invoke the Scheff\'e's theorem \cite[Corollary 2.3]{ref:vaart} to conclude the proof.



\ifCLASSOPTIONcaptionsoff
 \newpage
\fi



\bibliographystyle{IEEEtran}
\bibliography{./references.bib}

\begin{thebibliography}{100}
\providecommand{\url}[1]{#1}
\csname url@samestyle\endcsname
\providecommand{\newblock}{\relax}
\providecommand{\bibinfo}[2]{#2}
\providecommand{\BIBentrySTDinterwordspacing}{\spaceskip=0pt\relax}
\providecommand{\BIBentryALTinterwordstretchfactor}{4}
\providecommand{\BIBentryALTinterwordspacing}{\spaceskip=\fontdimen2\font plus
\BIBentryALTinterwordstretchfactor\fontdimen3\font minus
  \fontdimen4\font\relax}
\providecommand{\BIBforeignlanguage}[2]{{%
\expandafter\ifx\csname l@#1\endcsname\relax
\typeout{** WARNING: IEEEtran.bst: No hyphenation pattern has been}%
\typeout{** loaded for the language `#1'. Using the pattern for}%
\typeout{** the default language instead.}%
\else
\language=\csname l@#1\endcsname
\fi
#2}}
\providecommand{\BIBdecl}{\relax}
\BIBdecl

\bibitem{ref:paul}
D.~Paul, ``Asymptotics of sample eigenstructure for a large dimensional spiked
  covariance model,'' \emph{Stat. Sin.}, vol.~17, no.~4, pp. 1617--1642, 2007.

\bibitem{ref:bGeorges}
F.~Benaych-Georges and R.~R. Nadakuditi, ``The eigenvalues and eigenvectors of
  finite, low rank perturbations of large random matrices,'' \emph{Adv. Math.},
  vol. 227, no.~1, pp. 494--521, 2011.

\bibitem{ref:bloemendal}
A.~Bloemendal, A.~Knowles, H.~T. Yau, and J.~Yin, ``On the principal components
  of sample covariance matrices,'' \emph{Probab. Theory Relat. Fields}, vol.
  164, no. 1-2, pp. 459--552, 2016.

\bibitem{ref:haokai}
H.~Xi, F.~Yang, and J.~Yin, ``{Convergence of eigenvector empirical spectral
  distribution of sample covariance matrices},'' \emph{Ann. Stat.}, vol.~48,
  no.~2, pp. 953--982, 2020.

\bibitem{ref:bao}
Z.~{Bao}, X.~{Ding}, J.~{Wang}, and K.~{Wang}, ``{Statistical inference for
  principal components of spiked covariance matrices},'' \emph{arXiv:2008.11903
  [math.ST]}, 2020.

\bibitem{ref:wWang}
W.~Wang and J.~Fan, ``\BIBforeignlanguage{English (US)}{Asymptotics of
  empirical eigenstructure for high dimensional spiked covariance},''
  \emph{\BIBforeignlanguage{English (US)}{Ann. Stat.}}, vol.~45, no.~3, pp.
  1342--1374, 2017.

\bibitem{ref:vantrees}
H.~L. Van~Trees, \emph{Optimum Array Processing: Part {IV} of Detection,
  Estimation, and Modulation Theory}.\hskip 1em plus 0.5em minus 0.4em\relax
  John Wiley \& Sons, 2004.

\bibitem{ref:harryLee}
H.~Lee and F.~Li, ``An eigenvector technique for detecting the number of
  emitters in a cluster,'' \emph{IEEE Trans. Signal Process.}, vol.~42, no.~9,
  pp. 2380--2388, 1994.

\bibitem{ref:denis}
D.~Igambi, X.~Yang, and B.~Jalal, ``Robust adaptive beamforming based on
  desired signal power reduction and output power of spatial matched filter,''
  \emph{IEEE Access}, vol.~6, pp. 50\,217--50\,228, 2018.

\bibitem{ref:fengFeng}
F.~Chen, F.~Shen, and J.~Song, ``Robust adaptive beamforming using
  low-complexity correlation coefficient calculation algorithms,''
  \emph{Electron. Lett.}, vol.~51, no.~6, pp. 443--445, 2015.

\bibitem{ref:luYan}
L.~Yan, S.~Piao, F.~Xu, and J.~Yang, ``Improved {OP} approach utilising
  correlated projection and eigenspace processing for robust adaptive
  beamforming,'' \emph{Electron. Lett.}, vol.~55, pp. 1170--1172, 2019.

\bibitem{ref:johnPaul}
I.~M. Johnstone and D.~Paul, ``{PCA} in high dimensions: {A}n orientation,''
  \emph{Proc. IEEE}, vol. 106, no.~8, pp. 1277--1292, 2018.

\bibitem{ref:arthurLu}
I.~M. Johnstone and A.~Y. Lu, ``On consistency and sparsity for principal
  components analysis in high dimensions,'' \emph{J. Am. Stat. Assoc.}, vol.
  104, no. 486, pp. 682--693, 2009.

\bibitem{ref:huiZou}
H.~Zou and L.~Xue, ``A selective overview of sparse principal component
  analysis,'' \emph{Proc. IEEE}, vol. 106, no.~8, pp. 1311--1320, 2018.

\bibitem{ref:liaoMahoney}
Z.~Liao, R.~Couillet, and M.~Mahoney, ``Sparse quantized spectral clustering,''
  in \emph{9th Int. Conf. Learn. Represent. ({ICLR})}, Virtual Only, France,
  2021.

\bibitem{ref:johnsonMestre}
B.~Johnson, Y.~Abramovich, and X.~Mestre, ``{MUSIC}, {G-MUSIC}, and
  maximum-likelihood performance breakdown,'' \emph{IEEE Trans. Signal
  Process.}, vol.~56, no.~8, pp. 3944--3958, 2008.

\bibitem{ref:yangEdgar}
F.~Yang, S.~Liu, E.~Dobriban, and D.~P. Woodruff, ``How to reduce dimension
  with {PCA} and random projections?'' \emph{IEEE Trans. Inf. Theory}, vol.~67,
  no.~12, pp. 8154--8189, 2021.

\bibitem{ref:fumito}
F.~Tagashira, T.~Obuchi, and T.~Tanaka, ``Sharp asymptotics of matrix sketching
  for a rank-one spiked model,'' in \emph{IEEE Int. Symp. Inf. Theory
  ({ISIT})}.\hskip 1em plus 0.5em minus 0.4em\relax IEEE, 2021, pp. 250--255.

\bibitem{ref:donohoGavishJohn}
D.~L. Donoho, M.~Gavish, and I.~M. Johnstone, ``Optimal shrinkage of
  eigenvalues in the spiked covariance model,'' \emph{Ann. Stat.}, vol.~46,
  no.~4, pp. 1742--1778, 2018.

\bibitem{ref:donohoBehrooz}
D.~L. {Donoho} and B.~{Ghorbani}, ``{Optimal Covariance Estimation for
  Condition Number Loss in the Spiked Model},'' \emph{arXiv:1810.07403
  [math.ST]}, 2018.

\bibitem{ref:remiMonasson}
R.~Monasson and D.~Villamaina, ``Estimating the principal components of
  correlation matrices from all their empirical eigenvectors,'' \emph{{EPL}},
  vol. 112, no.~5, p. 50001, 2015.

\bibitem{ref:bunPotters}
J.~Bun, R.~Allez, J.-P. Bouchaud, and M.~Potters, ``Rotational invariant
  estimator for general noisy matrices,'' \emph{IEEE Trans. Inf. Theory},
  vol.~62, no.~12, pp. 7475--7490, 2016.

\bibitem{ref:donohoeGavish}
M.~Gavish and D.~L. Donoho, ``Optimal shrinkage of singular values,''
  \emph{IEEE Trans. Inf. Theory}, vol.~63, no.~4, pp. 2137--2152, 2017.

\bibitem{ref:huiJiaLi}
H.-J. Li, Z.~Wang, J.~Pei, J.~Cao, and Y.~Shi, ``Optimal estimation of low-rank
  factors via feature level data fusion of multiplex signal systems,''
  \emph{IEEE Trans. Knowl. Data Eng.}, pp. 1--12, 2020.

\bibitem{ref:rish}
R.~Dudeja, M.~Bakhshizadeh, J.~Ma, and A.~Maleki, ``Analysis of spectral
  methods for phase retrieval with random orthogonal matrices,'' \emph{IEEE
  Trans. Inf. Theory}, vol.~66, no.~8, pp. 5182--5203, 2020.

\bibitem{ref:rCoulliet}
R.~Couillet and W.~Hachem, ``Fluctuations of spiked random matrix models and
  failure diagnosis in sensor networks,'' \emph{IEEE Trans. Inf. Theory},
  vol.~59, no.~1, pp. 509--525, 2012.

\bibitem{ref:oliver}
O.~Ledoit and S.~P{\'e}ch{\'e}, ``Eigenvectors of some large sample covariance
  matrix ensembles,'' \emph{Probab. Theory Relat. Fields}, vol. 151, no.~1, pp.
  233--264, 2011.

\bibitem{ref:nadler}
B.~Nadler, ``{Finite sample approximation results for principal component
  analysis: A matrix perturbation approach},'' \emph{Ann. Stat.}, vol.~36,
  no.~6, pp. 2791--2817, 2008.

\bibitem{ref:nobel}
A.~A. Shabalin and A.~B. Nobel, ``Reconstruction of a low-rank matrix in the
  presence of {G}aussian noise,'' \emph{J. Multivar. Anal.}, vol. 118, pp.
  67--76, 2013.

\bibitem{ref:rajPCInfo}
R.~R. Nadakuditi, ``When are the most informative components for inference also
  the principal components?'' \emph{arXiv:1302.1232 [math.ST]}, 2013.

\bibitem{ref:leeb}
W.~Leeb and E.~Romanov, ``Optimal spectral shrinkage and {PCA} with
  heteroscedastic noise,'' \emph{IEEE Trans. Inf. Theory}, vol.~67, no.~5, pp.
  3009--3037, 2021.

\bibitem{ref:rajDenoise}
R.~R. Nadakuditi, ``Optshrink: {A}n algorithm for improved low-rank signal
  matrix denoising by optimal, data-driven singular value shrinkage,''
  \emph{IEEE Trans. Inf. Theory}, vol.~60, no.~5, pp. 3002--3018, 2014.

\bibitem{ref:gavishSing}
M.~Gavish and D.~L. Donoho, ``The optimal hard threshold for singular values is
  $4/\sqrt{3} $,'' \emph{IEEE Trans. Inf. Theory}, vol.~60, no.~8, pp.
  5040--5053, 2014.

\bibitem{ref:baoWang}
Z.~Bao, X.~Ding, , and K.~Wang, ``Singular vector and singular subspace
  distribution for the matrix denoising model,'' \emph{Ann. Stat.}, vol.~49,
  no.~1, pp. 370--392, 2021.

\bibitem{ref:leebDeno}
W.~E. Leeb, ``Matrix denoising for weighted loss functions and heterogeneous
  signals,'' \emph{SIAM J. Math. Data Sci.}, vol.~3, no.~3, pp. 987--1012,
  2021.

\bibitem{ref:wolf}
O.~Ledoit and M.~Wolf, ``Shrinkage estimation of large covariance matrices:
  {K}eep it simple, statistician?'' \emph{J. Multivar. Anal.}, vol. 186, p.
  104796, 2021.

\bibitem{ref:tuladhar}
S.~R. Tuladhar, J.~R. Buck, and K.~E. Wage, ``Random matrix theory model for
  mean notch depth of the diagonally loaded minimum variance distortionless
  response beamformer for a single interferer case,'' \emph{J. Acoust. Soc.
  Am.}, vol. 135, no.~4, pp. 2359--2359, 2014.

\bibitem{ref:wage}
K.~E. Wage and J.~R. Buck, ``\BIBforeignlanguage{Undetermined}{Snapshot
  performance of the dominant mode rejection beamformer},''
  \emph{\BIBforeignlanguage{Undetermined}{IEEE J. Oceanic Eng.}}, vol.~39,
  no.~2, pp. 212--225, 2014.

\bibitem{ref:sandeep}
S.~Gogineni, P.~Setlur, M.~Rangaswamy, and R.~R. Nadakuditi, ``Passive radar
  detection with noisy reference channel using principal subspace similarity,''
  \emph{IEEE Trans. Aerosp. Electron. Syst.}, vol.~54, no.~1, pp. 18--36, 2018.

\bibitem{ref:gog1}
------, ``Passive radar detection with noisy reference signal using measured
  data,'' in \emph{Proc. IEEE Radar Conf.}, Seattle, WA, USA, 2017, pp.
  858--861.

\bibitem{ref:gog2}
------, ``Comparison of passive radar detectors with noisy reference signal,''
  in \emph{Proc. IEEE Stat. Signal Process. Workshop}, Palma de Mallorca,
  Spain, 2016, pp. 1--5.

\bibitem{ref:gog3}
------, ``Random matrix theory inspired passive bistatic radar detection of
  low-rank signals,'' in \emph{Proc. IEEE Int. Radar Conf.}, Washington, DC,
  2015, pp. 1656--1659.

\bibitem{ref:rokhlin}
V.~Rokhlin, A.~Szlam, and M.~Tygert, ``A randomized algorithm for principal
  component analysis,'' \emph{SIAM J. Matrix Anal. Appl.}, vol.~31, no.~3, pp.
  1100--1124, 2010.

\bibitem{ref:halkoa}
N.~Halko, P.-G. Martinsson, Y.~Shkolnisky, and M.~Tygert, ``An algorithm for
  the principal component analysis of large data sets,'' \emph{SIAM J. Sci.
  Comput.}, vol.~33, no.~5, pp. 2580--2594, 2011.

\bibitem{ref:halkob}
N.~Halko, P.~G. Martinsson, and J.~A. Tropp, ``Finding structure with
  randomness: {P}robabilistic algorithms for constructing approximate matrix
  decompositions,'' \emph{SIAM Rev. Soc. Ind. Appl. Math.}, vol.~53, no.~2, pp.
  217--288, 2011.

\bibitem{ref:raskutti}
G.~Raskutti and M.~W. Mahoney, ``A statistical perspective on randomized
  sketching for ordinary least-squares,'' \emph{J. Mach. Learn. Res.}, vol.~17,
  no. 213, pp. 1--31, 2016.

\bibitem{ref:liuDobriban}
S.~Liu and E.~Dobriban, ``Ridge regression: {S}tructure, cross-validation, and
  sketching,'' in \emph{8th Int. Conf. Learn. Represent. ({ICLR})}, Addis
  Ababa, Ethiopia, 2020.

\bibitem{ref:drineas}
P.~Drineas, M.~W. Mahoney, and S.~Muthukrishnan, ``Sampling algorithms for
  $\ell_2$ regression and applications,'' in \emph{Proc. Annu. ACM-SIAM Symp.
  Discrete Algorithms}, ser. SODA '06.\hskip 1em plus 0.5em minus 0.4em\relax
  USA: Society for Industrial and Applied Mathematics, 2006, pp. 1127--1136.

\bibitem{ref:wangGittens}
S.~Wang, A.~Gittens, and M.~W. Mahoney, ``Sketched ridge regression:
  {O}ptimization perspective, statistical perspective, and model averaging,''
  in \emph{34th {ICML}}, ser. Proceedings of Machine Learning Research,
  vol.~70.\hskip 1em plus 0.5em minus 0.4em\relax PMLR, 2017, pp. 3608--3616.

\bibitem{ref:cannings}
T.~I. Cannings and R.~J. Samworth, ``Random-projection ensemble
  classification,'' \emph{J. R. Stat. Soc. Series B Stat. Methodol.}, vol.~79,
  no.~4, pp. 959--1035, 2017.

\bibitem{ref:pilancia}
M.~Pilanci and M.~J. Wainwright, ``Randomized sketches of convex programs with
  sharp guarantees,'' \emph{IEEE Trans. Inf. Theory}, vol.~61, no.~9, pp.
  5096--5115, 2015.

\bibitem{ref:pilancib}
------, ``Iterative hessian sketch: {F}ast and accurate solution approximation
  for constrained least-squares,'' \emph{J. Mach. Learn. Res.}, vol.~17,
  no.~53, pp. 1--38, 2016.

\bibitem{ref:pilancic}
------, ``Newton sketch: {A} near linear-time optimization algorithm with
  linear-quadratic convergence,'' \emph{SIAM J. Optim.}, vol.~27, no.~1, pp.
  205--245, 2017.

\bibitem{ref:andersonIntroMultivariate}
T.~W. Anderson, \emph{An Introduction to Multivariate Statistical
  Analysis}.\hskip 1em plus 0.5em minus 0.4em\relax John Wiley \& Sons, 1958.

\bibitem{ref:muirhead}
R.~J. Muirhead, \emph{Aspects of Multivariate Statistical Theory}.\hskip 1em
  plus 0.5em minus 0.4em\relax John Wiley \& Sons, 2009, vol. 197.

\bibitem{ref:johnstoneStatisticalChallenges}
I.~M. Johnstone and D.~M. Titterington, ``Statistical challenges of
  high-dimensional data,'' \emph{Philos. Trans. R. Soc. London, Ser. A}, vol.
  367, no. 1906, pp. 4237--4253, 2009.

\bibitem{ref:anderson}
G.~W. Anderson, A.~Guionnet, and O.~Zeitouni, \emph{An Introduction to Random
  Matrices}.\hskip 1em plus 0.5em minus 0.4em\relax Cambridge University Press,
  2010.

\bibitem{ref:bai}
Z.~Bai and J.~W. Silverstein, \emph{Spectral Analysis of Large Dimensional
  Random Matrices}.\hskip 1em plus 0.5em minus 0.4em\relax Springer, 2010,
  vol.~20.

\bibitem{ref:edelman}
A.~Edelman, ``Eigenvalues and condition numbers of random matrices,''
  \emph{SIAM J. Matrix Anal. Appl.}, vol.~9, no.~4, pp. 543--560, 1988.

\bibitem{ref:jack1}
J.~W. Silverstein, ``Some limit theorems on the eigenvectors of large
  dimensional sample covariance matrices,'' \emph{J. Multivar. Anal.}, vol.~15,
  no.~3, pp. 295--324, 1984.

\bibitem{ref:jack2}
------, ``On the eigenvectors of large dimensional sample covariance
  matrices,'' \emph{J. Multivar. Anal.}, vol.~30, no.~1, pp. 1--16, 1989.

\bibitem{ref:jack3}
------, ``Weak convergence of random functions defined by the eigenvectors of
  sample covariance matrices,'' \emph{Ann. Probab.}, vol.~18, no.~3, pp.
  1174--1194, 1990.

\bibitem{ref:zanella}
A.~Zanella and M.~Chiani, ``{On the distribution of the
  $\ell\textsuperscript{th}$ largest eigenvalue of spiked complex {W}ishart
  matrices},'' \emph{Acta Phys. Pol. B}, vol.~51, no.~7, pp. 1687--1705, 2020.

\bibitem{ref:mehta}
M.~L. Mehta, \emph{Random Matrices}, 3rd~ed.\hskip 1em plus 0.5em minus
  0.4em\relax New York: Academic Press, 2004.

\bibitem{ref:forresterLogGases}
P.~J. Forrester, \emph{Log-Gases and Random Matrices (LMS-34)}.\hskip 1em plus
  0.5em minus 0.4em\relax Princeton, NJ: Princeton University Press, 2010.

\bibitem{ref:potters}
M.~Potters and J.~P. Bouchaud, \emph{A First Course in Random Matrix Theory:
  For Physicists, Engineers and Data Scientists}.\hskip 1em plus 0.5em minus
  0.4em\relax Cambridge University Press, 2020.

\bibitem{ref:akemannOxfordHandbook}
G.~Akemann, J.~Baik, and P.~Di~Francesco, \emph{The Oxford Handbook of Random
  Matrix Theory}.\hskip 1em plus 0.5em minus 0.4em\relax Oxford University
  Press, 2011.

\bibitem{ref:couillet}
R.~Couillet and M.~Debbah, \emph{Random Matrix Methods for Wireless
  Communications}.\hskip 1em plus 0.5em minus 0.4em\relax Cambridge University
  Press, 2011.

\bibitem{ref:telatar}
E.~Telatar, ``Capacity of multi-antenna {Gaussian} channels,'' \emph{Eur.
  Trans. Telecommun.}, vol.~10, no.~6, pp. 585--595, 1999.

\bibitem{ref:tulino}
A.~Tulino and S.~Verd{\'u}, ``\BIBforeignlanguage{English (US)}{Random matrix
  theory and wireless communications},'' \emph{\BIBforeignlanguage{English
  (US)}{Found. Trends Commun. Inf. Theory}}, vol.~1, no.~1, pp. 1--182, 2004.

\bibitem{ref:asendorf}
N.~Asendorf and R.~R. Nadakuditi, ``Improved detection of correlated signals in
  low-rank-plus-noise type data sets using informative canonical correlation
  analysis ({ICCA}),'' \emph{IEEE Trans. Inf. Theor.}, vol.~63, no.~6, pp.
  3451--3467, 2017.

\bibitem{ref:fan}
J.~Fan, Y.~Liao, and M.~Mincheva, ``High-dimensional covariance matrix
  estimation in approximate factor models,'' \emph{Ann. Stat.}, vol.~39, no.~6,
  pp. 3320--3356, 2011.

\bibitem{ref:ontaski}
A.~Onatski, ``Testing hypotheses about the number of factors in large factor
  models,'' \emph{Econometrica}, vol.~77, no.~5, pp. 1447--1479, 2009.

\bibitem{ref:ke}
Z.~T. Ke, Y.~Ma, and X.~Lin, ``Estimation of the number of spiked eigenvalues
  in a covariance matrix by bulk eigenvalue matching analysis,'' \emph{J. Am.
  Stat. Assoc.}, pp. 1--19, 2021.

\bibitem{ref:onatskiSignal}
A.~Onatski, M.~J. Moreira, and M.~Hallin, ``{Signal detection in high
  dimension: The multispiked case},'' \emph{Ann. Stat.}, vol.~42, no.~1, pp.
  225--254, 2014.

\bibitem{ref:johnstone}
I.~M. Johnstone, ``On the distribution of the largest eigenvalue in principal
  components analysis,'' \emph{Ann. Stat.}, vol.~29, no.~2, pp. 295--327, 2001.

\bibitem{ref:baikPhaseTrans}
J.~Baik, G.~B. Arous, and S.~Péché, ``{Phase transition of the largest
  eigenvalue for nonnull complex sample covariance matrices},'' \emph{Ann.
  Probab.}, vol.~33, no.~5, pp. 1643--1697, 2005.

\bibitem{ref:tracyLevelSp}
C.~A. Tracy, ``Level-spacing distributions and the {A}iry kernel,''
  \emph{Commun. Math. Phys.}, vol. 159, pp. 151--174, 1994.

\bibitem{ref:tracySymplectic}
------, ``On orthogonal and sympletic matrix ensembles,'' \emph{Commun. Math.
  Phys.}, vol. 177, pp. 727--754, 1996.

\bibitem{ref:marchenko}
V.~A. Marchenko and L.~A. Pastur, ``Distribution of eigenvalues for some sets
  of random matrices,'' \emph{Sb. Math.}, vol. 114, no.~4, pp. 507--536, 1967.

\bibitem{ref:thomas}
J.~Thomas, L.~Scharf, and D.~Tufts, ``The probability of a subspace swap in the
  {SVD},'' \emph{IEEE Trans. Signal Process.}, vol.~43, no.~3, pp. 730--736,
  1995.

\bibitem{ref:tuft}
D.~Tufts, A.~Kot, and R.~Vaccaro, ``The threshold analysis of {SVD}-based
  algorithms,'' in \emph{Proc. IEEE ICASSP}.\hskip 1em plus 0.5em minus
  0.4em\relax Los Alamitos, CA, USA: IEEE Computer Society, 1988, pp.
  2416--2419.

\bibitem{ref:paulRandomMtrx}
D.~Paul and A.~Aue, ``Random matrix theory in statistics: A review,'' \emph{J.
  Stat. Plann. Inference}, vol. 150, pp. 1--29, 2014.

\bibitem{ref:viragLimits2}
A.~Bloemendal and B.~Vir{\'a}g, ``{Limits of spiked random matrices II},''
  \emph{Ann. Probab.}, vol.~44, no.~4, pp. 2726--2769, 2016.

\bibitem{ref:viragLimits1}
------, ``Limits of spiked random matrices {I},'' \emph{Probab. Theory Relat.
  Fields}, vol. 156, no. 3-4, pp. 795--825, 2013.

\bibitem{ref:peche}
S.~P{\'e}ch{\'e}, ``The largest eigenvalue of small rank perturbations of
  {Hermitian} random matrices,'' \emph{Probab. Theory Relat. Fields}, vol. 134,
  no.~1, pp. 127--173, 2006.

\bibitem{ref:nadakuditi}
R.~R. Nadakuditi and J.~W. Silverstein, ``Fundamental limit of sample
  generalized eigenvalue based detection of signals in noise using relatively
  few signal-bearing and noise-only samples,'' \emph{IEEE J. Sel. Top. Signal
  Process.}, vol.~4, no.~3, pp. 468--480, 2010.

\bibitem{ref:dWang}
Z.~Bao and D.~Wang, ``Eigenvector distribution in the critical regime of {BBP}
  transition,'' \emph{Probab. Theory Relat. Fields}, pp. 1--81, 2021.

\bibitem{ref:hoyle}
D.~C. Hoyle and M.~Rattray, ``Statistical mechanics of learning multiple
  orthogonal signals: Asymptotic theory and fluctuation effects,'' \emph{Phys.
  Rev. E}, vol.~75, p. 016101, 2007.

\bibitem{ref:georgeFob}
F.~Benaych-Georges and R.~R. Nadakuditi, ``The singular values and vectors of
  low rank perturbations of large rectangular random matrices,'' \emph{J.
  Multivar. Anal.}, vol. 111, pp. 120--135, 2012.

\bibitem{ref:capitaineFo}
M.~Capitaine, ``Limiting eigenvectors of outliers for spiked
  information-plus-noise type matrices,'' in \emph{S{\'e}minaire de
  Probabilit{\'e}s XLIX}, ser. Lecture Notes in Mathematics.\hskip 1em plus
  0.5em minus 0.4em\relax Springer International Publishing, 2018, pp.
  119--164.

\bibitem{ref:dingFo}
X.~Ding, ``High dimensional deformed rectangular matrices with applications in
  matrix denoising,'' \emph{Bernoulli}, vol.~26, no.~1, pp. 387--417, 2020.

\bibitem{ref:bourChi}
P.~Bourgade and H.-T. Yau, ``The eigenvector moment flow and local quantum
  unique ergodicity,'' \emph{Commun. Math. Phys.}, vol. 350, no.~1, pp.
  231--278, 2017.

\bibitem{ref:knowChi}
A.~Knowles and J.~Yin, ``Eigenvector distribution of {W}igner matrices,''
  \emph{Probab. Theory Relat. Fields}, vol. 155, no.~3, pp. 543--582, 2013.

\bibitem{ref:taoChi}
T.~Tao and V.~Vu, ``Random matrices: {U}niversal properties of eigenvectors,''
  \emph{Random Matrices: Theory Appl.}, vol.~01, no.~01, p. 1150001, 2012.

\bibitem{ref:beniChi}
L.~Benigni, ``Fermionic eigenvector moment flow,'' \emph{Probab. Theory Relat.
  Fields}, vol. 179, no.~3, pp. 733--775, 2021.

\bibitem{ref:marChi}
J.~Marcinek and H.-T. Yau, ``High dimensional normality of noisy
  eigenvectors,'' \emph{arXiv:2005.08425 [math.PR]}, 2020.

\bibitem{ref:capitaine}
M.~Capitaine and C.~Donati-Martin, ``Non universality of fluctuations of
  outlier eigenvectors for block diagonal deformations of wigner matrices,''
  \emph{{ALEA} Lat. Am. J. Probab. Math. Stat.}, vol.~18, no.~1, pp. 129--165,
  2021.

\bibitem{ref:baoChi}
Z.~Bao, X.~Ding, J.~Wang, and K.~Wang, ``Principal components of spiked
  covariance matrices in the supercritical regime,'' \emph{arXiv:1907.12251
  [math.ST]}, 2019.

\bibitem{ref:terry}
P.~B. Denton, S.~J. Parke, T.~Tao, and X.~Zhang, ``Eigenvectors from
  eigenvalues: {A} survey of a basic identity in linear algebra,'' \emph{Bull.
  Amer. Math. Soc.}, vol.~59, no.~1, pp. 31--58, 2021.

\bibitem{ref:feldman}
M.~J. Feldman, ``Spiked singular values and vectors under extreme aspect
  ratios,'' \emph{arXiv:2104.15127 [math.ST]}, 2021.

\bibitem{ref:johnstoneRoys}
I.~M. Johnstone and B.~Nadler, ``{Roy’s largest root test under rank-one
  alternatives},'' \emph{Biometrika}, vol. 104, no.~1, pp. 181--193, 2017.

\bibitem{ref:prathapRoyRoot}
P.~Dharmawansa, B.~Nadler, and O.~Shwartz, ``Roy’s largest root under
  rank-one perturbations: The complex valued case and applications,'' \emph{J.
  Multivariate Anal.}, vol. 174, p. 104524, 2019.

\bibitem{ref:wang}
D.~Wang, ``The largest eigenvalue of real symmetric, {Hermitian} and
  {Hermitian} self-dual random matrix models with rank one external source,
  part {I},'' \emph{J. Statist. Phys.}, vol. 146, no.~4, pp. 719--761, 2012.

\bibitem{ref:ghur1}
T.~Wirtz and T.~Guhr, ``Distribution of the smallest eigenvalue in the
  correlated {W}ishart model,'' \emph{Phys. Rev. Lett.}, vol. 111, p. 094101,
  2013.

\bibitem{ref:ghur2}
------, ``Distribution of the smallest eigenvalue in complex and real
  correlated {Wishart} ensembles,'' \emph{J. Phys. A: Math. Theor.}, vol.~47,
  no.~7, p. 075004, 2014.

\bibitem{ref:akemann}
G.~Akemann and P.~Vivo, ``Compact smallest eigenvalue expressions in
  {W}ishart{\textendash}{L}aguerre ensembles with or without a fixed trace,''
  \emph{J. Stat. Mech. Theory Exp.}, vol. 2011, no.~05, p. P05020, 2011.

\bibitem{ref:AOSing}
A.~Onatski, ``The {T}racy–{W}idom limit for the largest eigenvalues of
  singular complex {W}ishart matrices,'' \emph{Ann. Appl. Probab.}, vol.~18,
  no.~2, pp. 470--490, 2008.

\bibitem{ref:srivastava}
M.~Srivastava, ``Singular {W}ishart and multivariate beta distributions,''
  \emph{Ann. Stat.}, vol.~31, no.~5, pp. 1537--1560, 2003.

\bibitem{ref:ratnarajahSing}
T.~Ratnarajah and R.~Vaillancourt, ``Complex singular wishart matrices and
  applications,'' \emph{Comput. Math. with Appl.}, vol.~50, no.~3, pp.
  399--411, 2005.

\bibitem{ref:ratnarajahMulti}
T.~Ratnarajah, ``Complex singular wishart matrices and multiple-atenna
  systems,'' in \emph{13th IEEE SSP}, 2005, pp. 1032--1037.

\bibitem{ref:mallik}
R.~Mallik, ``The pseudo-{W}ishart distribution and its application to {MIMO}
  systems,'' \emph{IEEE Trans. Inf. Theory}, vol.~49, no.~10, pp. 2761--2769,
  2003.

\bibitem{ref:shafi}
P.~Smith, S.~Roy, and M.~Shafi, ``Capacity of {MIMO} systems with
  semicorrelated flat fading,'' \emph{IEEE Trans. Inf. Theory}, vol.~49,
  no.~10, pp. 2781--2788, 2003.

\bibitem{ref:garcia}
J.~A. Dı\`az-Garcı\`a, R.~G. J\`aimez, and K.~V. Mardia, ``{W}ishart and
  {P}seudo-{W}ishart distributions and some applications to shape theory,''
  \emph{J. Multivar. Anal.}, vol.~63, no.~1, pp. 73--87, 1997.

\bibitem{ref:goodall}
C.~R. Goodall and K.~V. Mardia, ``Multivariate aspects of shape theory,''
  \emph{Ann. Stat.}, vol.~21, no.~2, pp. 848--866, 1993.

\bibitem{ref:uhlig}
H.~Uhlig, ``On singular {W}ishart and singular multivariate beta
  distributions,'' \emph{Ann. Stat.}, vol.~22, no.~1, pp. 395--405, 1994.

\bibitem{ref:james}
A.~T. James, ``Distributions of matrix variates and latent roots derived from
  normal samples,'' \emph{Ann. Math. Statist.}, vol.~35, no.~2, pp. 475--501,
  1964.

\bibitem{ref:takemura}
A.~Takemura, \emph{Zonal Polynomials}, ser. Lecture Notes-Monograph.\hskip 1em
  plus 0.5em minus 0.4em\relax IMS, 1984, vol.~4.

\bibitem{ref:harish}
Harish-Chandra, ``Differential operators on a semisimple {L}ie algebra,''
  \emph{Am. J. Math.}, vol.~79, no.~1, pp. 87--120, 1957.

\bibitem{ref:itzykson}
C.~Itzykson and J.~Zuber, ``The planar approximation {II},'' \emph{J. Math.
  Phys.}, vol.~21, no.~3, pp. 411--421, 1980.

\bibitem{ref:peterJ}
P.~J. Forrester, ``Probability densities and distributions for spiked and
  general variance {W}ishart $\beta$-ensembles,'' \emph{Random Matrices: Theory
  Appl.}, vol.~02, no.~04, p. 1350011, 2013.

\bibitem{ref:peter22}
------, ``Rank $1$ perturbations in random matrix theory -- a review of exact
  results,'' \emph{arXiv:2201.00324 [math-ph]}, 2022.

\bibitem{ref:szego}
G.~Szeg{\"o}, \emph{Orthogonal Polynomials}, ser. American Math. Soc:
  Colloquium publ.\hskip 1em plus 0.5em minus 0.4em\relax American Mathematical
  Society, 1975.

\bibitem{ref:gradshteyn}
I.~Gradshteyn and I.~Ryzhik, \emph{Table of Integrals, Series, and Products},
  7th~ed.\hskip 1em plus 0.5em minus 0.4em\relax Boston: Academic Press, 2007.

\bibitem{ref:erdelyi}
A.~Erdélyi, \emph{Higher Transcendental Functions}.\hskip 1em plus 0.5em minus
  0.4em\relax New York: McGraw-Hill, 1953, vol.~1.

\bibitem{ref:prathapSIAM}
P.~Dharmawansa, M.~R. McKay, and Y.~Chen, ``Distributions of {D}emmel and
  related condition numbers,'' \emph{SIAM J. Matrix Anal. Appl.}, vol.~34,
  no.~1, pp. 257--279, 2013.

\bibitem{ref:vaart}
A.~W. Van~der Vaart, \emph{{Asymptotic Statistics}}, ser. Cambridge Series in
  Statistical and Probabilistic Mathematics.\hskip 1em plus 0.5em minus
  0.4em\relax Cambridge: Cambridge University Press, 2000.

\bibitem{ref:onatskiSphe}
A.~Onatski, M.~J. Moreira, and M.~Hallin, ``Asymptotic power of sphericity
  tests for high-dimensional data,'' \emph{Ann. Stat.}, vol.~41, no.~3, pp.
  1204--1231, 2013.

\bibitem{ref:onatskiSpiked}
I.~M. Johnstone and A.~Onatski, ``Testing in high-dimensional spiked models,''
  \emph{Ann. Stat.}, vol.~48, no.~3, pp. 1231--1254, 2020.

\bibitem{ref:dobriban}
E.~Dobriban, ``Sharp detection in {PCA} under correlations: {A}ll eigenvalues
  matter,'' \emph{Ann. Stat.}, vol.~45, no.~4, pp. 1810--1833, 2017.

\bibitem{ref:chiani}
M.~Chiani, M.~Win, and A.~Zanella, ``On the capacity of spatially correlated
  {MIMO} {R}ayleigh-fading channels,'' \emph{IEEE Trans. Inf. Theory}, vol.~49,
  no.~10, pp. 2363--2371, 2003.

\bibitem{ref:dharma}
P.~Dharmawansa and I.~M. Johnstone, ``Joint density of eigenvalues in spiked
  multivariate models,'' \emph{Stat}, vol.~3, no.~1, pp. 240--249, 2014.

\bibitem{ref:batemanTransforms}
H.~Bateman, \emph{Tables of Integral Transforms [Volumes I \& II]}.\hskip 1em
  plus 0.5em minus 0.4em\relax McGraw-Hill Book Company, 1954, vol.~1.

\bibitem{ref:erdelyiGer}
A.~Erd{\'e}lyi, ``Beitrag zur theorie der konfluenten hypergeometrischen
  funktionen von mehreren ver{\"a}nderlichen,''
  \emph{H{\"o}lder-Pichler-Tempsky in Komm.}, 1937.

\bibitem{ref:milton}
M.~Abramowitz and I.~A. Stegun, \emph{Handbook of Mathematical Functions With
  Formulas, Graphs, and Mathematical Tables}.\hskip 1em plus 0.5em minus
  0.4em\relax US Government Printing Office, 1964, vol.~55.

\bibitem{ref:khatri}
C.~G. Khatri, ``On the moments of traces of two matrices in three situations
  for complex multivariate normal populations,'' \emph{Sankhyā}, vol.~32,
  no.~1, pp. 65--80, 1970.

\bibitem{ref:prathapJMVA}
P.~Dharmawansa, ``Some new results on the eigenvalues of complex non-central
  {W}ishart matrices with a rank-1 mean,'' \emph{J. Multivariate Anal.}, vol.
  149, pp. 30--53, 2016.

\bibitem{ref:SCNArxiv}
P.~Dissanayake, P.~Dharmawansa, and Y.~Chen, ``Distribution of the scaled
  condition number of single-spiked complex {Wishart} matrices,''
  \emph{arXiv:2105.05307 [math.PR]}, 2021.

\end{thebibliography}
\end{document}